\documentclass{article}

\usepackage{ifthen}

\newboolean{isMAA}
\setboolean{isMAA}{false} 

\ifthenelse{\boolean{isMAA}}
{
\usepackage{maa-monthly}
}
{
\usepackage[a4paper]{geometry}
\usepackage{amssymb} 
\usepackage{amsmath}
\usepackage{amsthm}
}

\newboolean{isCoda}
\setboolean{isCoda}{false}

\usepackage[utf8]{inputenc}
\usepackage[english]{babel}
\usepackage{graphicx}
\usepackage{xcolor}
\usepackage{xspace}
\usepackage{paralist}
\usepackage{mathrsfs}
\usepackage{enumitem}
\usepackage{tikz}
\usepackage{multicol}
\usepackage[
colorlinks=true,
urlcolor=blue,
linkcolor=black
]{hyperref}

\usepackage[T1]{fontenc}
\usepackage[font=small,labelfont=bf]{caption}

\DeclareCaptionLabelFormat{andtable}{#1~#2  \&  \tablename~\thetable}

\usepackage{picinpar,moresize,xfrac,graphpap,dcolumn,wrapfig,graphicx}
\usepackage{graphpap}

\usepackage{tikz}

\newcommand{\charliex}[1]{#1}

\newcommand{\russell}[1]{{\color[rgb]{0.3,0,0.8}\textbf{russell: }{#1\xspace}\color{black}}}

\newcommand{\removex}[1]{}

\newcommand{\myboldhead}[1]{\vspace{0in}\hspace{-.22in}\textbf{#1.}}

\newcommand{\RP}[1]{\mathbb{R}P^{#1}\xspace}

\newcommand{\face}{\ensuremath{\mathcal{F}}\xspace}

\newcommand{\mtx}[1]{\mathbf{{#1}}} 
\newcommand{\ptq}[1]{\mtx{#1}}
\newcommand{\plq}[1]{\mtx{#1}^*}
\newcommand{\ptc}[1]{\mtx{#1}}   

\newcommand{\VC}{\ptq{V}}

\newcommand{\axis}{\Lambda}     
\newcommand{\spear}{\lambda}     

\theoremstyle{definition}
\newtheorem{example}{Example}[section]
\newtheorem{theorem}{Theorem}

\newtheorem{lemma}[theorem]{Lemma}

\theoremstyle{definition}
\newtheorem{definition}[theorem]{Definition}
\newtheorem{remark}[theorem]{Remark}

\definecolor{b1}{rgb}{0.158099,0.313781,0.636957} 
\definecolor{b2}{rgb}{0.525367,0.691857,0.998936} 
\definecolor{g1}{rgb}{0.256000,0.640000,0.576000} 
\definecolor{g2}{rgb}{0.559573,0.800781,0.760580} 
\definecolor{p1}{rgb}{0.416000,0.192000,0.640000} 
\definecolor{p2}{rgb}{0.680880,0.561880,0.799881} 
\definecolor{r1}{rgb}{0.800000,0.200000,0.000000} 
\definecolor{r2}{rgb}{1.000000,0.702595,0.603461} 
\definecolor{k1}{gray}{0}
\definecolor{k2}{gray}{0.7} 

\title{Penrose's eight-conic theorem}
\author{Russell Arnold\footnote{Mathematical Section, Goetheanum, Dornach, Switzerland} \and
Albert Chern\footnote{Department of Computer Science and Engineering, University of California San Diego, La Jolla, USA} \and Morten Eide\footnote{Independent, Jevnaker, Norway} \and Charles Gunn\footnote{Independent, Falkensee, Germany (Corresponding author: projgeom@gmail.com) } \and Thomas Neukirchner\footnote{Free Waldorf School, Karlsruhe, Germany} \and Roger Penrose\footnote{Mathematical Institute, University of Oxford, UK}}
\date{\today}

\begin{document}


\maketitle
\charliex{
\begin{abstract}
    This article proves the following theorem, first enunciated by Roger Penrose about 70 years ago but never published: In $\RP{2}$, if conics are assigned to seven of the vertices of a combinatorial cube such that (i) conics connected by an edge are in double contact, and (ii) the chords of contact associated to a cube face meet in a common point, then there exists an eighth conic such that the completed cube satisfies (i) and (ii). The theorem turns out to be a remarkable generalization of many well-known theorems of projective geometry --  Pappus, Desargues, Pascal, Brianchon and Monge are the best-known ones. This archetypal principle provides a unifying framework in which the myriad specializations of the theorem and their interrelationships can be grasped as an organic whole, enriching the field of projective geometry and opening new vistas for research. The article begins with a series of motivational examples. It then gives a geometric proof assuming that the conics are regular, followed by an algebraic one that removes this restriction. The geometric proof is obtained as a slice of an analogous theorem for quadrics in $\RP{3}$; the algebraic one is based on the determinants of a special matrix associated to the configuration of conics.    
\end{abstract}
}
\removex{
\begin{abstract}
    This article proves the following theorem, first enunciated by Roger Penrose about 70 years ago: In $\RP{2}$, if regular conics are assigned to seven of the vertices of a combinatorial cube such that (i) conics connected by an edge are in double contact, and (ii) the chords of contact associated to a cube face meet in a common point, then there exists an eighth conic such that the completed cube satisfies (i) and (ii). This conic is unique if not all the common points of condition (ii) are the same. 
    The proof is based on the following analogous theorem, which is also proved: In $\RP{3}$, if regular quadrics are assigned to seven of the vertices of a combinatorial cube such that (i) quadrics connected by an edge are in ring contact, and (ii) the ring planes associated to a cube face meet in a common axis, then there exists an eighth quadric such that the completed cube satisfies (i) and (ii). This quadric is unique if not all the common axes of condition (ii) are the same.
\end{abstract}
}
\section{Introduction}
\begin{figure}[h]
    \centering
    \setlength{\unitlength}{1pt}
    \begin{picture}(400,110)
        \put(0,0){
        \includegraphics[width=200pt]{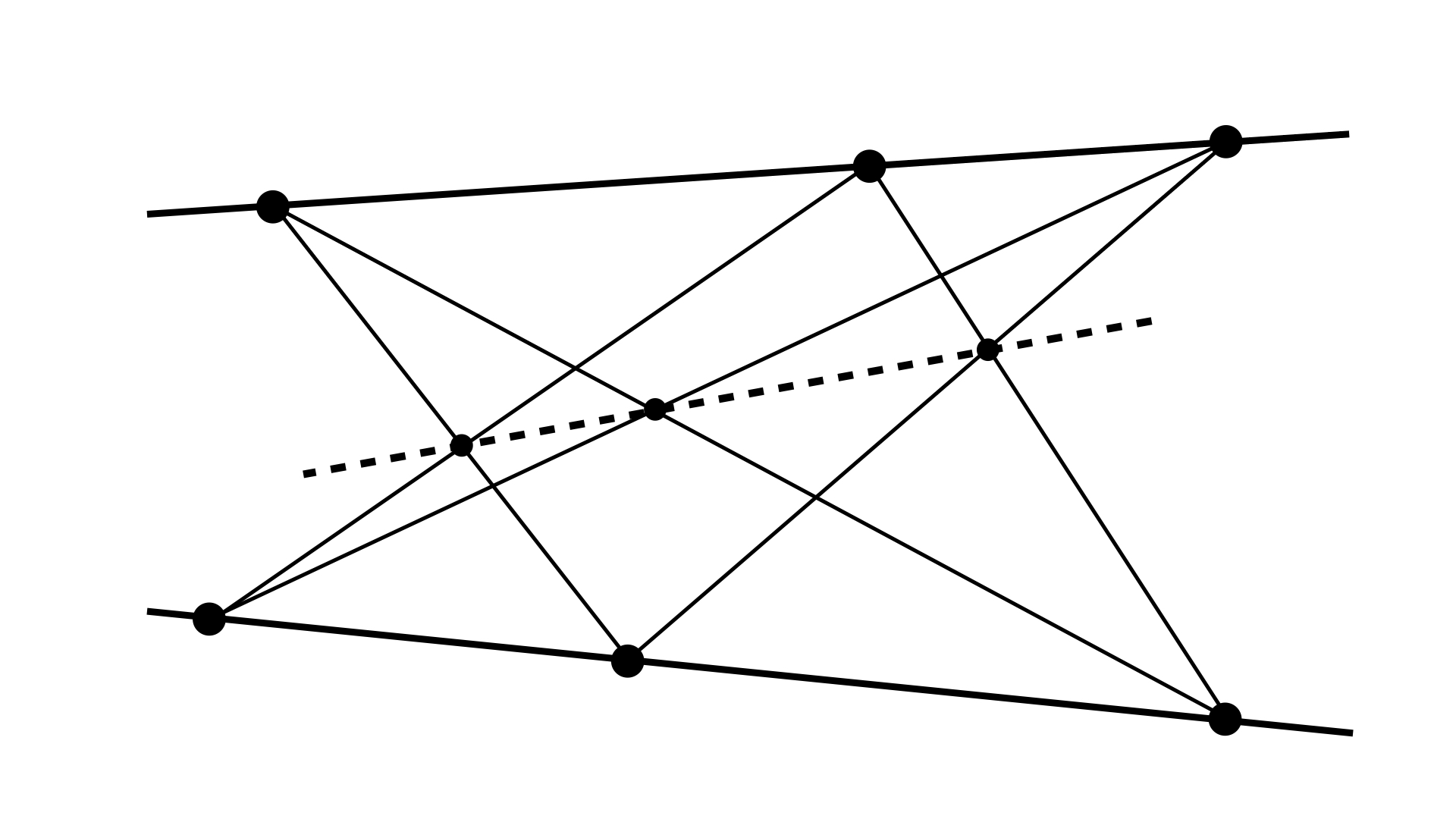}
        }
        \put(200,0){
        \includegraphics[width=200pt]{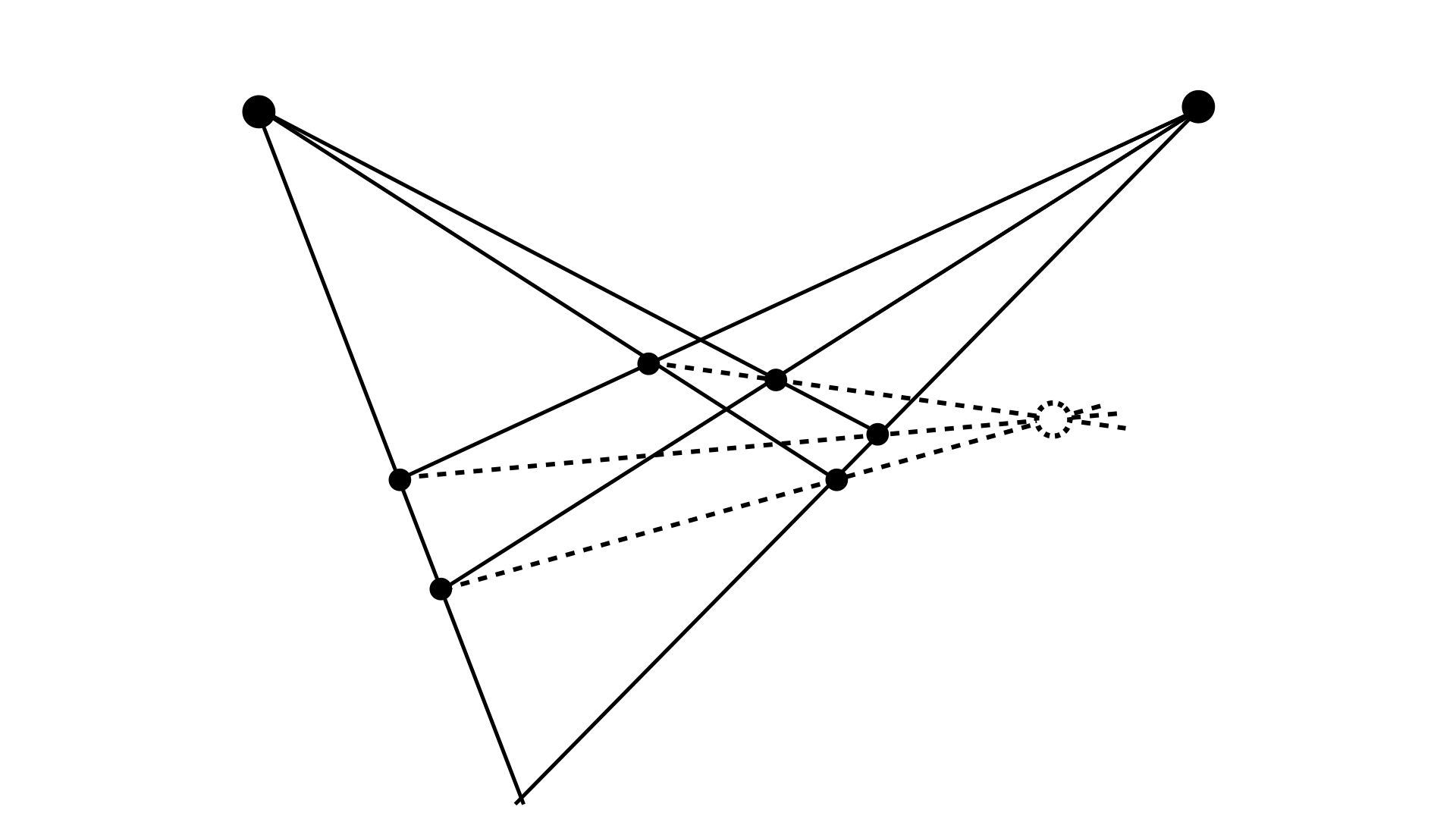}
        }
        \put(12,80){\small\(\ell\)}
        \put(12,27){\small\(\ell'\)}
        \put(35,90){\small\(A\)}
        \put(120,95){\small\(B\)}
        \put(170,97){\small\(C\)}
        \put(27,15){\small\(A'\)}
        \put(87,10){\small\(B'\)}
        \put(170,2){\small\(C'\)}
        \put(227,95){\small\(L\)}
        \put(373,95){\small\(L'\)}
        \put(240,70){\small\(a\)}
        \put(260,73){\small\(b\)}
        \put(270,83){\small\(c\)}
        \put(320,80){\small\(c'\)}
        \put(330,67){\small\(b'\)}
        \put(350,72){\small\(a'\)}
    \end{picture}
    \caption{Pappus's theorem (340 AD) (left) and its dual theorem (right).}
    \label{fig:PappusInIntro}
\end{figure}

Let us start with a classical theorem from 340 AD by Pappus of Alexandria (290--350 AD) \cite{jones2013pappus} and its dual counterpart (Figure~\ref{fig:PappusInIntro}). 

    \begin{multicols}{2}
    \emph{
    Consider a pair of lines in the plane (call them \(\ell\) and \(\ell'\)). Assign three arbitrary points on each line (\(A,B,C\) on \(\ell\) and \(A',B',C'\) on \(\ell'\)).  Join these three pairs of points in a crossing pattern, giving three meeting points (namely, the meet of lines \(AB'\) and \(A'B\), the meet of lines \(BC'\) and \(B'C\), and the meet of lines \(CA'\) and \(C'A\)).  Then these three points lie on a line.}
    
    \columnbreak
    \emph{
    Consider a pair of points in the plane (call them \(L\) and \(L'\)).  Draw three lines through each point (\(a,b,c\) through \(L\) and \(a',b',c'\) through \(L'\)).  Meet these three pairs of lines in a diagonal pattern, giving three lines (namely the join of points \(ab'\) and \(a'b\), the join of points \(bc'\) and \(b'c\), and the join of points \(ca'\) and \(c'a\)).  Then these three lines pass through a point.}
\end{multicols}

A fascinating quality of Pappus's theorem is that the statement of the theorem requires very little geometric structure, simply that two distinct points can be joined with a unique line, and two distinct lines meet in a unique point (possibly at the line at infinity), reminiscent of the practice of perspective drawing.  No lengths, angles, or parallelisms are involved.  This is a theorem of \emph{projective geometry}.
Note that one may swap ``points'' and ''lines'' and swap ``meet'' and ``join'' in a planar projective geometric theorem to obtain the dual theorem, as shown above.

It is fair to say that neither Pappus's theorem nor its dual are obvious.  A natural question is whether there is a general pattern to theorems like Pappus's that could shed light on Pappus's theorem, making it more intuitive and easier to discover.

\begin{figure}
    \centering
    \setlength{\unitlength}{1pt}
    \begin{picture}(400,110)
        \put(0,0){
        \includegraphics[width=200pt]{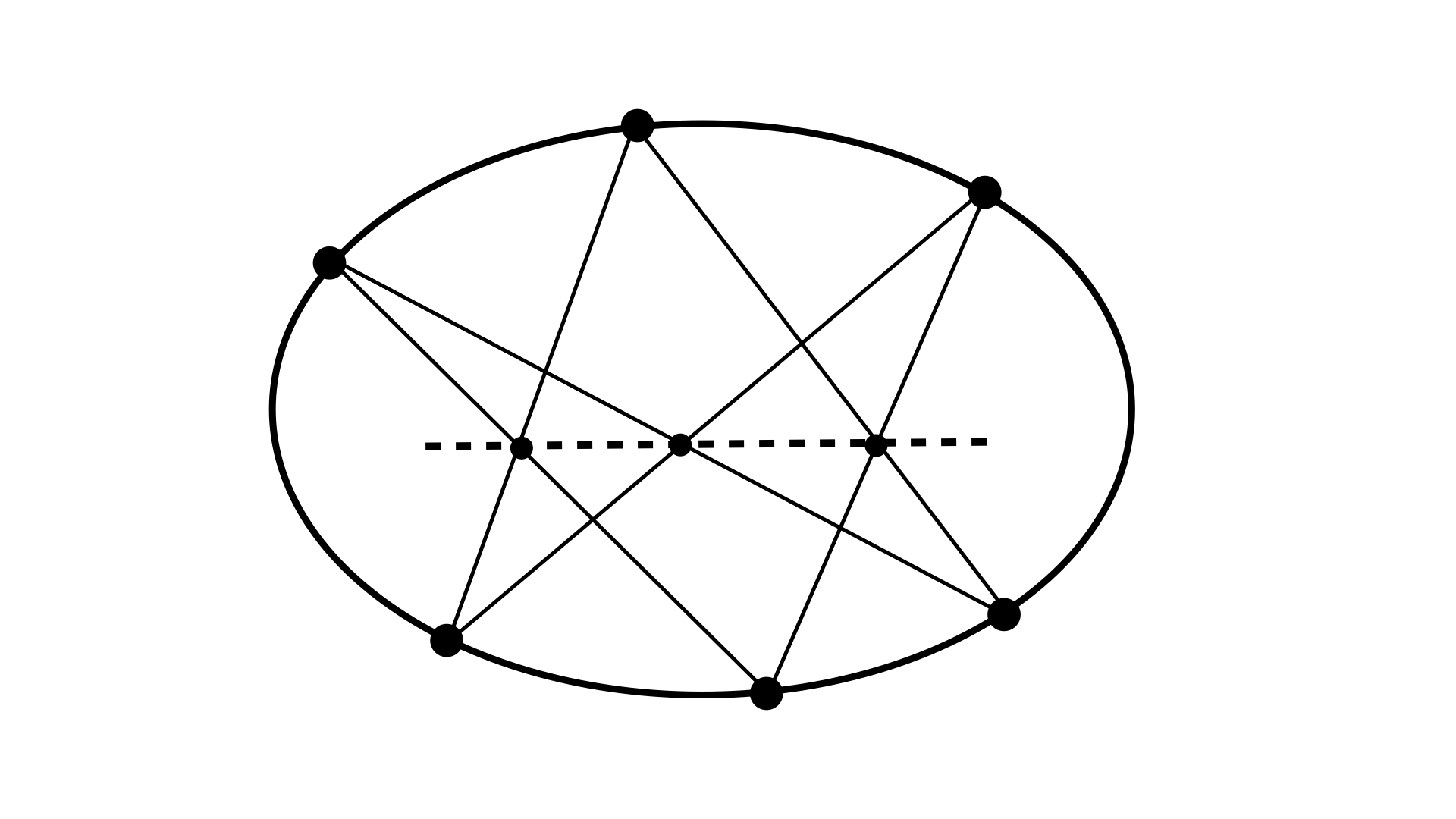}
        }
        \put(200,0){
        \includegraphics[width=200pt]{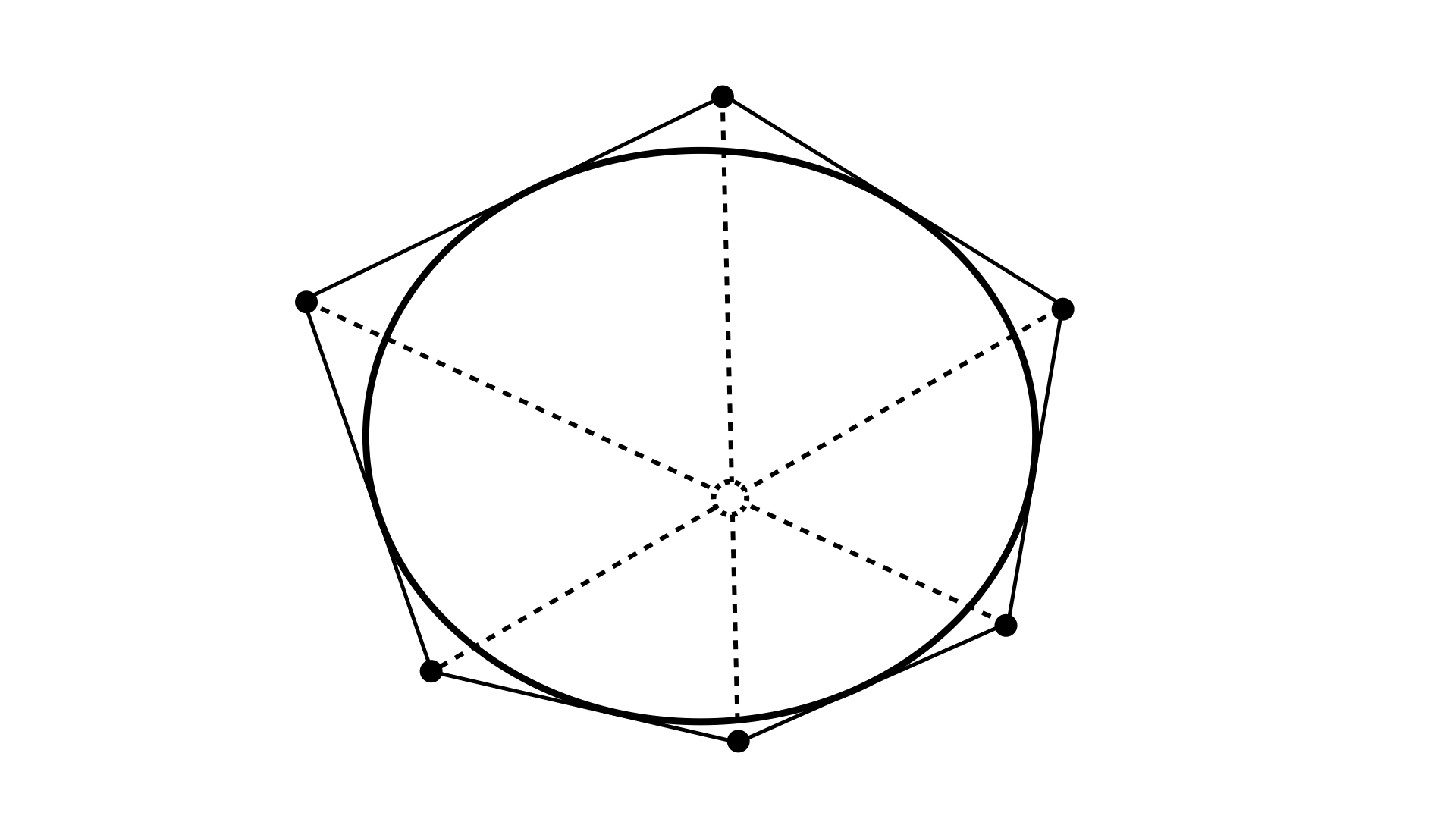}
        }
        \put(35,75){\small\(A\)}
        \put(90,100){\small\(B\)}
        \put(140,90){\small\(C\)}
        \put(50,15){\small\(A'\)}
        \put(105,5){\small\(B'\)}
        \put(145,17){\small\(C'\)}
        \put(270,90){\small\(a\)}
        \put(325,87){\small\(b'\)}
        \put(350,45){\small\(c\)}
        \put(322,8){\small\(a'\)}
        \put(279,5){\small\(b\)}
        \put(242,40){\small\(c'\)}
    \end{picture}
    \caption{Pascal's theorem (1640) (left) and its dual theorem---Brianchon's theorem (1807) (right).}
    \label{fig:PascalInIntro}
\end{figure}

Blaise Pascal (1623--1662), at age 16, discovered a remarkable generalization of Pappus's theorem, known as Pascal's mystic hexagram \cite{del2020pascal}. The initial line pair of Pappus's theorem, as well as the initial point pair of the dual Pappus theorem, is a degeneration of a conic section.  After replacing the pair of lines with an arbitrary conic, the rest of the statement of the theorem remains true (Figure~\ref{fig:PascalInIntro}).

\begin{quote}
    \begin{multicols}{2}
    \emph{
    Consider a conic in the plane. Choose six points on the conic (labeled \(A,B,C,A',B',C'\)).  Join these three pairs of points in a crossing pattern, giving three meeting points (namely, the meet of lines \(AB'\) and \(A'B\), the meet of lines \(BC'\) and \(B'C\), and the meet of lines \(CA'\) and \(C'A\)).  Then these three points lie on a line.}
    
    \columnbreak
    \emph{
    Consider a conic in the plane.  Draw six lines tangent to the conic (labeled \(a,b,c,a',b',c'\)).  Meet these three pairs of lines in a diagonal manner, giving three joining lines (namely the join of points \(ab'\) and \(a'b\), the join of points \(bc'\) and \(b'c\), and the join of points \(ca'\) and \(c'a\)).  Then these three lines go through a point.}
\end{multicols}
\end{quote}
The dual theorem (right-hand side) is known as Brianchon's Theorem, named after Charles Julien Brianchon (1783--1864).

\begin{figure}
    \centering
    \setlength{\unitlength}{1pt}
    \begin{picture}(400,110)
        \put(0,0){
        \includegraphics[width=200pt]{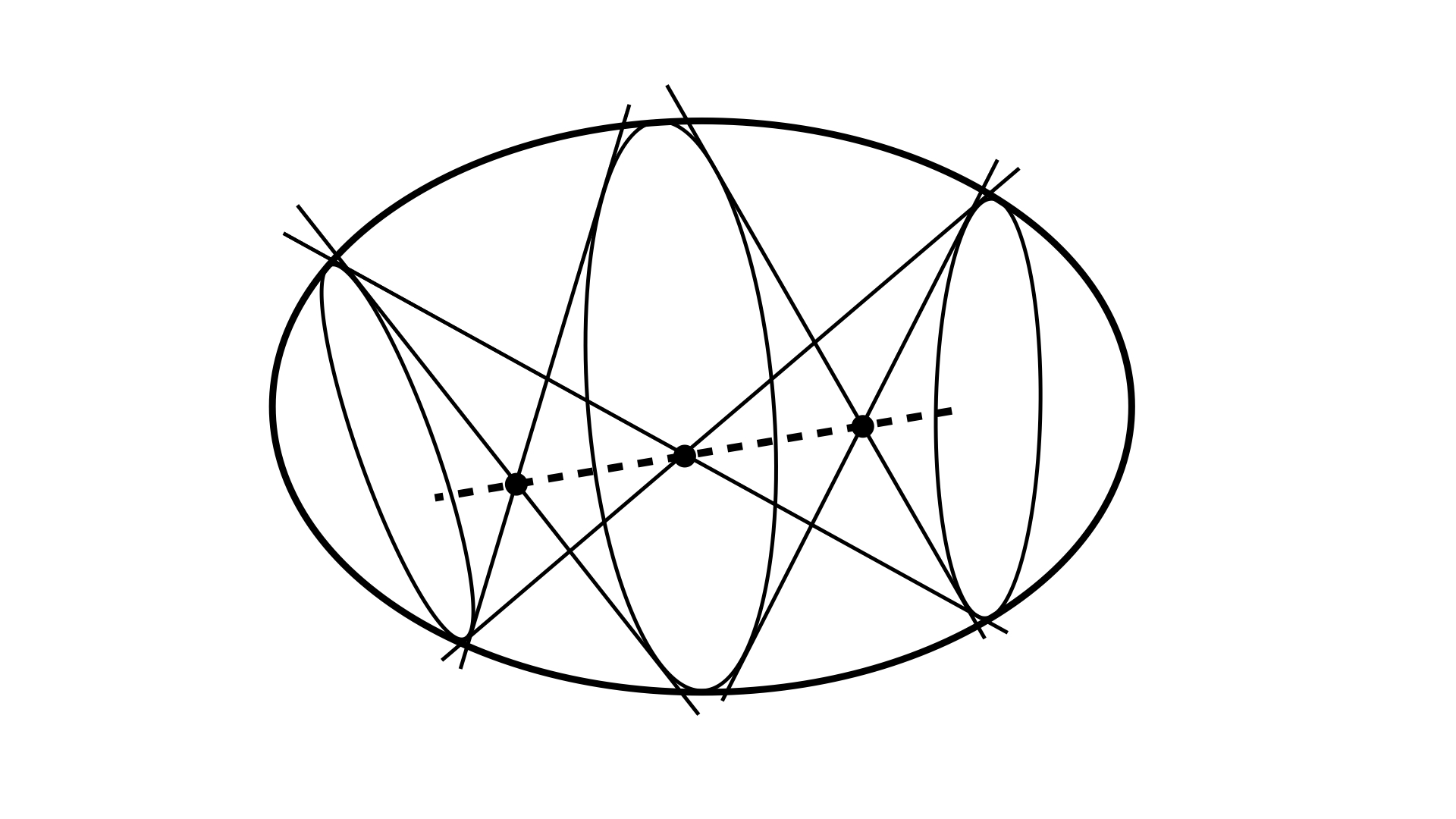}
        }
        \put(200,0){
        \includegraphics[width=200pt]{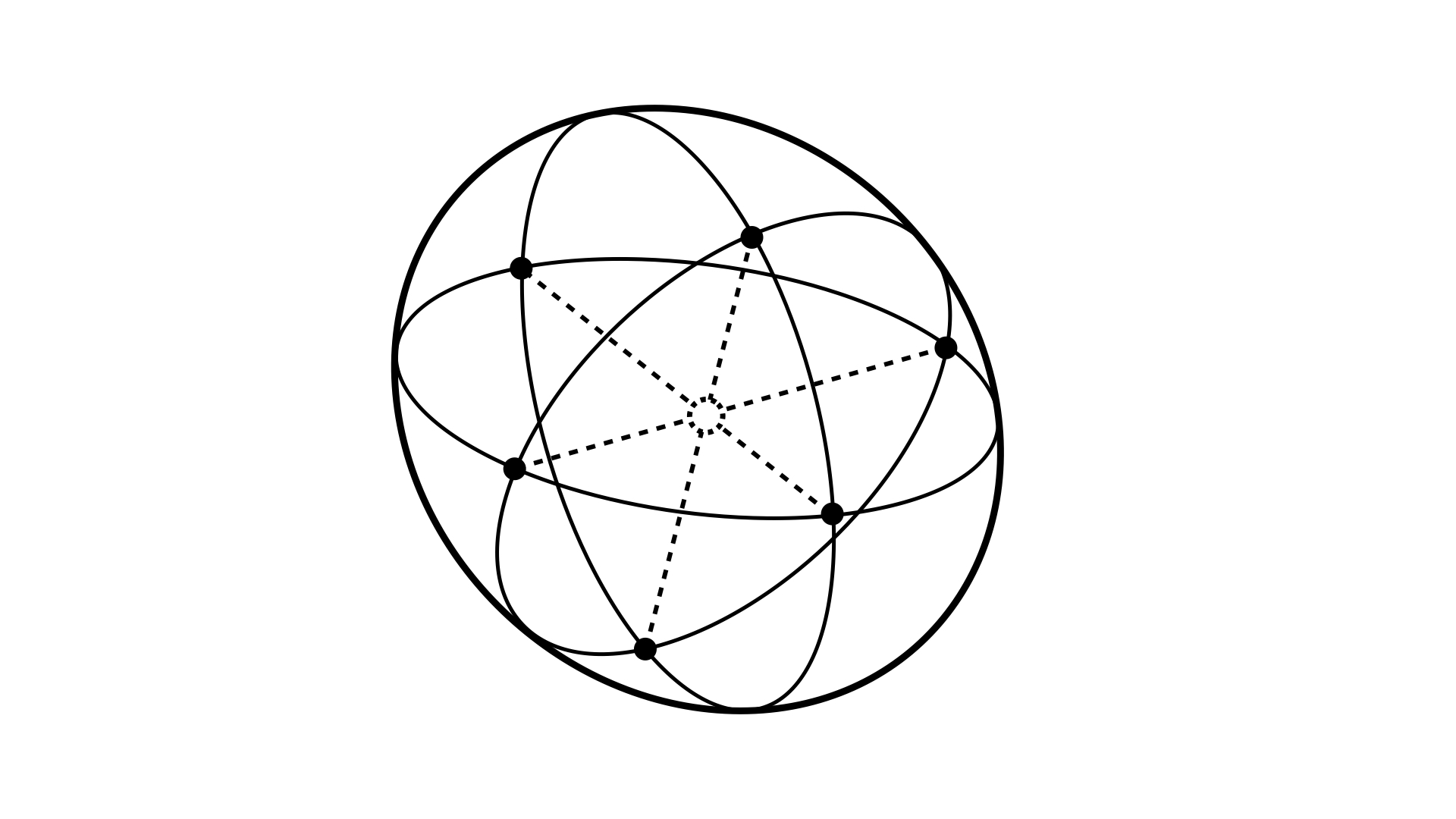}
        }
        \put(30,50){\small\(\ptq{S}_0\)}
        \put(50,50){\small\(\ptq{T}^1\)}
        \put(92,60){\small\(\ptq{T}^2\)}
        \put(135,60){\small\(\ptq{T}^3\)}
        \put(250,80){\small\(\ptq{S}_0\)}
        \put(262,60){\small\(\ptq{T}^1\)}
        \put(285,85){\small\(\ptq{T}^2\)}
        \put(320,72){\small\(\ptq{T}^3\)}
    \end{picture}
    \caption{Salmon's theorem (1848) and its dual theorem.}
    \label{fig:SalmonInIntro}
\end{figure}

In 1848, George Salmon published a textbook on conic sections, which includes a further generalization of the Pascal--Brianchon theorem \cite[\S 264]{salmon1917}, rediscovered by Evelyn, Money-Coutts, and Tyrrell in 1974 \cite[\S 2.3]{evelyn1974seven}.  The three point pairs \(AA'\), \(BB'\), \(CC'\) of Pascal's theorem are degenerations of thin ellipses that are in double contact with the initial conic.  And the lines \(AB'\) and \(A'B\) are the limit of a pair of common tangents between the two thin ellipses collapsing to \(AA'\) and \(BB'\).  After replacing the three point pairs with three conics in double contact with a fourth, the remainder of Pappus's theorem still holds (Figure~\ref{fig:SalmonInIntro}).

\begin{quote}
    \begin{multicols}{2}
    \emph{
    If three conics \(\ptq{T}^1,\ptq{T}^2,\ptq{T}^3\) are each in double contact with a fourth \(\ptq{S}_0\), then the three intersection points of pairs of common tangents (namely, between \(\ptq{T}^1\ptq{T}^2\), between \(\ptq{T}^2\ptq{T}^3\), and between \(\ptq{T}^3\ptq{T}^1\)), lie on a line.}
    
    \columnbreak
    \emph{
    If three conics \(\ptq{T}^1,\ptq{T}^2,\ptq{T}^3\) are each in double contact with a fourth \(\ptq{S}_0\), then a common chord between \(\ptq{T}^1\ptq{T}^2\), a common chord between \(\ptq{T}^2\ptq{T}^3\), and a common chord between \(\ptq{T}^3\ptq{T}^1\), all pass through a point.}\footnote{The original statement is ``If three conics are each in double contact with a fourth, then six of their chords of intersection form a complete quadrilateral.''}
\end{multicols}
\end{quote}

The reader may have observed a pattern of this progression of generalizing Pappus's theorem.  Line pairs and point pairs are replaced by conic sections.  The incidence relations of point pairs on a conic or line pairs tangent to a conic are both generalized by two conics in double contact.

\begin{figure}
    \centering
    \setlength{\unitlength}{1pt}
    \begin{picture}(400,110)
        \put(0,0){
        \includegraphics[width=200pt]{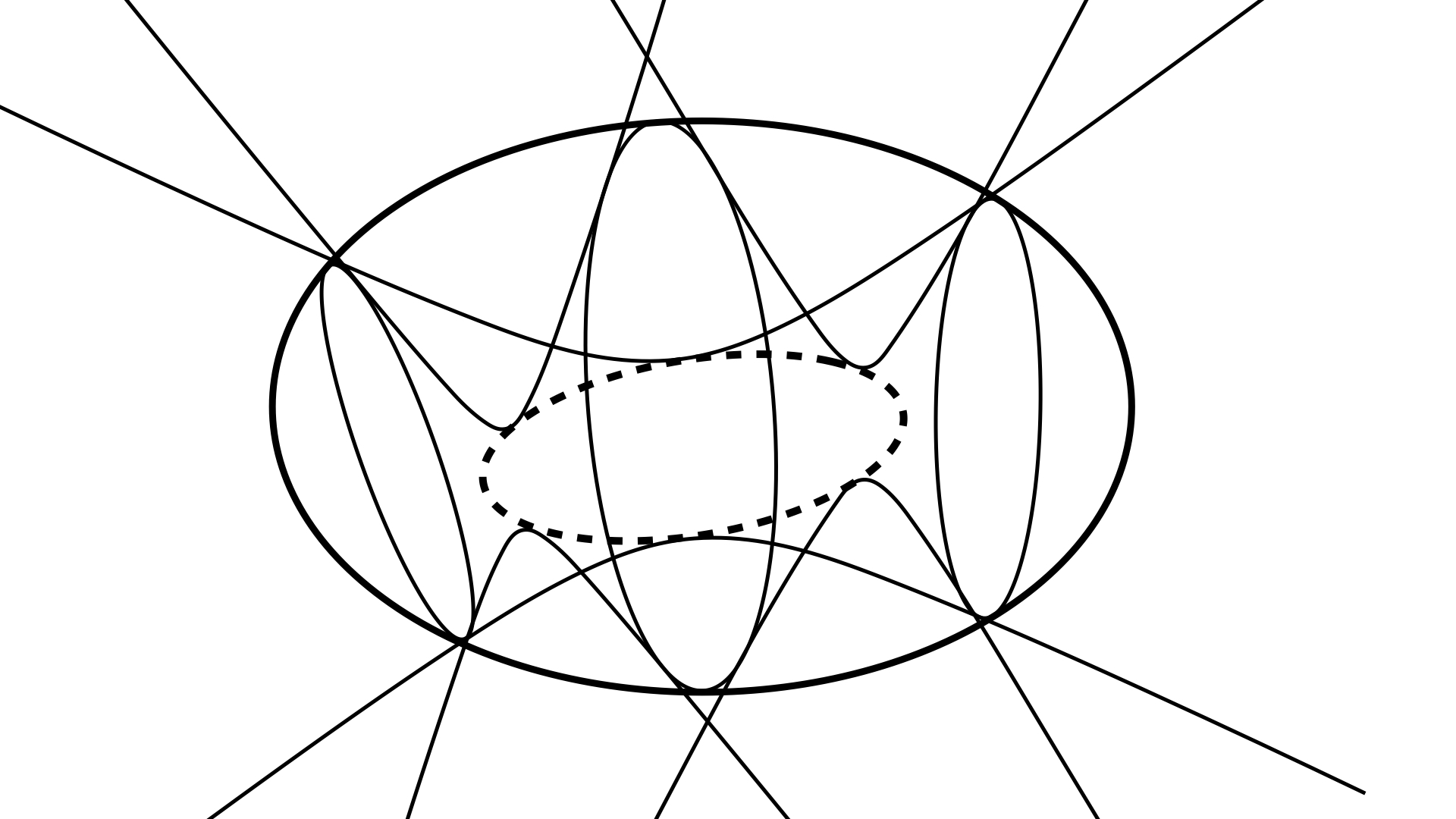}
        }
        \put(200,0){
        \includegraphics[width=200pt]{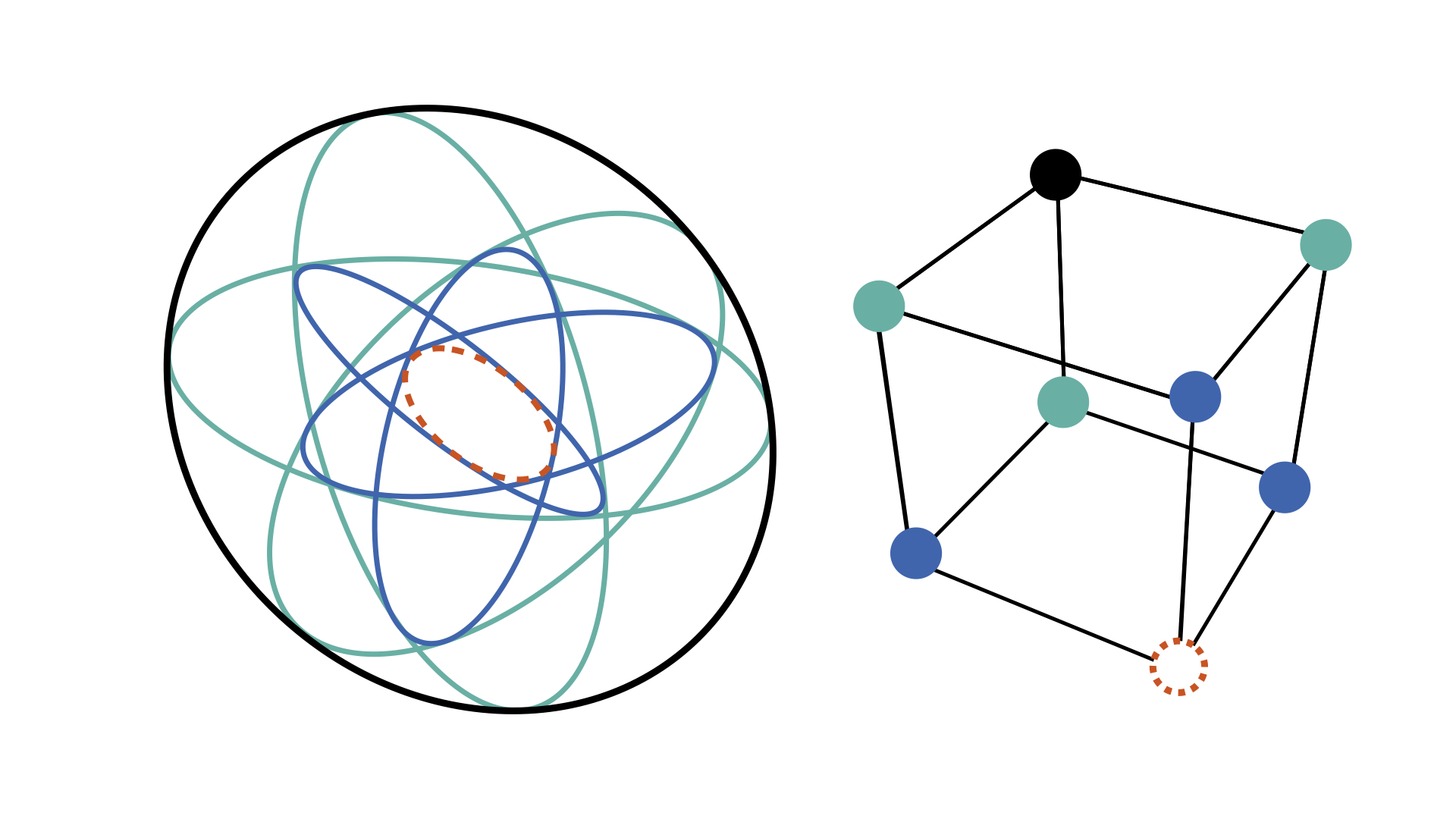}
        }
    \end{picture}
    \caption{The main theorem---the eight-conic theorem---first enunciated by Roger Penrose in around 1950.  The right figure illustrates an instance of the configuration where all 8 conics are ellipses, shown next to the corresponding cube graph of conics and double contact relationships.}
    \label{fig:PenroseInIntro}
\end{figure}

\charliex{
\section{Overview}

This article presents a ``final'' generalization of this progression of projective geometric theorems, first stated by Roger Penrose, the last author of this article, in around 1950 during his undergraduate years.  We refer to it \emph{Penrose's eight-conic theorem} (Figure~\ref{fig:PenroseInIntro}, right):
\begin{quote}
    \emph{
    If conics are assigned to seven of the vertices of a combinatorial cube such that conics connected by an edge are in double contact, and the chords of contact in each face share a common point, then there exists an  eighth conic which completes the cube under the same conditions.
    }
\end{quote}
That is, if three conics \(\ptq{T}^1, \ptq{T}^2, \ptq{T}^3\) each in double contact with a fourth \(\ptq{S}_0\), and another three conics \(\ptq{S}_1, \ptq{S}_2, \ptq{S}_3\) are, respectively, in double contact with \(\ptq{T}^2,\ptq{T}^3\), with \(\ptq{T}^3,\ptq{T}^1\), and with \(\ptq{T}^1,\ptq{T}^2\) (so that the chords of contact of a face share a point), then there exists a conic \(\ptq{T}^0\) that has double contact with \(\ptq{S}_1, \ptq{S}_2, \ptq{S}_3\). We will show below in Section \ref{sec:facefigure} that the condition on the chords of contact is equivalent to choosing the correct branch among 3 possible families of double contact conics.

The article first presents a geometric proof of the theorem, under the assumption that the conics are regular. After reviewing background material on projective geometry (Section \ref{sec:background_and_notation}), it proceeds by proving an analogous ``8-quadric'' theorem, in which conics in double contact are replaced with 3D quadrics in ring contact (Section \ref{sec:the_main_theorems}). The 8-conic theorem is proved from the 8-quadric theorem by showing that 1) the slice of a regular 8-conic configuration with a generic plane is a regular 8-conic configuration, and 2) a regular 7-conic configuration can be extruded to a regular 7-quadric configuration.  Proving the regular 8-conic theorem in this way has the advantage that in 3D there is extra geometric structure that simplifies the proof.  The idea for this proof originated with Roger Penrose during his graduate studies and was reconstructed for this article. 

The geometric proof is followed by a discussion of special 3D configurations that yield, by slicing, well-known theorems in the projective plane (Section \ref{sec:examples_of_8_quadric_configurations}), introducing the possibility of non-regular conics in the cube.

Equipped with a geometric understanding of the content of the 8-conic theorem in this weaker form, the article proceeds in Section \ref{sec:alg-proof} to an algebraic proof that applies to all conics, regular or not; for example, the line pairs and point pairs of the example theorems above, or even double lines and double points. The proof is carried out on the basis of a series of parametrized formulas relating the neighboring conics of the cube.  Section \ref{sec:det-structure} shows that all these formulas can be obtained as subdeterminants of a special $4 \times 4$ symmetric matrix, some of whose entries are linear or quadratic polynomials.  The symmetric subdeterminants correspond to the conics of the Penrose cube, while the other subdeterminants contain accessory information, such as the chords of contact, or even secondary theorems not contained in the 8-conic theorem itself. 

Section \ref{sec:rere_history} of the article gives a brief account of previous research and of the genesis of this article. 
Finally, Section \ref{sec:outlook} provides a recapitulation of the article integrated with a discussion of directions for future research. 
\ifthenelse{\boolean{isCoda}}
{Section \ref{sec:coda} closes the orticle with reflections on the archetypal nature of the theorem and its significance for projective geometry as a whole.
}{}
}
\section{Background and notation}
\label{sec:background_and_notation}

\subsection{Conics and the projective plane}

We recall that a point $P$ in the projective plane $\RP{2}$ can be represented by an equivalence class of vectors in a 3-dimensional real\footnote{We work primarily over the real field but allow complex coefficients when they arise through algebraic closure.} vector space $\mathtt{V}$:
\begin{equation}
P=\{cV:  c \in \mathbb{R}\setminus \{0\}\}
\end{equation}
for $V \in \mathtt{V} \setminus \{0\}$. A line $p$ can be described as its set of incident points by means of a vector $v \in \mathtt{V}^* \setminus \{0\}$ in the dual vector space of $\mathtt{V}$, the space of linear functions $\mathtt{V} \rightarrow \mathbb{R}$: 
\begin{equation}
p=\{P \in \RP{2}: v(V)=0, V \in P\}.
\end{equation}
Points are said to be ``collinear'' if they lie on the same line; lines are ``copunctual'' if they pass through the same point. 

The lines also naturally form a projective $2$-space called the projective dual of ${\RP{2}}^*$ by representing them as equivalence classes of dual vectors:
\begin{equation}\label{lines}
p=\{cv:  c \in \mathbb{R}\setminus \{0\}\}.
\end{equation}
for $v \in \mathtt{V}^* \setminus \{0\}$. We could equally well have started by describing lines according to Eq. \eqref{lines} and proceeded to points by describing them as the sets of incident lines:
\begin{equation}
P=\{p \in \RP{2}: v(V)=0, v \in p\}
\end{equation}
according to the fact that for a finite dimensional vector space, the dual space of the dual space is isomorphic to the original space according to the mapping $v \mapsto f_V$:
\[
f_V(v)=v(V).
\]
thus, the roles of the two are completely interchangeable: the principle of duality.  The symmetric relationship of the two can be expressed with the notation
\[
\langle v, V \rangle=v(V).
\]

A point-line correlation is an invertible map $\RP{2} \mapsto {\RP{2}}^*$ induced by a linear transformation $\mathtt{V} \mapsto \mathtt{V}^*$.  The fundamental theorem of projective geometry states that this is the most general type of mapping between the two, which always sends any three collinear points to three copunctual lines.  A line-point correlation is a similarly defined map ${\RP{2}}^* \mapsto {\RP{2}}$.  

Given an ordered basis $\mathcal{B}=(V_0,V_1,V_2)$ of $\mathtt{V}$, representing points $(P_0,P_1,P_2)$, the dual basis $\mathcal{B}^*=(v_0,v_1,v_2)$,
\[
\langle v_j,V_k \rangle = \begin{cases} 
0 & \text{if } k \neq j \\
1 & \text{if } k=j,
\end{cases}
\]
represents the lines $(p_0,p_1,p_2)$ constituting the sides of the triangle $(P_0,P_1,P_2)$ with $p_j$ opposite to $P_j$.  For simplicity of notation and to avoid significant repetition, we shall often identify the vectors $V$ or $v$ with the geometric elements $P$ or $p$, keeping in mind that specific representing vectors are always implied since otherwise the basis is not well defined.  The point represented by $aP_0+bP_1+cP_2$ is often denoted by the notation $[a:b:c]$ referred to as {\it homogeneous coordinates}.

The linear transformation defining a point-line correlation can then be represented by a matrix $A$ relative to the bases $\mathcal{B}$ and $\mathcal{B}^*$, and vice versa for a line-point correlation.  The matrix $cA$ for $c \neq 0$ clearly represents the same correlation (note, however, that changing the relative scaling of the vectors representing $P_j$ leads to a different matrix).  A point-line correlation has two naturally associated line-point correlations: its inverse, which has matrix $A^{-1}$, and the mapping of a line $p$ to the intersection point of any two lines mapped to by two points lying on $p$, which is well-defined and has matrix $\mathrm{cof} A$, the cofactor matrix of $A$.  By Cramer's rule,
\begin{equation}\label{cramer}
A^{-1}=\frac 1 {\det A} \mathrm{cof} A^T
\end{equation}
where the superscript $^T$ denotes the transpose, therefore these two associated line-point correlations coincide exactly when the point-line correlation is represented by a symmetric matrix $A=A^T$ which is known as a {\it polarity}. A polarity can also be characterized as a correlation that is also an involution.  A line and point which are mapped to each other by a polarity are referred to as  \emph{polar to} each other.

If $A$ is the matrix of a polarity and $P$ is a point, which we identify with the column vector of a representative in the basis $(P_0,P_1,P_2)$, then the line $p$ polar to $P$ is represented by a column vector, which we again identify with $p$ itself, in the basis $(p_0,p_1,p_2)$:
\[
p=AP
\]
The condition that a point is incident to its polar line may thus be expressed by the matrix equation
\[
p^TP=P^T A P=0
\]
The locus of such points is a {\it regular point-wise conic} which is uniquely defined by the polarity and vice versa, both having a unique representative symmetric matrix (up to scaling) in a given basis.  The line to which a point of a point-wise conic is both polar and incident is its {\it tangent line}.   Dual to this, the locus of lines incident to their polar points is called a {\it regular line-wise conic} and we refer to these points as the {\it turning points} of the line-wise conic.  When convenient and no ambiguity can result, we sometimes identify the matrix representing a conic in a given basis with the conic itself.

Given a regular point-wise conic, there is a uniquely associated regular line-wise conic, the lines of the latter being the tangent lines of the points of the former. 
If the first conic is represented by the matrix $A$, the second is represented by the inverse matrix $A^{-1}$.
A {\it degenerate point-wise conic} is a conic defined in the same way but by a non-invertible symmetric matrix.  In this case, three collinear points will still map to three copunctual lines when all images are defined, however, for any point $P$ whose representatives lie in the kernel of the defining matrix $A$, the image of $P$ under the {\it degenerate polarity} is geometrically undefined, since the zero vector does not represent any line in ${\RP{2}}^*$, moreover, if $Q$ is any other point, then
\[
A(cP+Q)=AQ
\]
for any scalar $c$, so any point lying on a line through $P$ will be mapped to the same line.  
\begin{itemize}
\item If $\mathrm{ker} A$ has dimension $1$ then it constitutes one equivalence class in $\RP{2}$ so there is a unique such point.  In this case, the points of a degenerate point-wise conic constitute a {\it line pair} with $P$ as their meeting point, called the {\it carrying point} of the line pair.  If these two lines are called $p$ and $q$ then the matrix is given by
\[
A=pq^T+qp^T
\]
\item If $\mathrm{ker} A$ has dimension $2$ then it constitutes the representatives of the points on a line in the projective plane.  In this case the points of a degenerate point-wise conic constitute a line, referred to for algebraic reasons as a {\it double line}.  If this line is called $p$ then the matrix is given by
\[
A=pp^T
\]
\end{itemize}
Dual to this, degenerate line-wise conics can be defined:
\begin{itemize}
\item If $\mathrm{ker} A$ has dimension $1$ then it constitutes one equivalence class in ${\RP{2}}^*$ so there is a unique line $p$ represented by vectors in the kernel.  In this case the lines of a degenerate line-wise conic constitute a {\it point pair} with $p$ as their joining line, called the {\it carrying line} of the point pair.  If these two points are called $P$ and $Q$ then the matrix is given by
\[
A=PQ^T+QP^T
\]
\item If $\mathrm{ker} A$ has dimension $2$ then it constitutes the representatives of the lines through a point in the projective plane.  In this case, the lines of a degenerate line-wise conic constitute a point, referred to for algebraic reasons as a {\it double point}.  If this point is called $P$ then the matrix is given by
\[
A=PP^T
\]
\end{itemize}

\subsubsection{Pencils of conics}
The linear span of two conics $A$ and $B$ is called a {\it pencil} of conics. A pencil that contains a regular conic is called a {\it regular} pencil. We will need the following lemma in Sect. \ref{sec:facefigure}:
\begin{lemma}[Regular pencil] \label{lemma:regpencil}
A regular pencil of conics contains at most 3 distinct degenerate conics.
\end{lemma}
\begin{proof}
   Assume $A$ is regular. The pencil, with the exception of $A$, can be written $A+tB, t\in \mathbb{R}$. $A+tB$ is a degenerate conic exactly when $|A+tB|= 0$, where vertical bars around a matrix denote the determinant. The left-hand side is a cubic polynomial in $t$, hence has at most 3 roots.   
\end{proof}
\subsubsection{Double contact}
For a regular point-wise conic $A$ and a double line $pp^T$, the polarity defined by the conic
\begin{equation}\label{dcp}
B=A+cpp^T
\end{equation}
acts the same on the points of $p$, independently of the parameter $c$:
\begin{equation}\label{cc_polarity}
BP=AP+cpp^TP=AP
\end{equation}
when $P$ is incident to $p$.  In particular, the intersection points of $B$ with $A$ are always the intersection points of $p$ with $A$ and the tangent lines there are always the same, independently of $c$.  Thus, $A$ and $B$ are said to be in {\it double contact} with the {\it chord of contact} $p$, and the pencil of conics defined by Eq. \eqref{dcp} is called a {\it double contact } pencil.
  The property of being in double contact is clearly self-dual. Indeed, two point-wise conics share tangent lines at two common points if and only if the two corresponding line-wise conics share turning points at two common lines. Furthermore, Eq. \eqref{dcp} is self-dual, that is, the line-wise conics of a double contact pencil can be expressed as the linear span of a regular conic and a double point, given by the meet of the two lines, as can be checked by taking cofactor matrices.  We call this double point the {\it point of contact}.  Therefore, the polarity property expressed by Eq. \eqref{cc_polarity} can be stated as
  \begin{lemma}\label{dc_polarity}
  The polarities on two conics in double contact act the same when restricted to a mapping between the points of the chord of contact and the lines of the point of contact.
  \end{lemma}

\subsection{Quadrics and projective space}
\label{quadrics_intro}
Similar definitions hold in projective 3-space $\RP{3}$, using a 4-dimensional vector space $\mathtt{V}$ where now the dual space ${\RP{3}}^*$ is the space of planes.  A correlation is defined by the property that 3 collinear points map to 3 coaxial planes (three planes going through the same line: we use the word "axis" to emphasize the dual role that a line plays as the intersection of two planes in contrast to its more familiar role of connecting two points, for which we use the term ``spear'').  A polarity is defined the same way, but we now refer to the surface of points incident to their polar planes as a point-wise {\it quadric}.  There are now 3 different cases to consider:
\begin{itemize}
\item If the kernel $\mathrm{ker} A$ has dimension $1$, it represents the vertex of a cone that constitutes the points of the quadric.
\item If the kernel $\mathrm{ker} A$ has dimension $2$, it represents the axis of intersection of a plane pair.  If these planes are called $p$ and $q$ then,
\[
A=pq^T+qp^T.
\]
\item If the kernel $\mathrm{ker} A$ has dimension $3$, it is a (double) plane.  If this plane is called $p$ then
\[
A=pp^T
\]
\end{itemize}
and for plane-wise quadrics
\begin{itemize}
\item If the kernel $\mathrm{ker} A$ has dimension $1$, it represents the plane of a conic curve, the planes going through a tangent line to this curve constitute the planes of the quadric.
\item If the kernel $\mathrm{ker} A$ has dimension $2$, it represents the joining line of a point pair.  If these points are called $P$ and $Q$ then
\[
A=PQ^T+QP^T.
\]
\item If the kernel $\mathrm{ker} A$ has dimension $3$, it is a (double) point.  If this point is called $P$ then
\[
A=PP^T
\]
\end{itemize}

The set of geometric elements represented by vectors in the kernel of $A$ is known as the {\it degenerate set}. For point-wise quadrics this is: the vertex of a cone, the points of the axis of a plane pair, the points of a double plane; and for plane-wise quadrics: the plane of a conic, the planes of the joining line of a point pair, the planes of a double point.  For degenerate point-wise quadrics, the polarity mapping is only defined for points not in the degenerate set to planes that go through the degenerate set, and vice versa for plane-wise quadrics.
For point-wise quadrics:
\begin{itemize}
    \item The polarity on a cone maps a point which is not the vertex of the cone to a plane through the vertex, and it maps the vertex to $0$ which doesn't represent any geometric element.
    \item The polarity on a plane pair maps a point not on the axis of the plane pair to a plane through the axis and any point on the axis to $0$.
    \item The polarity on a double plane maps any point not on the plane to the plane and any point on the plane to $0$.
\end{itemize}
For plane-wise quadrics:
\begin{itemize}
    \item The polarity on a conic maps a plane which is not the plane of the conic to a point on the plane of the conic (in particular, the point polar to the intersection line of this plane with the plane of the conic), and it maps the plane of the conic to $0$.
    \item The polarity on a point pair maps a plane that does not go through the connecting line of the point pair to a point on the connecting line, and any plane through the connecting line to $0$.
    \item The polarity on a double point maps any plane not through the point to the point and any plane through the point to $0$.
\end{itemize}
Taking this into account, it is still possible to define the binary relation that {\it a point and a plane are polar to each other with respect to a degenerate point-wise quadric} if the point is mapped to the plane by the polarity and {\it a plane and a point are polar to each other with respect to a degenerate plane-wise quadric} if the plane is mapped to the point by the polarity. 

\subsubsection{Ring contact}
The condition analogous to double contact is {\it ring contact}: the quadrics
\begin{equation}\label{rcp}
B=A+cpp^T
\end{equation}
meet in the conic where $p$ intersects $A$ and they have the same tangent planes at all these shared points.  We call $p$ the {\it ring plane} and  refer to its polar point as the {\it ring points}.  By the same reasoning as in the conic case, we have
\begin{equation}
BP=AP+cpp^TP=AP
\end{equation}
when $P$ is incident to $p$ and thus

\begin{lemma}\label{rc_polarity}
  The polarities on two quadrics in ring contact act the same when restricted to a mapping between the points of the ring plane and the planes of the ring points.
  Conversely, if two polarities coincide on a plane, then the corresponding quadrics are in ring contact with respect to this plane.
\end{lemma}
  
\begin{remark}
    In the case where the ring point is incident to the ring plane, that is, the ring plane is a tangent plane and the ring point is its point of contact, the conic of the ring contact is a (double) pair of lines.  This case can be distinguished as {\it X-contact} although we include it under the umbrella term ``ring contact''.
\end{remark}

\subsubsection{Some contact lemmas}

We provide a few lemmas that will be useful in the proof of the 8-quadric theorem.
By {\it ring contact triangle} we mean a set of 3 ring contact pencils that pairwise have one conic in common.

\begin{lemma}
\label{noringcontacttriangles}
If 3 quadrics are mutually in ring contact and at least one is regular, then all 3 are in the same ring contact pencil, i.e., there are no ring contact triangles.
\end{lemma}

\begin{proof}
If $A, B$ and $C$ are in ring contact, then there exist planes $p,q$ and $r$ such that
\[
B=A+pp^T, \qquad C=B+qq^T,\quad \mathrm{and} \quad \qquad A = \alpha  C+rr^T
\]
where $\alpha$ is a scalar (note that we are not entitled to arbitrarily scale it to $1$ because we already fixed the relative scaling of $A$ and $C$ with the first two equations).

Then
\[
C=A+pp^T+qq^T
\]
implying that
\[
(1-\alpha)C=pp^T+qq^T+rr^T.
\]
However, the right hand side has rank at most 3: it is either a cone with the intersection point of $p, q$ and $r$ as its vertex, or they meet in an axis through which it is a plane pair, or all three planes coincide.  But we may, without loss of generality, choose $C$ to be the regular quadric and thus have rank $4$, implying that $\alpha=1$.  Thus the left side vanishes and then, because 3 double planes may only be dependent if they coincide, $p, q$ and $r$ are necessarily the same plane.
\end{proof}

We remark that the same cannot be said of conics in double contact since rank 3 conics are nondegenerate, i.e., double contact triangles can exist, for instance
\[
A=\begin{pmatrix}
    -1 & 0 & 0\\
    0 & 1 & 0\\
    0 & 0 & 1
\end{pmatrix}, \qquad
B=\begin{pmatrix}
    -1 & 0 & 0\\
    0 & 2 & 0\\
    0 & 0 & 2
\end{pmatrix}, \qquad
C=\begin{pmatrix}
    -1 & 0 & 0\\
    0 & 2 & 0\\
    0 & 0 & 1
\end{pmatrix}.
\]

On the other hand, if the chords of contact $p,q,r$ all meet in a point, then the span of $pp^T,qq^T,rr^T$  
contains only line pairs through this common point, so we have the following lemma: 
\begin{lemma}
\label{noconcdoublecontacttraingles}
If 3 conics are pairwise in double contact with concurrent chords of contact, then all 3 are in the same double contact pencil.
\end{lemma}

\subsection{Notation}
\newcommand{\myln}[1]{\mathbf{#1}}
\newcommand{\abpr}[2]{\myln{#1#2}}   
\newcommand{\abdb}[1]{\myln{#1}^2}

Before we turn to the proof of the theorems, we introduce some notation used in the sequel. 
First, an alternative notation for degenerate conics, used in Sec. \ref{sec:facefigure} and Sec. \ref{sec:alg-proof}. Secondly, notation for so-called \emph{complete} quadrics, required to represent non-regular conics and quadrics, and used in Theorem. \ref{thm:refinement}.

\subsubsection{Abbreviated notation}
\label{sec:abbr-nota}

Until now, we have written a line pair formed by the two lines $p$ and $q$ directly as a matrix via the symmetrized outer project $pq^T+qp^T$, and the double line of $p$ as $pp^T$. 
$p$ and $q$ can also be interpreted as the zero sets of the linear polynomials 
\begin{equation*}
    p(x,y,z) := ax+by+cz,  ~~~~~~~ q(x,y,z) := a'x+b'y+c'z 
\end{equation*}Their product is the quadratic polynomial
\begin{align*}
r(x,y,z):&= p(x,y,z)q(x,y,z) \\
&=aa'x^2 + bb'y^2+cc'z^2 +(ab'+a'b) xy + (ac'+a'c) xz + (bc'+b'c) yz  \\
&= (x,y,z)^T(pq^T+qp^T)(x,y,z)
\end{align*}
The final equality shows that $r(x,y,z)=0$ represents the same line pair as $(pq^T+qp^T)$. Setting $q=p$ yields an expression for the double line $p^2$ which agrees with the matrix $pp^T$.  We introduce the expressions $\abpr{p}{q}$ and $\abdb{p}$ as the {\it abbreviated notation} for the line pair, resp., the double line. 
We also write $\myln{p}$ for the line itself (not a double line); although it is not a conic, it will be used in Section \ref{sec:det-structure} to generate conics within the ring of polynomials in the variables $\{x,y,z\}$. We continue to use $p$ to denote a line in $\RP{2}$ when the membership in the polynomial ring is not highlighted.

Let $A$ and $B$ be two points. Similar remarks lead to the abbreviated notation $\abpr{A}{B}$, resp., $\abdb{A}$ for the point pair formed by $A$ and $B$, resp., the double point $A$. 

Finally, in the sequel we also use $\ptc{X}$ to represent the conic or quadric with matrix $X$. This ensures that bold-face font signals a conic, whether as $\ptc{X}$ or the abbreviated forms $\abpr{p}{q}$ and $\abdb{p}$.  

The advantage of the abbreviated notation is that it gives access to the algebraic structure of the ring of polynomials when evaluating linear combinations of conics involving degenerate conics, see for example Sec. \ref{sec:facefigure}. It was introduced by Salmon \cite{salmon1917}.

\subsubsection{Notation for complete quadrics}
\label{sec:nota-for-comp-conics}

Every quadric has a point-wise, or primary, aspect as well as a plane-wise, or dual, aspect. 
We write the primary (dual) form of a quadric as $\ptq{Q}$ ($\plq{Q}$). When does the pair $(\ptq{Q}, \plq{Q})$ represent a quadric?  When either is regular, then so is the other, and the relationship is given by $\ptq{Q}\plq{Q} = \lambda Id$ for $\lambda \neq 0$ by equation \eqref{cramer}, where by $\ptq{Q}$ and $\plq{Q}$, we mean some matrices representing these quadrics. Otherwise, both are degenerate and satisfy  the relationship $\ptq{Q}\plq{Q} = \lambda Id$ for $\lambda = 0$, which can be seen by taking limits of regular quadrics. Then a necessary condition for $(\ptq{Q}, \plq{Q})$ to represent a quadric is $\ptq{Q}\plq{Q} = \lambda Id$. If in addition both are not zero, we call the pair a \emph{complete} quadric. In discussions limited to regular quadrics, we usually refer to a quadric simply in its point-wise form as $\ptq{Q}$, since $\plq{Q}$ is uniquely determined. When degenerate quadrics play a role, as in Sec. \ref{thm:refinement}, we employ the full notation introduced here.  

All of the above also applies to conics, if \emph{plane-wise} is replaced by \emph{line-wise}.

\section{The main theorems}
\label{sec:the_main_theorems}

We state the 8-conic and 8-quadric theorems in their full generality now.

\begin{theorem}[8-conic theorem] \label{thm:8conic}
In $\RP{2}$, if conics are assigned to seven of the vertices of a combinatorial cube such that (i) conics connected by an edge are in double contact, and (ii) the chords of contact associated to a cube face meet in a common point, then there exists an eighth conic such that the completed cube satisfies (i) and (ii). This conic is unique if not all the common points of condition (ii) are the same. 
\end{theorem}

\begin{theorem}[8-quadric theorem]   
In $\RP{3}$, if quadrics are assigned to seven of the vertices of a combinatorial cube such that (i) quadrics connected by an edge are in ring contact, and (ii) the ring planes associated to a cube face meet in a common axis, then there exists an eighth quadric such that the completed cube satisfies (i) and (ii). This quadric is unique if not all the common axes of condition (ii) are the same.  \label{thm:8quadric}
\end{theorem}

We begin with a geometric proof under the weaker condition that the conics are regular. This will be followed in Section \ref{sec:alg-proof} with an algebraic proof that removes this restriction. We first prove Thm. \ref{thm:8quadric} and use it to prove Thm \ref{thm:8conic}.

\subsection{Notation and terminology}

The 8-conic and 8-quadric theorems are concerned with assigning conics (or quadrics) to (some or all of) the vertices of a combinatorial cube graph so that they satisfy certain edge and face relations. We introduce common terminology and notation for these configurations.
\begin{definition}
    A \emph{Penrose $n$-configuration} is an assignment of $n$ conics (quadrics) to the vertices of a cube.
    The configuration is \emph{regular} if the conics (quadrics) are regular. 
\end{definition}

Both the 8-conic and 8-quadric theorems assume that a Penrose 7-configuration is given and assert a unique completion to an 8-configuration. 

We can immediately count that a Penrose 7-configuration has 4 completed vertices, 9 completed edges, and 3 completed faces. 

\begin{figure}[h!]
    \centering
    \includegraphics[width=0.33\textwidth]{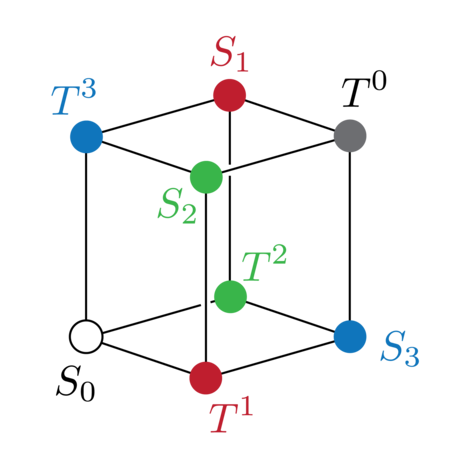}
    \hspace{1cm}
    \includegraphics[width=0.44\textwidth]{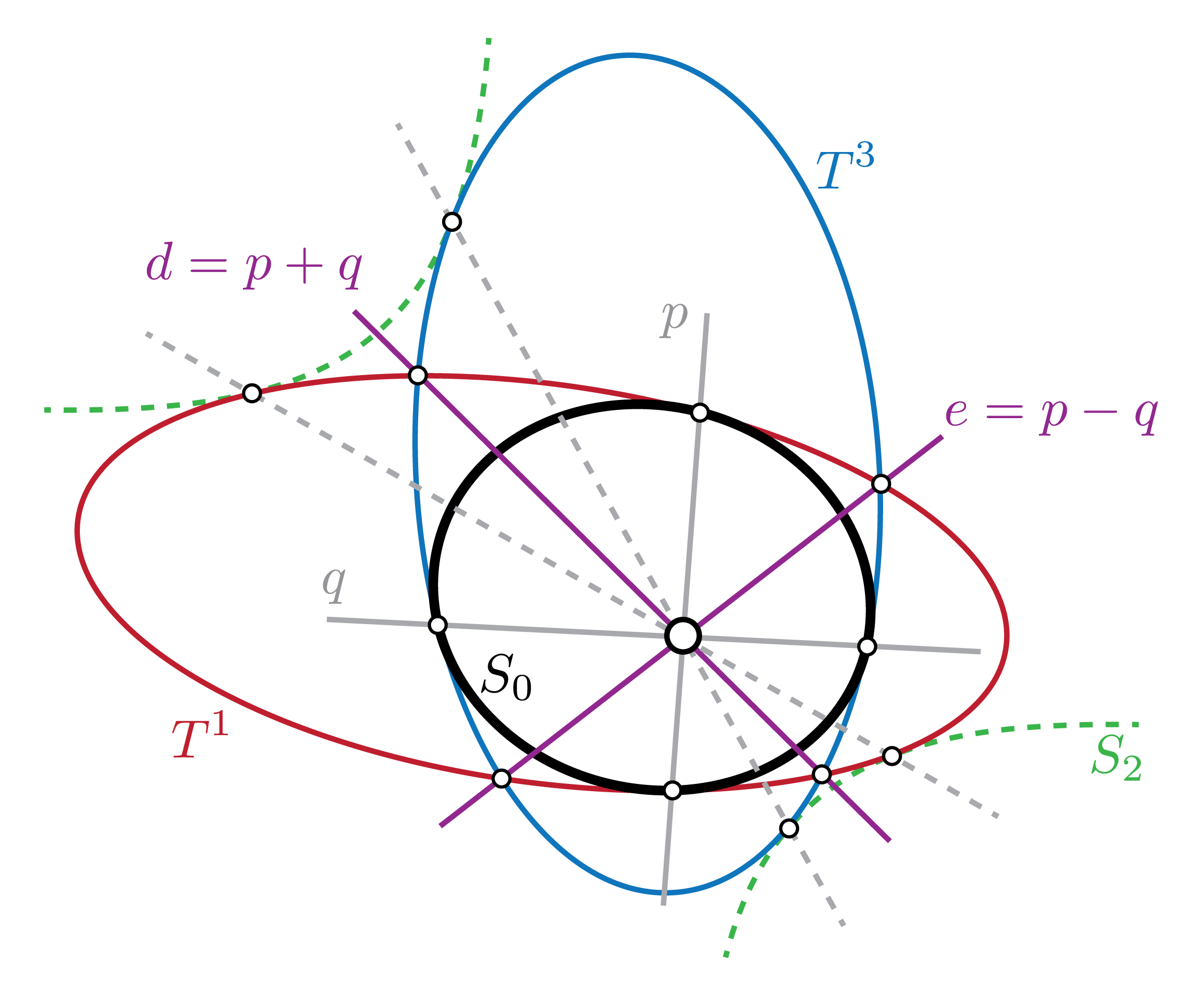}
    \caption{The cube graph of the 8-conic theorem (left) and an example of 4 conics $\ptq{S}_0,\ptq{S}_2,\ptq{T}^1,\ptq{T}^3$ satisfying condition (ii) (right).}
    \label{fig:cube+facecondition}
\end{figure}

Next, we introduce notation for the vertices and faces of the cube. 
Consult Figure \ref{fig:cube+facecondition}, left.    
We decompose the cube vertices into two interlocked tetrahedra and label the first group of vertices $\ptq{S}_i, i\in(0,1,2,3)$ and the second $\ptq{T}^i$, where $\ptq{T}^i$ is opposite to $\ptq{S}_i$ with respect to the center of the cube.
Write the face of the cube with corners $\ptq{T}^i\ptq{S}_k\ptq{T}^j\ptq{S}_m$ as $\face^{ij}_{km}$ where $ijkm$ is a permutation of $0123$. The common point of the chords of contact in this face we write as $P^{ij}_{km}$, the \emph{face point} of the face.  Since the edges of the cube represent double contact pencils, each such pencil contains a unique double line $\VC^{i}_j=(\myln{p}^i_j)^2$, where $\myln{p}^i_j$ the chord of contact of $\ptc{T}^i$ and $\ptc{S}_j$.

When working with the 8-quadric theorem, we call the common axis of the ring planes of the face $\face^{ij}_{km}$ the \emph{face axis} $\axis^{ij}_{km}$.  Furthermore, in what follows, we assume that the empty vertex of a Penrose 7-configuration is $\ptq{T}^0$.


\subsection{The face figure of a Penrose conic configuration} \label{sec:facefigure}

In Figure \ref{fig:cube+facecondition} (left), consider two conics corresponding to opposite nodes in a face of the cube, for example, $\ptq{T}^1$ and $\ptq{T}^3$. Then by \eqref{dcp} they can be written as $\ptq{T}^1=\ptq{S}_0+\abdb{p}$ and $\ptq{T}^3=\ptq{S}_0+\abdb{q}$ with $\myln{p}$ and $\myln{q}$ being the chords of contact between $\ptq{S}_0$ and $\ptq{T}^i$. Using abbreviated notation: 
$$\ptq{T}^1-\ptq{T}^3=\abdb{p} - \abdb{q} = (\myln{p} + \myln{q})(\myln{p} - \myln{q}) =: \abpr{d}{e}$$
where $\myln{d} :=  \myln{p} + \myln{q}$ and $\myln{e} = \myln{p} - \myln{q}$ is a line pair in the pencil spanned by $\ptq{T}^1$ and $\ptq{T}^3$. Hence $\myln{p},\myln{q},\myln{d}$ and $\myln{e}$ have a common point.  Consult Figure \ref{fig:cube+facecondition}, right.  

By Lemma \ref{lemma:regpencil}, there are at most 3 line pairs in a regular pencil of conics. Hence, the set of conics that are in double contact with $\ptq{T}^1$ and $\ptq{T}^3$ splits into at most 3 1-parameter families such that condition (ii) of the 8-conic theorem holds automatically for all conics $\ptq{S}_0, \ptq{S}_2$ in the same family. Each 1-parameter family can be explicitly written in terms of $\myln{p}$, $\myln{q}$, and $\ptq{S}_0$ (see \cite[\S 287]{salmon1917}). Given $\ptc{S}_1$ in such a family, its parameter value can be obtained directly from this explicit form.

\subsection{Penrose's 8-conic theorem}

The connection of the 8-conic theorem and the 8-quadric theorem is based on the following observation: 
Given two regular quadrics in ring contact and a plane that is not tangent to either quadric, slicing the quadrics with the plane yields two regular conics in double contact. In this slice, the common axis of the ring planes is cut in the common point of the chords of contact of the conics. Conversely, if two regular conics in double contact are given, it is possible to extrude them to obtain two regular quadrics in ring contact. (Similar remarks apply to the plane-wise view.)

The extrusion needed for the 7-configurations encountered in the following proofs is provided by the following lemma, whose proof appears below in Sect. \ref{sec:extrusion}. 

\begin{lemma}[Extrusion] \label{lemma:extrusion}
Given a regular Penrose 7-conic  configuration, it is possible to extrude it to a regular Penrose 7-quadric configuration. 
\end{lemma}

 
\begin{proof}[Proof of Penrose-8 conic theorem]

Apply Lemma  \ref{lemma:extrusion} to the given regular Penrose 7-conic configuration to obtain an extrusion to a regular Penrose 7-quadric configuration. Then apply Thm. \ref{thm:8quadric} to complete this to a Penrose 8-quadric configuration. The intersection of the plane of the original 7-conic configuration with the complete 8-quadric configuration is the desired regular Penrose 8-conic configuration. Face points of the conic configuration are distinct exactly when the corresponding face axes of the quadric configuration are distinct. Hence the eighth conic is unique exactly when  the face points are not all the same.
\end{proof}

\myboldhead{Remark}  Due to factors explained in the next section, the proof presented below for the 8-quadric theorem does not cover the case that all the face axes, resp., points, are identical, hence the 8-conic theorem proof above does not cover the case that all chords of contact of the configuration pass through a single point.

\subsection{Penrose's 8-quadric theorem}
 \label{sec:proof8quadric}

We begin with three lemmas.


Let $i$ and $j \in \{0,1,2,3\}$, $m$ and $n$ in $\{1,2,3\}$ and $P^m_i$ denote the ring point for $\ptq{S}_i$ and $\ptq{T}^m$, $p^m_i$ denote the corresponding ring plane, $\axis^{mn}_{ij}$ the common axis of the ring planes for the face $\mathcal{F}^{mn}_{ij}$ and $\spear^{mn}_{ij}$ the common spear (point-wise line) of the  ring points.  We call these the face axes, resp., the face spears.  By assumption, they are well defined.  Equivalently, given any one of the three complete faces of the 7-configuration, its four ring planes are not all identical.  Consider, for instance, the face $\face^{12}_{03}$.  By Lemma \ref{noringcontacttriangles}, we immediately conclude that $\ptq{S}_0$ is not in ring contact with $\ptq{S}_3$ and $\ptq{T}^1$ is not in ring contact with $\ptq{T}^2$.  By similar reasoning applied to the other two complete faces, we get the following lemma:

\begin{lemma}\label{lemma:7noextraring}
    $\ptq{S}_i$ is not in ring contact to $\ptq{S}_j$ and $\ptq{T}^m$ is not in ring contact to $\ptq{T}^n$ for any of the seven quadrics of a regular Penrose 7-configuration.
\end{lemma}

Since a face axis is determined by any two ring planes on a face, and every face of the 7-configuration has at least two ring planes, the axes and spears of the incomplete faces are also uniquely determined.  In particular, the shared axis of the ring planes for the incomplete face $\mathcal{F}_{ij}^{0n}$ must go through the ring planes for the two edges already given, so we set
\begin{equation}\label{axis0nij_def}
\axis^{0n}_{ij}:=p^n_i \wedge p^n_j
\end{equation}
and dually
\begin{equation}
\label{spear0nij_def}
\spear^{0n}_{ij}:=P^n_i \vee P^n_j.
\end{equation}
These will both be non-zero since otherwise, if we had $P^n_i=P^n_j$ or $p^n_i=p^n_j$, $\ptq{S}_i$ and $\ptq{S}_j$ would have to be in ring contact along the common conic of intersection with $\ptq{T}^n$, contradicting Lemma \ref{lemma:7noextraring}.

We now demonstrate a global condition satisfied by a regular quadric 7-configuration which has no analogy in the case of conics. 
\begin{lemma}[Concurrent axes and spears]\label{lem:faconcur}
    For a regular Penrose 7-quadric configuration: If two of the face axes at $\ptq{S}_0$ are distinct, then at each complete vertex, all three face axes are distinct. In this case, all ring planes and face axes are concurrent in a point $O$. Otherwise, there is a single face axis $\Omega$ common to all ring planes. Dually, all ring points and face spears are concurrent in a plane $o$ or there is a single face spear $\omega$ common to all ring points. $O$ (resp., $\Omega$) and $o$ (resp., $\omega$) are a polar pair with respect to all quadrics of the configuration.
\end{lemma}
 \begin{proof}

Assume $\axis^{23}_{01}=\axis^{13}_{02}$. Since each axis contains two of the ring planes at $\ptq{S}_0$, all three of these ring planes pass through this axis so that it is also identical to $\axis^{12}_{03}$.  

This contradicts the assumption that two of the three face axes $\axis^{ij}_{0k}$ are distinct. Hence, all three are distinct if any two are.  

We show now that the $\axis^{ij}_{0k}$ cannot all lie in a plane. Indeed, if they do, then that plane must be identical to the three ring planes $p^j_0$ but by the regularity of $\ptq{S}_0$ this would imply that $\ptq{T}^i$ is in ring contact to $\ptq{T}^j$ for $i,j=1,2,3$ contradicting Lemma \ref{lemma:7noextraring}.

Hence, the three axes at $\ptq{S}_0$ are distinct and do not lie in a plane, and the ring planes are therefore distinct, and have a common point $O$, which is therefore also a common point of the axes.
        
Apply the same argument to the other complete vertices $\ptq{T}^i$ to conclude that the axes of the faces meeting at $\ptq{T}^i$ have a common point $O'$. We show that $O'=O$. Indeed, two of these axes belong to $\{\axis^{ij}_{0k}\}$, which were shown above to meet in $O$. So $O' = O$.   Since this applies to all 3 incomplete faces, all 6 axes are concurrent in $O$. 

We now consider the case where two of the axes at $\ptq{S}_0$ are the same. Then by the above,  all three axes $\axis^{ij}_{0k}$ are the same, call it $\Omega$. By reasoning similar to the above, we can propagate this condition to the three remaining face axes, to show that they also have to be equal to $\Omega$, and consequently that all ring planes pass through $\Omega$.  

By regularity and duality the analogous result can be obtained for the ring points and the face spears. Furthermore, by Lemma \ref{rc_polarity}, $O$ (resp., $\Omega$) and $o$ (resp., $\omega$) are a polar pair with respect to all quadrics of the configuration. 

 \end{proof}

We will also need this standard geometric result (see, e.g. \cite[\S 4.3]{Ziegler}): 
\begin{lemma}[Conical Octahedron Lemma]
\label{lem:conical-octahedron}
Three conics in $\RP{3}$, without a common point, that intersect pairwise in two points (including the possibility that the two points coincide and the conics thus share a tangent line), lie on a unique quadric.
\end{lemma}
Note that the conics cannot lie in a plane since two of such conics have 4 points in common.  The reader is also referred to the ``coordinate version'' of the 8-quadric proof below, which provides an independent proof of the lemma.

\subsubsection{The proof}

We are now in a position to prove the 8-quadric theorem, subject to the following restrictions.

First, we assume that two of the face axes at $\ptq{S}_0$  are distinct. By Lemma \ref{lem:faconcur}, the alternative is that all six face axes are the same. This case is handled in Appendix \ref{app:binique}.

We also assume for the moment that $o$ is not incident to $O$.  This is equivalent to assuming that the 7 quadrics of the configuration do not all share a common contact element.  We postpone to Sect. \ref{Ooincidnetproof} the proof that this condition is not necessary, since it depends on elements of the proof of Lemma \ref{lemma:extrusion}.

Figure \ref{fig:Penrose3D} depicts the completed cube configuration that we must construct. There are 6 axes corresponding to condition (ii), one for each face. Each face axis is the common line of the 4 ring planes on the edges of the face.  As mentioned above, we restrict our attention to the first case in Lemma \ref{lem:faconcur}, whereby the 3 axes of the complete faces of the 7-configuration are not all identical and have a common point $O$, as shown in Figure \ref{fig:ContactRings}.  The dual plane $o$ is taken as the plane at infinity in this figure, thus the cones of ring contact are all cylinders through the conics $\ptq{R}^k_j$ where $\ptq{S}_j$ and $\ptq{T}^k$ make ring contact. The 4 ring points on a face lie on six distinct spears in $o$.

\begin{proof}[Proof of the Penrose 8-quadric theorem]

Under the assumption that two of the faces axes at $\ptq{S}_0$ are distinct,
it suffices to construct a unique quadric $\ptq{T}^0$ such that: for each $j \in \{1,2,3\}$,  (i) $\ptq{T}^{0}$ is in ring contact to $\ptq{S}_{j}$; and (ii) all 4 ring planes for the face $\face^{0j}_{km}$ have a common axis.

\begin{figure}
    \centering
    \includegraphics[width=0.45\textwidth]{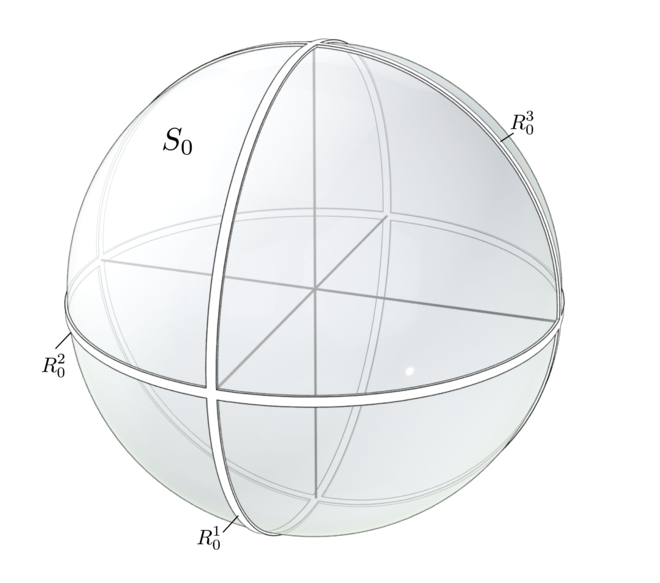}
    \includegraphics[width=0.45\textwidth]{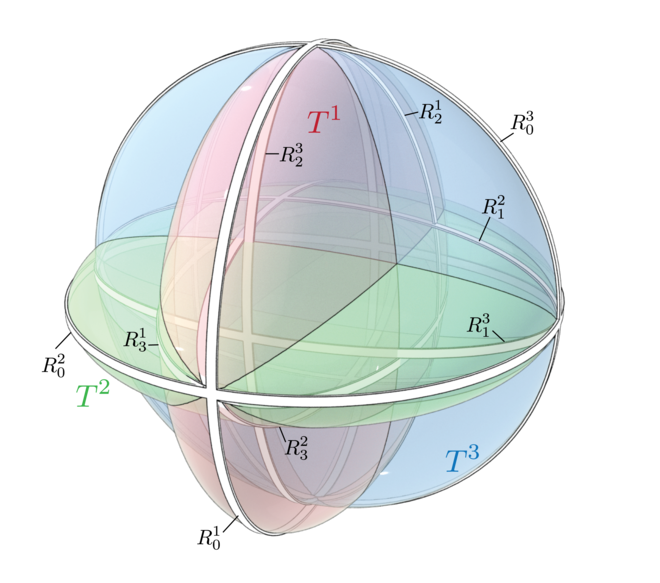}
    \includegraphics[width=0.45\textwidth]{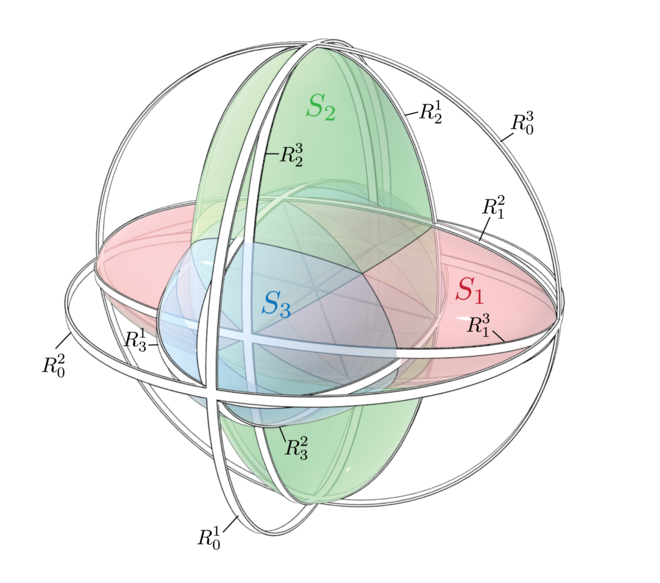}
    \includegraphics[width=0.45\textwidth]{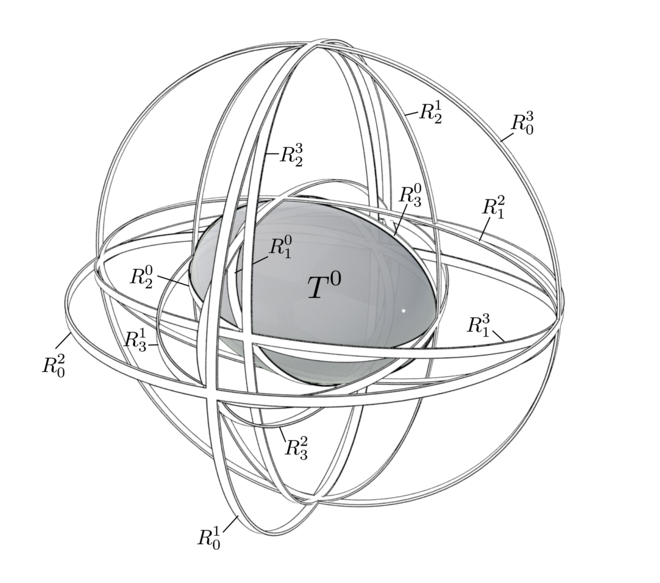}
    \linebreak
    

    \includegraphics[width=0.6\textwidth]{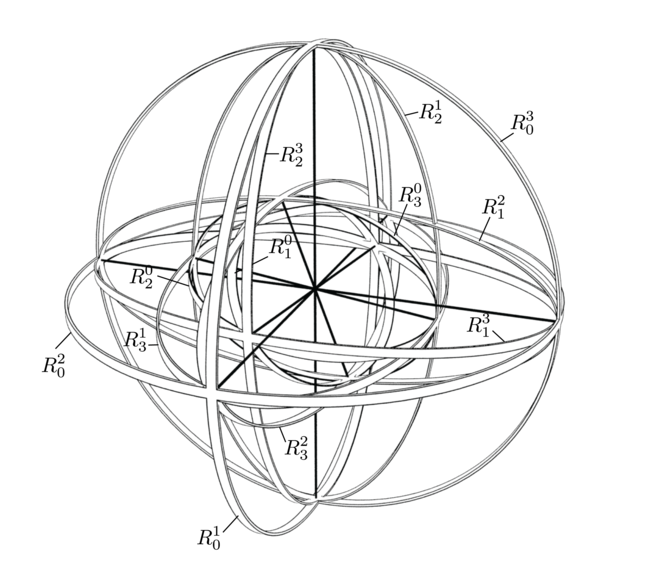}
    \caption{Eight quadrics in ring contact and the structure of all 12 contact rings $\ptq{R}_i^k$ together with the 6 axes $\axis_{ij}^{kn}$, sharing the point $O$. On each axis the four rings corresponding to a face of the cube diagram meet.}
    \label{fig:Penrose3D}
    \label{fig:ContactRings}
\end{figure}

To begin with, we consider the planes and points that must necessarily be the ring planes and ring points for such a quadric to complete the cube.  The ring plane for $\ptq{S}_i$ and $\ptq{T}^0$ will have to go through these axes for both faces adjacent to the desired edge representing ring contact:
\begin{equation}
\label{p0i_def}
p^0_i:=\axis^{0n}_{ij} \vee \axis^{0j}_{in}
\end{equation}
(note: the joining plane of the two lines is defined because they meet in a point) and dually
\begin{equation}
\label{P0i_def}
P^0_i:=\spear^{0n}_{ij} \wedge \spear^{0j}_{in}.
\end{equation}
To see that neither of these vanish, assume the contrary, that $\axis^{0n}_{ij}= \axis^{0j}_{in}$ (or its dual).  Then by definition \eqref{axis0nij_def}, $\axis^{0n}_{ij}$ would be a shared line common to all four of the planes $p_i^n$, $p_j^n$, $p_n^j$ and $p_i^j$ which would then at the same time coincide with $\axis^{ij}_{0n}$, $\axis^{jn}_{i0}$ and $\axis^{in}_{0j}$, contradicting the assumption that the latter are not all identical.

The possibility still remains that  these final 3 axes or 3 spears are not linearly independent.
Define the four symbols $\{++,+-,-+,--\}$ as follows: $+$ in the first place means the 3 axes $\axis^{0n}_{ij}$ do not have a common plane, while $-$ means they do; while $+$ in the second place means that associated 3 spears do not have a common point, while $-$ means they do.


We first  dispose of case $(--)$.
Here, the shared plane of the $\axis^{0n}_{ij}$, doubled, is the final point-wise quadric and the shared common point of the $\spear^{0n}_{ij}$, doubled, is the final quadric plane-wise.  We now show that this point and plane are incident.

Claim: In case $(--)$, the common point $P=P_1^0=P_2^0=P_3^0$ of the $\axis^{ij}_{0n}$ is incident to the common plane $p=p_1^0=p_2^0=p_3^0$ of the $\spear^{ij}_{0n}$, therefore $p^2$ (and plane-wise, $P^2$) is in $X$ contact with $\ptq{S}_1, \ptq{S}_2$ and $\ptq{S}_3$.
\begin{proof}[Proof of claim]
This would be shown if p were tangent to any of $\ptq{S}_1, \ptq{S}_2$ and $\ptq{S}_3$ because $P$ is polar to $p$ for all 3 of them and so if $p$ were tangent to $\ptq{S}_j$, $P$ would be the point of tangency and thus incident to $p$.  Therefore, for the sake of contradiction, we assume that $p$ is not incident to $P$ which implies that $\ptq{R}_j=p \cap \ptq{S}_j$ is a nondegenerate conic.

Then $\ptq{R}_i$ and $\ptq{R}_j$ are in double contact with $\axis^{0n}_{ij}$ as chord of contact because
\begin{equation} \label{ringintersections}
\axis^{0n}_{ij} \cap \ptq{R}_i=\axis^{0n}_{ij} \cap \ptq{S}_i
=p^n_j \cap p^n_i \cap \ptq{S}_i
=p^n_j \cap p^n_i \cap \ptq{T}^n
=\axis^{0n}_{ij} \cap \ptq{T}^n
=p^n_i \cap p^n_j \cap \ptq{S}_j
=\axis^{0n}_{ij} \cap \ptq{R}_j
\end{equation}
The first equality follows from the fact that the intersection points of $\ptq{S}_i$ with $\axis^{0n}_{ij}$ necessarily occur on $\ptq{R}_i$ since as point sets $\axis^{0n}_{ij} \subset p^0_i$, the second equality follows from the definition \eqref{axis0nij_def} of $\axis^{0n}_{ij}$, the third from the ring contact of $\ptq{S}_i$ and $\ptq{T}^n$ in the plane $p^n_i$ and the rest of the equalities follow by symmetrical reasoning.  Both $\ptq{S}_i$ and $\ptq{S}_j$ are in contact with $\ptq{T}^n$ at these intersection points, so slicing the common tangent planes with $p$ yields the lines of double contact of $\ptq{R}_i$ and $\ptq{R}_j$.  Therefore, under the assumption that $P$ and $p$ are not incident, we obtain the conclusion that the three conics $\ptq{R}_i$ are mutually in double contact with chords of contact $\axis^{0n}_{ij}$, however these are concurrent in the point $O$, so by Lemma \ref{noconcdoublecontacttraingles} must be in the same double contact pencil which would imply that the 3 axes $\axis^{0n}_{ij}$ coincide which we've already seen cannot be the case.  Therefore $P$ is incident to $p$.
\end{proof}


The remainder of the proof assumes that the 3 axes are linearly independent (that is, the cases $(++)$ and $(+-)$. Dualizing the arguments then takes care of the remaining case $(-+)$.)  

When the 3 axes are linearly independent, the rings of contact are 
\begin{equation}
\ptq{R}^0_i:=\ptq{S}_i \cap p^0_i.
\end{equation}
The conics $\ptq{R}^0_i$ and $\ptq{R}^0_j$ either meet in two points on $\axis^{0n}_{ij}$ or $\axis^{0n}_{ij}$ is a shared tangent line, because $\axis^{0n}_{ij}$ is the intersection line of their respective planes and
\begin{equation*}
\axis^{0n}_{ij} \cap \ptq{R}^0_i=\axis^{0n}_{ij} \cap \ptq{R}^0_j
\end{equation*}
by the same reasoning as in Eq. \eqref{ringintersections}.

To complete the cube diagram and prove the theorem, we apply Lemma \ref{lem:conical-octahedron} to obtain a unique quadric containing all three conics that is in ring contact with $\ptq{S}_i$ along $\ptq{R}^0_i$ with ring point $P^0_i$ for $i=1,2,3$.   

Using this fact, we now provide a synthetic proof of the 8-quadric theorem, followed by a self-contained coordinate-based proof which includes a proof of Lemma \eqref{lem:conical-octahedron} as part of the construction.

Recall that to complete the proof we must show that for each $j \in \{1,2,3\}$,  (i) $\ptq{T}^{0}$ is in ring contact to $\ptq{S}_{j}$; and (ii) all 4 ring planes for the face $\face^{0j}_{km}$ have a common axis.   

Condition (ii) is trivially fulfilled since we have used the face axes $\axis^{0j}_{km}$ to construct the new ring planes.

We now show that condition (i) is also fulfilled. The polarity of $\ptq{R}^0_i$ in the plane $p^0_i$ maps a point $P$ to the line of intersection of $p^0_i$ and the polar plane of $P$ with respect to any quadric containing $\ptq{R}^0_i$: $\ptq{S}_i$ and $\ptq{T}^0$ in particular.  By equations \eqref{axis0nij_def} and \eqref{spear0nij_def}, and Lemma \ref{rc_polarity}, the polarities on $\ptq{S}_j$, $\ptq{S}_k$ and $\ptq{T}^0$ all behave the same way when restricted to points on $\axis^{0n}_{jk}$ and planes through $\spear^{0n}_{jk}$.  Combining these two facts and applying them to the cases $i=j$ and $i=k$, any point on $\axis^{0n}_{jk}$ is mapped under the respective polarities of $\ptq{R}^0_j$ and $\ptq{R}^0_k$ to two lines which together span the polar plane of $P$ with respect to any of the quadrics $\ptq{S}_j$, $\ptq{S}_k$, $\ptq{T}^0$.  Finally, considering the same scenario with $\axis^{0k}_{jn}$ in place of $\axis^{0n}_{jk}$, and appealing to the definitions of $P^0_j$ and $p^0_j$ in equations \eqref{p0i_def} and \eqref{P0i_def}, we conclude that the polarities on $\ptq{S}_j$ and $\ptq{T}^0$ agree when restricted to planes through $P^0_j$ and points on $p^0_j$ so $\ptq{T}^0$ and $\ptq{S}_j$ are in ring contact by Lemma \ref{rc_polarity}.
\end{proof}

Readers who are content with the above proof may skip ahead to Section \ref{sec:examples_of_8_quadric_configurations}. 
\begin{proof}[Coordinate-based proof of 8-quadric theorem]
Here we construct the point-wise matrix of the eighth quadric in terms of the point-wise matrices of the seven given ones. Dualizing yields the plane-wise matrix.

Let the matrices representing the quadrics $\ptq{S}_i$ and $\ptq{T}^m$ in the ordered basis $(O,X_1,X_2,X_3)$ where
\begin{equation}
X_n=\axis^{0n}_{ij} \wedge o.
\end{equation}
That this forms a basis follows from the independence of the $\axis^{0n}_{ij}$ since we are working in case (I), together with the fact that $O$ is not incident to $o$.  Because $o$ is polar to $O$ for all of the given quadrics, and $X_1,X_2$ and $X_3$ span $o$, their representing matrices in this basis may be chosen to have the form
\begin{equation}
\label{polarizing_matrix}
\begin{pmatrix}
1 & 0 & 0 & 0\\
0 & \ast & \ast & \ast\\
0 & \ast & \ast & \ast\\
0 & \ast & \ast & \ast
\end{pmatrix}.
\end{equation}
Let $\ptq{S}_1$ be represented by
\begin{equation}
A=\begin{pmatrix}
1 & 0 & 0 & 0\\
0 & a_{11} & a_{12} & a_{13}\\
0 & a_{12} & a_{22} & a_{23}\\
0 & a_{13} & a_{23} & a_{33}
\end{pmatrix},
\end{equation}
$\ptq{S}_2$ by
\begin{equation}
B=\begin{pmatrix}
1 & 0 & 0 & 0\\
0 & b_{11} & b_{12} & b_{13}\\
0 & b_{12} & b_{22} & b_{23}\\
0 & b_{13} & b_{23} & b_{33}
\end{pmatrix},
\end{equation}
and $\ptq{S}_3$ by
\begin{equation}
C=\begin{pmatrix}
1 & 0 & 0 & 0\\
0 & c_{11} & c_{12} & c_{13}\\
0 & c_{12} & c_{22} & c_{23}\\
0 & c_{13} & c_{23} & c_{33}
\end{pmatrix}.
\end{equation}
Then $\ptq{R}^0_1$ is in the plane $p^0_1$ spanned by $O, X_2$ and $X_3$ and is the restriction of $\ptq{S}_1$ thereto and thus can be represented in the ordered basis $(O,X_2,X_3)$ by the matrix
\begin{equation}
\bar A=\begin{pmatrix}
1 & 0 & 0 \\
0 &  a_{22} & a_{23}\\
0 &  a_{23} & a_{33}
\end{pmatrix}.
\end{equation}
The conics $\ptq{R}^0_2$ and $\ptq{R}^0_3$ can similarly be represented in the bases $(O,X_1,X_3)$ and $(O,X_1,X_2)$ by
\begin{equation}
\bar B=\begin{pmatrix}
1 & 0 & 0 \\
0 &  b_{11} & b_{13}\\
0 &  b_{13} & b_{33}
\end{pmatrix} \qquad \mathrm{and} \qquad
\bar C=\begin{pmatrix}
1 & 0 & 0 \\
0 &  c_{11} & c_{12}\\
0 &  c_{12} & c_{22}
\end{pmatrix}
\end{equation}
respectively.  If $\ptq{R}^0_2$ and $\ptq{R}^0_3$ are tangent then $b_{11}=c_{11}=0$, otherwise
\begin{equation}
\begin{pmatrix}
1 & 0  \\
0 &  c_{11} 
\end{pmatrix}
=
\begin{pmatrix}
1 & 0  \\
0 &  b_{11} 
\end{pmatrix}
\end{equation}
represents the pair of intersection points in the basis $(O,X_1)$ for $\axis^{01}_{23}$.  This implies in particular that
\begin{equation}
c_{11}=b_{11}
\end{equation}
and by similar reasoning
\begin{equation}
a_{22}=c_{22} \qquad \mathrm{and} \qquad a_{33}=b_{33}.
\end{equation}
Therefore
\begin{equation}
\label{t0matrix}
D=\begin{pmatrix}
1 & 0 & 0 & 0\\
0 & d_{11} & c_{12} & b_{13}\\
0 & c_{12} & d_{22} & a_{23}\\
0 & b_{13} & a_{23} & d_{33}
\end{pmatrix},
\end{equation}
where
\begin{equation}
d_{11}:=c_{11}=b_{11} , \qquad d_{22}:=a_{22}=c_{22} \qquad \mathrm{and} \qquad d_{33}:=a_{33}=b_{33},
\end{equation}
represents a quadric $T_0$ containing the three conics $\ptq{R}^0_1, \ptq{R}^0_2$ and $\ptq{R}^0_3$.  To show that $\ptq{T}^0$ satisfies the required ring contacts as well as conditions (ii) and (ii)', it is sufficient to show that it has $\axis^{0n}_{ij}$ polar to $\spear^{0n}_{ij}$ for all $\{i,j,n\}=\{1,2,3\}$ because then the intersection
\begin{equation}
P^0_i:=\spear^{0n}_{ij} \wedge \spear^{0j}_{in}
\end{equation}
is polar to $p^0_i$ with respect to both $\ptq{S}_i$ and $\ptq{T}^0$ and is thus the ring point corresponding to the ring of contact $\ptq{R}^0_i$.

Because $\ptq{S}_1$ and $\ptq{S}_2$, by virtue of ring contact with $\ptq{T}^3$, both have $\axis^{03}_{12}$ polar to $\spear^{03}_{12}$, both have $X_3$ polar to $x_3$ where
\begin{equation}\label{littlexdef}
x_n=\spear^{0n}_{ij} \vee O.
\end{equation}
Therefore, because the fourth column of each matrix gives the plane coordinates of the image of $X_3$:
\begin{equation}
a_{13}:a_{23}:a_{33}=b_{13}:b_{23}:b_{33}.
\end{equation}
By similar reasoning, 
\begin{equation}
a_{12}:a_{22}:a_{23}=c_{12}:c_{22}:c_{23}
\end{equation}
and
\begin{equation}
b_{11}:b_{12}:b_{13}=c_{11}:c_{12}:c_{13}.
\end{equation}
This implies that
\begin{equation}
b_{13}:a_{23}:d_{33}=b_{13}:b_{23}:b_{33},
\end{equation}
\begin{equation}
a_{12}:a_{22}:a_{23}=c_{12}:d_{22}:a_{23}
\end{equation}
and
\begin{equation}
d_{11}:c_{12}:b_{13}=c_{11}:c_{12}:c_{13}.
\end{equation}
Therefore, $\ptq{T}^0$ has $X_n$ polar to $x_n$ for $n=1,2,3$ completing the proof.
\end{proof}

\subsubsection{When the eighth quadric is not regular}
The case $(++)$, in which both triples of axes and spears are linearly independent, clearly produces regular matrices, hence the eighth quadric is regular. This section considers what kinds of degenerate quadrics arise in the other cases. We note that by Eq. \eqref{littlexdef} the $\spear^{0n}_{ij}$ meet at a point exactly when the $x_n$ share a common line, i.e., when the resulting matrix $D$ is degenerate.  It must also have rank $3$ because the $\spear^{0n}_{ij}$ do not coincide and $o$ is not incident to $O$ and thus is independent of the $x_n$.  Therefore, in case $(+-)$, if the $\spear^{0n}_{ij}$ are dependent, $\ptq{T}^0$ is a cone. In case $(-+)$, we dualize and obtain a conic.  The final point-wise quadric in ring contact is the plane of this conic, doubled.  In case $(--)$, the final point-wise quadric is again a double plane but now it does not arise as the plane of a conic but rather as a contact element together with a pair of lines: the lines of $X$-contact (note however that the $\ptq{S}_j$ are not in $X$-contact with each other: as discussed in, e.g. \cite{stachel2020}, there are several different contact relations for which two quadrics may meet in a pair of lines, $X$-contact being the special case where this is the entirety of the intersection; in the present case, this will hold for each $\ptq{S}_j$ with the double plane $\ptq{T}^0$ but not for the $\ptq{S}_j$ with each other).  From this it is clear that the full information of the quadrics is not contained in the point-wise or plane-wise description alone but requires knowledge of both (as well as possibly lines).  However, in the present paper we leave the precise formulation somewhat one-sided, leaving a more principled exposition to future work.

\subsubsection{Refinement of Penrose 8-quadric theorem}
The above discussion contains a proof of the following refinement of Thm. \ref{thm:8quadric}.
\begin{theorem} \label{thm:refinement}
Given a regular Penrose 7-quadric configuration. Let $(\ptq{Q}, \plq{Q})$ be the quadric that, by Thm. \ref{thm:8quadric}, completes the configuration. Consider the 3  face axes and 3 face spears of the cube faces meeting at the eighth vertex. 
Then:
    \begin{enumerate}
    \item If the axes do not lie in a plane and the spears do not pass through a point,  $\ptq{Q}$ and $\plq{Q}$ are regular.
    \item If the axes lie in a plane and the spears do not pass through a point, $\ptq{Q}$ is a  double plane and $\plq{Q}$ is a conic,
    \item If the axes do not lie in a plane and the spears pass through a point,  $\ptq{Q}$ is a cone and $\plq{Q}$ is a double point,
    \item If the axes lie in a plane and the spears pass through a point, $\ptq{Q}$ is a double plane and $\plq{Q}$ is an incident double point.
    \end{enumerate}
\end{theorem}

There is an analogous version of this theorem for a Penrose 7-conic configuration. Out of space considerations, we omit the details.

\subsubsection{Proof of Extrusion Lemma}
\label{sec:extrusion}
We now turn to the promised proof of Lemma \ref{lemma:extrusion}: \emph{Given a regular Penrose 7-conic  configuration, it is possible to extrude it to a regular Penrose 7-quadric configuration.}

\begin{proof}[Proof of extrusion lemma]

Let $(P_1,P_2,P_3)$ be an ordered basis of the plane containing a given conic configuration and $O$ a point not on this plane.  Then if we choose a plane $o$ with homogeneous coordinates \mbox{$[u_0:u_1:u_2:u_3]$} in the basis dual to $(O,P_1,P_2,P_3)$, any conic in the plane, represented by a matrix
\begin{equation}
A=\begin{pmatrix}
a_{11} & a_{12} & a_{13}\\
a_{12} & a_{22} & a_{23}\\
a_{13} & a_{23} & a_{33}
\end{pmatrix}
\end{equation}
in the basis $(P_1,P_2,P_3)$ can be extruded to a quadric
\begin{equation}
\bar A=\begin{pmatrix}
u_0 & u_1 & u_2 & u_3\\
u_1 & a_{11} & a_{12} & a_{13}\\
u_2 & a_{12} & a_{22} & a_{23}\\
u_3 & a_{13} & a_{23} & a_{33}
\end{pmatrix}
\end{equation}
going through the conic and having $O$ polar to $o$.

The quadric configuration resulting from this construction applied to the 7 conics of the hypothesis of the 8-conic theorem will automatically satisfy conditions (i) and (ii) 
of the 8-quadric theorem, provided that the chords of contact of condition (i) for the conic case are extruded to ring planes through $O$. That is, a regular Penrose 7-conic configuration produces a regular Penrose 7-quadric configuration.

Let $A$ and $B$ represent two conics in double contact.  Then
\begin{equation}
\mathrm{rank}(A+rB)=1
\end{equation}
for some scalar $r$ and $\bar A$ will be in ring contact with
\begin{equation}
\label{extrusion}
\bar B(t)=\begin{pmatrix}
u_0 & u_1 & u_2 & u_3\\
u_1 & tb_{11} & tb_{12} & tb_{13}\\
u_2 & tb_{12} & tb_{22} & tb_{23}\\
u_3 & tb_{13} & tb_{23} & tb_{33}
\end{pmatrix}
\end{equation}
with ring plane going through $O$ exactly when $t=r$.  Thus in order to perform the extrusion, it suffices to simultaneously scale all the matrices representing the 7 conics so that for any pair in double contact represented by the matrices $A$ and $B$ we have
\begin{equation}\label{double_contact_scaling}
\mathrm{rank}(A-B)=1.
\end{equation}

To begin with, the scaling of the matrix representing $\ptq{S}_0$ can be chosen arbitrarily and then the matrices representing $\ptq{T}^1, \ptq{T}^2$, and $\ptq{T}^3$ are uniquely defined by Eq. \eqref{double_contact_scaling}.  Then it only remains to scale $\ptq{S}_1, \ptq{S}_2$ and $\ptq{S}_3$ so that Eq. \eqref{double_contact_scaling} holds for both double contacts of each.  Let $A,B$ and $C$ be the appropriately scaled matrices representing $\ptq{S}_0, \ptq{T}^1$ and $\ptq{T}^2$ respectively:
\begin{equation} \label{dc1}
A-B=V^1_0, \qquad A-C=V^2_0
\end{equation}
where $V^j_i$ is a rank 1 matrix representing the doubled chord of contact.  Then let $D$ be the matrix representing $\ptq{S}_3$ scaled so that
\begin{equation}
D-B=V^1_3.
\end{equation}
Then all four matrices are scaled so the double contact of $\ptq{T}^2$ and $\ptq{S}_3$ in these scalings must be written
\begin{equation} \label{dc2}
D-rC=V^2_3
\end{equation}
for some determined scalar $r$.  However, by sequentially eliminating $A, B$ and $D$ from equations \eqref{dc1}-\eqref{dc2}, we obtain
\begin{equation}
(r-1)C=V^1_3+V^2_0-V^1_0-V^2_3.
\end{equation}
Because of condition (ii) of the hypothesis of the 8-conic theorem, all the chords of contact for the face $\mathcal{F}^{12}_{03}$ share a point whose representative in the basis $(P_1,P_2,P_3)$ is necessarily in the common kernel of $V^1_3$, $V^2_0$ ,$V^1_0$ and $V^2_3$. However, by the regularity assumption on $\ptq{T}^2$, $C$ has trivial kernel which is only consistent with Eq. \eqref{double_contact_scaling} in the case $r=1$. 

Similar reasoning applied to the other two faces yields the scaling needed to establish the desired result.

\end{proof}

\subsubsection{Proof of the 8-quadric theorem with a common contact element}
\label{Ooincidnetproof}

We can now give the promised proof of the 8-quadric theorem when there is a common contact element of all the quadrics.

\begin{proof}[Proof with common contact element]
Given a configuration of conics in double contact, Eq. \eqref{extrusion} exhausts all possible extensions of $B$ to a quadric polarizing $O$ and $o$ and having all these double contacts as slices of ring contacts, with $O$ as the point of concurrency of the contact axes and thus $o$ as the common plane of the ring points.  Therefore, to prove the case of the 8-quadric theorem where $o$ is incident to $O$, we apply the 8-conic theorem to a planar slice of the 7-quadric configuration and then extrude the resulting complete 8-conic via Eq. \eqref{extrusion}, choosing the parameter $t$ that extrudes the seven conics obtained from the slicing process back to the seven quadrics with which we started.  The eighth extruded quadric is then the desired conclusion of the 8-quadric theorem.
\end{proof}

\section{Examples of 8-quadric configurations}
\label{sec:examples_of_8_quadric_configurations}

The aim of this section is to illustrate the structure of the 8-quadric theorem by constructing several examples step by step and, by slicing them, to obtain all the famous 2D-porisms mentioned in the Introduction. The geometric proof in the next section provides a direct proof valid also for these examples. 

\begin{example}\label{ex:3-cones}
Consult Figure \ref{fig:3-cones}. We start with a sphere (green) as quadric $\ptq{S}_0$. Then we choose 3 planes $p_0^i$ intersecting the sphere in 3 (blue) rings $\ptq{R}_0^i$ such that each pair of rings meet in a pair of points (white). They form a conical octahedron, featured  above in Lemma \ref{lem:conical-octahedron}. The corresponding cones touching the sphere along those rings are taken as quadrics $\ptq{T}^i$ (yellow). Now for each pair of cones the intersection $\ptq{T}^i\cap \ptq{T}^j$ consists of two conics (red). These six conics meet three at a time in the 8 yellow points, the common intersection $\ptq{T}^1 \cap \ptq{T}^2 \cap \ptq{T}^3$ of all the three cones. We call this configuration a \emph{conical cube}, whereby the yellow points represent the vertices and each red conic represents a pair of opposite edges of the cube. 

We will choose one conic from each of the three pairs as a degenerate quadric $\ptq{S}_k, \{i,j,k\}=\{1,2,3\}$. This situation can be pictured as a limit of thin ellipsoids in ring contact with $\ptq{T}^i, \ptq{T}^j$ whose planes of ring contact $p_k^i, p_k^j$ finally converge to a single plane, which contains the chosen conic of the intersection.  

\begin{figure}[h!]
     \centering
     \includegraphics[width=0.5\textwidth]{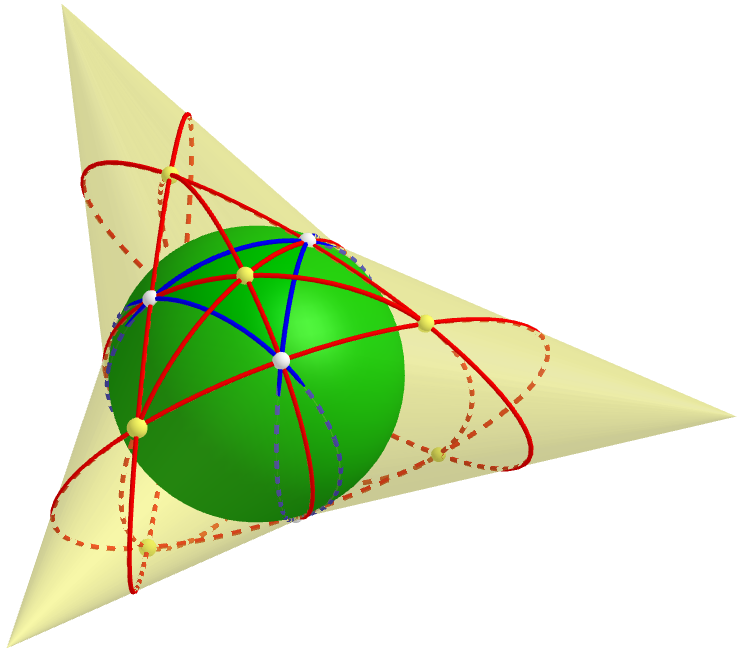}
     \includegraphics[width=0.45\textwidth]{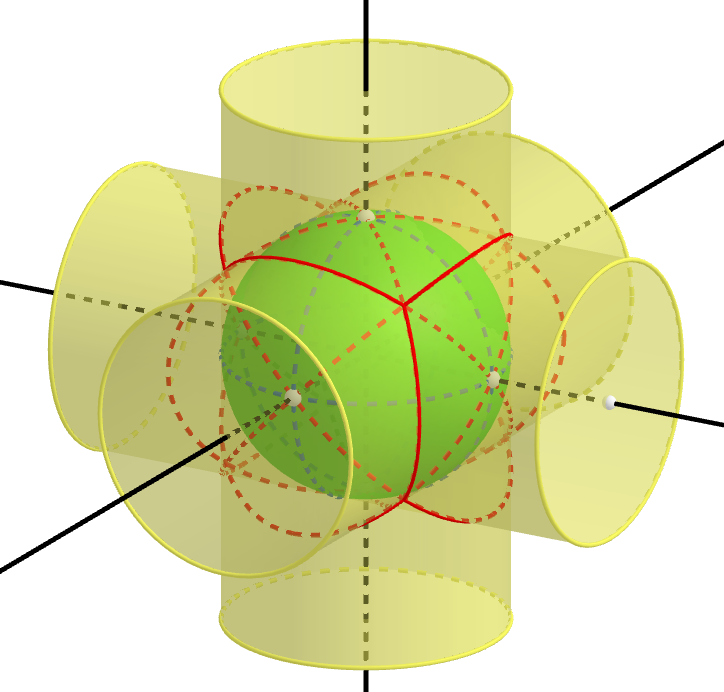}
     \caption{Three cones in ring contact with a sphere (left). The prototype of this situation are 3 cylinders with pairwise orthogonal axes (right). The 3 rings (blue) form a conical octahedron, whereas the 6 conics of the intersections (red) form a conical cube. }
     \label{fig:3-cones}
 \end{figure}

There are now two possibilities for such a triple of conics $\{\ptq{S}_1,\ptq{S}_2,\ptq{S}_3\}$, one leading to a regular $T_0$ and one not: Either the 3 conics meet in a point pair, 
or they form a conical octahedron. In Figure \ref{fig:3-cones}, left, the first possibility is shown by the three red conics that meet in the yellow point in the middle surrounded by the three blue ring conics (with a second meeting point on the opposite side of the sphere), while the second possibility is shown by the remaining three red conics, that surround the point in the middle and extend farther out along the 3 yellow cones and meet pairwise in two points.  

In the latter case the final quadric $\ptq{T}^0$ is regular (directly given by Lemma \ref{lem:conical-octahedron}). Since each quadric $\ptq{S}_i$ degenerates to a conic, all 3 rings $\ptq{R}_i^k$  on it also collapse to this conic and each axis $\axis_{ij}^{0k}$ is the intersection of two ring planes. In the former case, on the other hand, the final quadric $\ptq{T}^0$ is a point pair.  This case, however, cannot as easily be subsumed into the logic of section \ref{sec:proof8quadric} since the final 3 axes are identical: $\axis^{01}_{23}=\axis^{02}_{13}=\axis^{03}_{12}$. To picture this in terms of a limit, we consider a sequence of shrinking ellipsoids, the $\ptq{S}_i$ getting flatter and $\ptq{T}^0$ getting thinner and finally collapsing to a line segment connecting the two points, all planes through either point ultimately being tangent planes enveloping the ellipsoid in the limit.

Although these configurations are quite elementary in 3D their slices with a plane yield nontrivial theorems in 2D. In particular, consider the plane $p$ through the ring points of the 3 cones: It intersects the cones in 6 lines tangential to a circle being the intersection of the sphere $\ptq{S}_0$ with this plane. 
In the case of three concurrent conics $\ptq{S}_1,\ptq{S}_2,\ptq{S}_3$, they intersect the plane $p$ in 6 points which form a hexagon circumscribing the circle of intersection with the sphere. This leads to a Brianchon configuration.
In the case of three non-concurrent conics $\ptq{S}_1,\ptq{S}_2,\ptq{S}_3$ they intersect the plane $p$ in 6 points which split up into two triangles circumscribing the circle of intersection. This leads to a Poncelet configuration. 
A slice with a general plane gives the double contact theorem and the 4-conic theorem from \cite{evelyn1974seven}.
Note that a single 3D-configuration contains Brianchon figures in 4 different ways such that the Brianchon diagonals form a complete 4-point.

\begin{figure}[h!]
     \centering
     \includegraphics[width=0.4\textwidth]{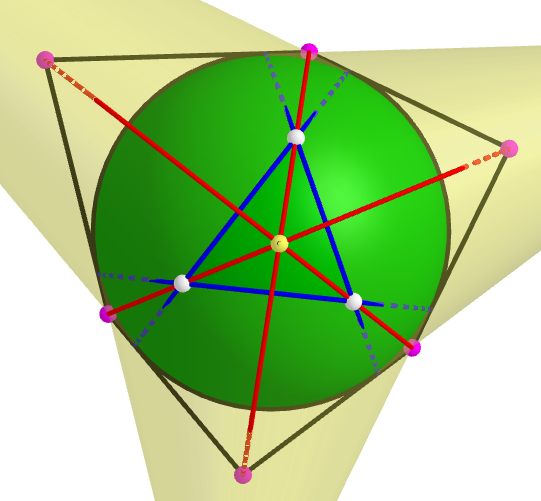}
     \hspace{0.5cm}
     \includegraphics[width=0.4\textwidth]{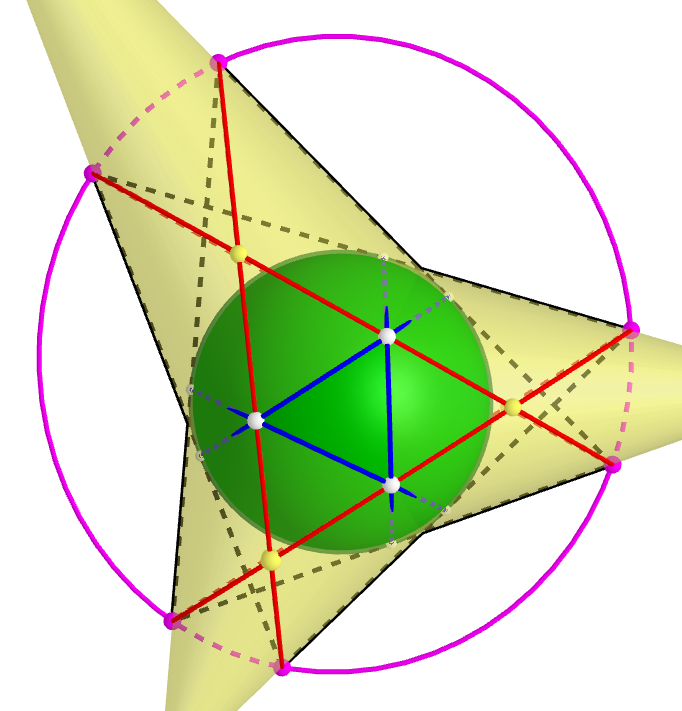}
     \caption{Left: Three concurrent conics of the conical cube induce a hexagon circumscribing a circle in the plane $p$ through the ring points of the cones and thus proves the Brianchon theorem. Right: In the case of  three non-concurrent conics the 6 points split into two triangles. This figure generates a Poncelet configuration. }
     \label{fig:Brianchon-Poncelet}
 \end{figure}
\end{example}

\begin{figure}[h!]
     \centering
     \includegraphics[width=0.55\textwidth]{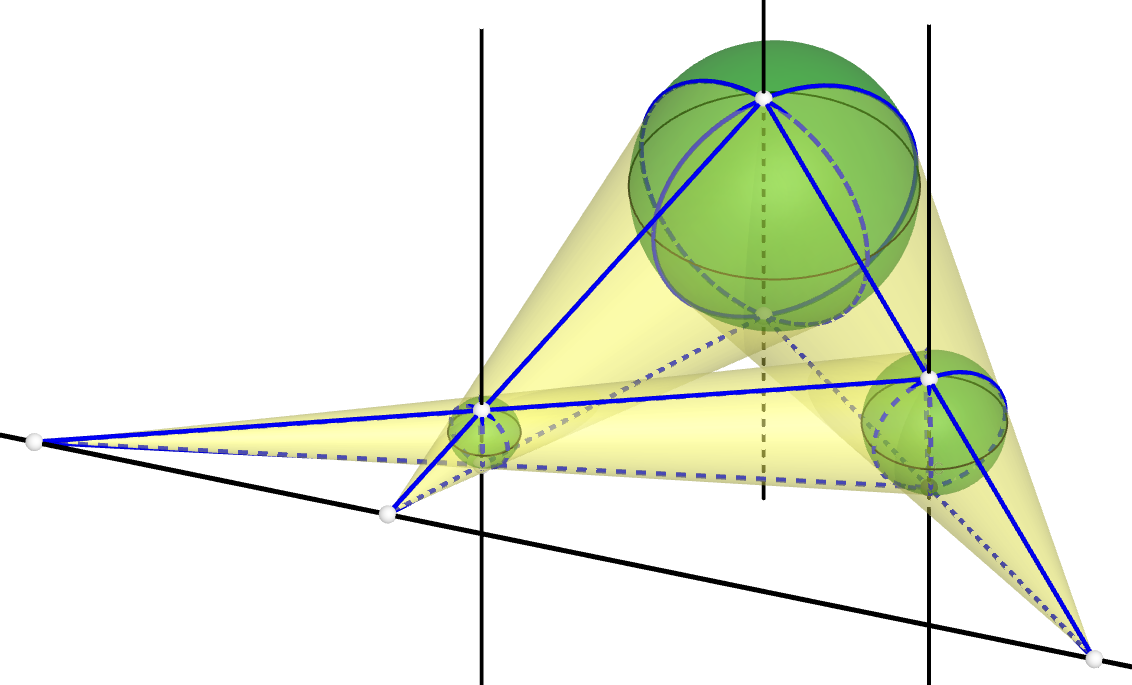}
    \includegraphics[width=0.4\textwidth]{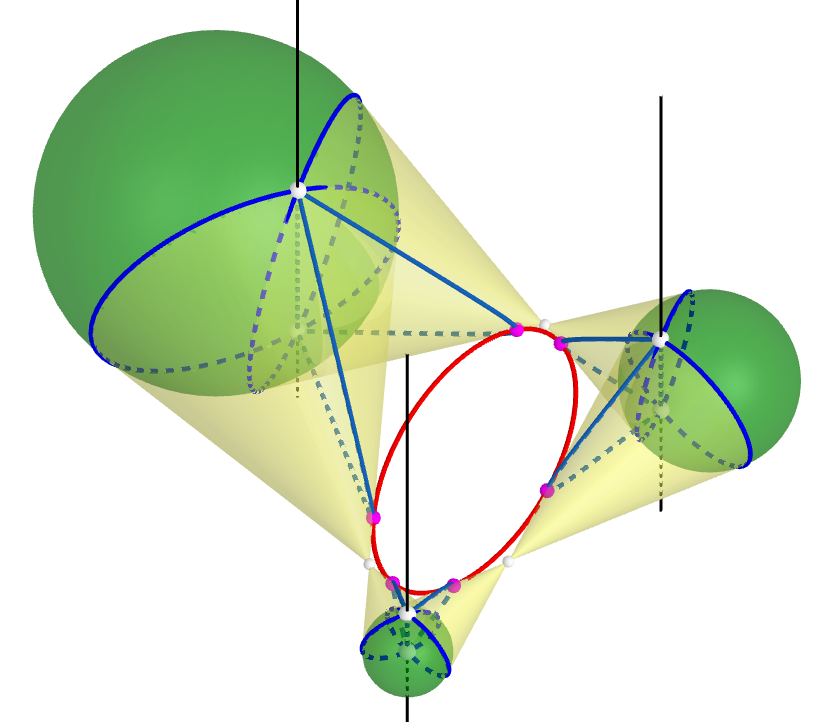}
     \caption{Left: Original proof of G. Monge of his theorem about the external tangents of three circles. Right: A theorem about the internal tangents}.
     \label{fig:Monge}
 \end{figure}

\begin{example}\label{ex:3-spheres}
In the next example, we consider three spheres $\ptq{T}^1,\ptq{T}^2,\ptq{T}^3$ such that no sphere is completely contained in another. All spheres intersect the plane at infinity in the same imaginary circle $\ptq{S}_0$, the absolute of Euclidean geometry. Thus, each sphere $\ptq{T}^i$ is in ring contact with $\ptq{S}_0$, where the conic $\ptq{S}_0$ as an imaginary degenerate quadric represents at the same time the rings $\ptq{R}_0^i$. We could proceed as in the previous example by choosing as quadrics $\ptq{S}_k$ the intersection of the spheres $\ptq{T}^i \cap \ptq{S}_0$, which are circles (besides the absolute imaginary circle at infinity) and by slicing with a plane through the centers of the 3 spheres we would obtain again a Brianchon-like theorem, which turns out to be the theorem of the radical center of 3 circles. 

Instead we choose cones enveloping each pair of spheres.  If the spheres are disjoint there are 6 such cones in total (dual to the 6 conics in Example \ref{ex:3-cones}) and again there are two cases for the choice of a triplet of cones. In one case the cones share 2 tangent planes, which are also common tangent planes of the 3 spheres.  See Figure \ref{fig:Monge}, left. Thus the final quadric $\ptq{T}^0$ is this pair of planes and the rings $\ptq{R}_i^0$ are 3 line pairs. This 3D configuration was already used by G. Monge (according to \cite{wells1992}{p.154}) in order to prove a theorem today bearing his name: The three pairs of external tangents of three circles meet in three points on a line.

In the other case the cones do not admit common tangent planes. See Figure \ref{fig:Monge}, right. Now the rings $\ptq{R}_i^0$ are obtained by intersecting each cone $\ptq{S}_i$ with the corresponding plane $p_i^0$ which is the join of a pair of vertical axes. The final quadric is again guaranteed by Lemma \ref{lem:conical-octahedron}. A planar slice through the centers of the spheres shows that 3 pairs of internal tangents of 3 circles are circumscribed to a conic. This theorem seems to be known only for the case, where the final conic degenerates into a point pair, see \cite{Baker_1936},\cite{salmon1917}.

Analogous to the Example \ref{ex:3-cones} there are 4 ways to assign a pair of common tangent planes to 3 spheres. The 6 ring points of the cones form a complete 4-line.

\end{example}

\begin{example}\label{ex:3-contact-elements}
Consider 3 planes $p_0^i$ tangent to a quadric $\ptq{S}_0$ at points $P_0^i$. In order to see the rings $\ptq{R}_0^i$ as real objects we choose here a regulus as quadric $\ptq{S}_0$, thus each ring is a line pair contained in the regulus. Consult Figure \ref{fig:Pascal}. This is the case referred to as X-contact in section \ref{quadrics_intro}, whereby the  ring point is incident to the ring plane. We want to picture these line pairs as `quadrics' $\ptq{T}^i$, although the situation is more complex.  They are in fact double planes when considered point-wise and double points when considered plane-wise, but in an X-contact pencil we can  identify them with the line pairs themselves. We can imagine these `quadrics' as the limit of a sequence of reguli all containing the line pair $\ptq{S}_0 \cap p_0^i$, as they get flatter and flatter. 

Each pair of line pairs $\ptq{T}^i,\ptq{T}^j$ meet in a pair of points on $\ptq{S}_0$ and are contained in a pair of tangent planes since $\ptq{T}^i,\ptq{T}^j \subset \ptq{S}_0$. We choose the plane pairs to be the quadrics $\ptq{S}_k$. Then the eighth quadric is the double plane spanned by the three points $P_0^1,P_0^2,P_0^3$ we started with, and the rings $\ptq{R}_i^0$ are the 3 sides of this triangle, each side considered as a double line, see Figure \ref{fig:Pascal}.

The configuration described above is well known. The 3 plane pairs $\ptq{S}_1,\ptq{S}_2,\ptq{S}_3$ form a hexahedron which is inscribed around a regulus such that 6 of the edges (a skew hexagon) lie in the regulus. Thus a planar slice of this configuration yields a hexagon inscribed in a conic such that opposite sides of the hexagon (corresponding to the plane pairs $\ptq{S}_i$) meet in 3 points which lie on a line (the intersection of the slicing plane with the plane $\ptq{T}^0$), the Pascal line of the hexagon. This 3D proof of Pascal's Theorem was already found in 1826 by G. Dandelin \cite{Dandelin1826}.  

\begin{figure}[h!]
     \centering
     \includegraphics[width=0.66\textwidth]{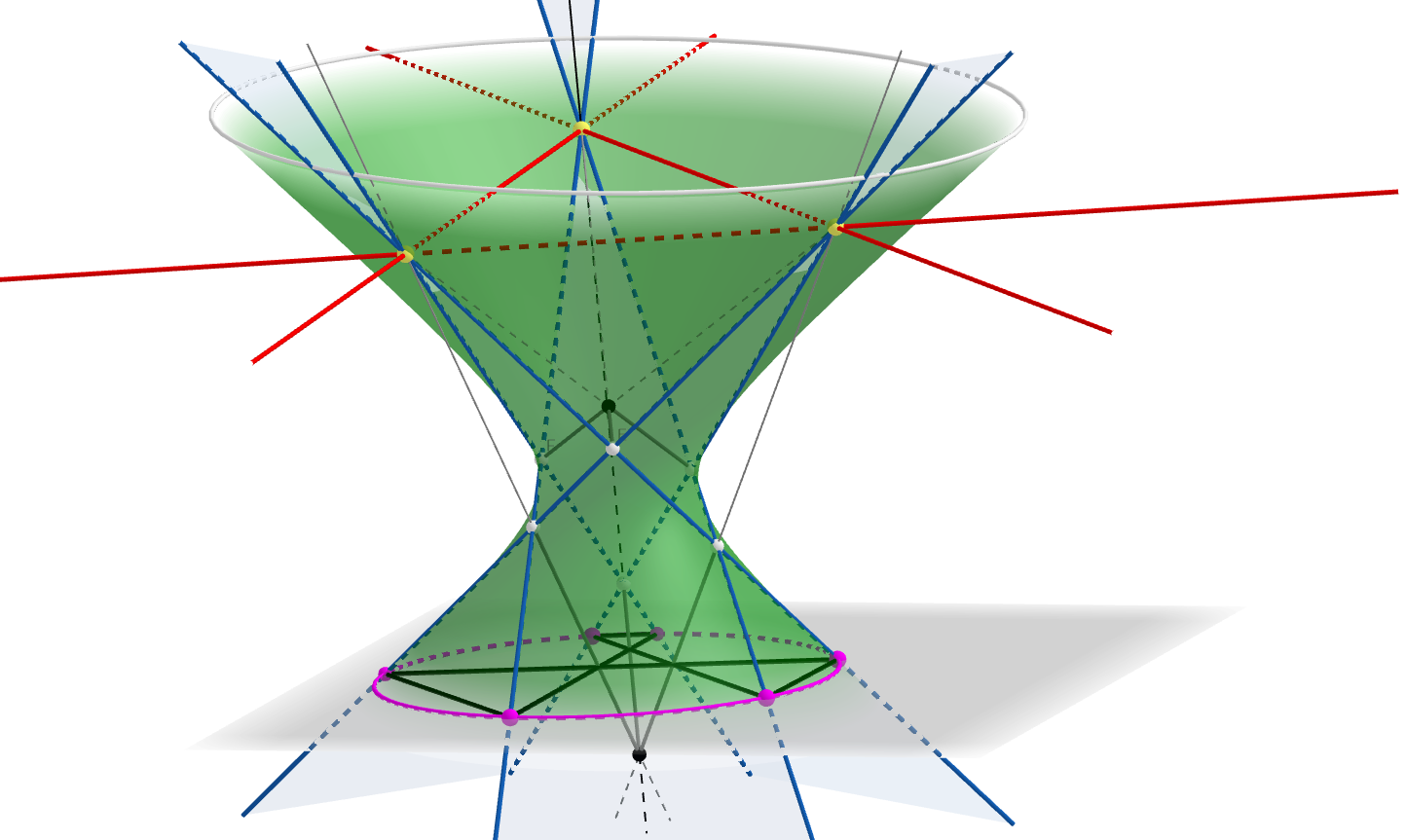}
     \caption{The proof of Pascal's Theorem using a hexahedron inscribed around a regulus. The Pascal axis is the intersection of the planar slice (bottom) with the ideal plane of the hexahedron (top).}
     \label{fig:Pascal}
 \end{figure}

Finally we will look at the same configuration but construct it in reverse order, starting with a double plane $\ptq{S}_0$ followed by three plane pairs $\ptq{T}^i$, each plane pair meeting in a line contained in $\ptq{S}_0$. Each pair of plane pairs intersects in two line pairs. In this special case the conical cube of the 6 line pairs are the opposite edges of a usual projective cube. As in the examples before a triplet of chosen line pairs might be concurrent or not. In the first case the 3 line pairs obviously split into two triples of lines, belonging to opposite vertices in the cube, see Figure \ref{fig:Desargues} left. A planar slice yields Desargues' theorem. For the other choice, the 3 line pairs form a closed zigzag and the planar slice yields Pascal's theorem. In both cases the axis in the theorem is the intersection with the plane $\ptq{S}_0$.

\begin{figure}[h!]
     \centering
     \includegraphics[width=0.25\textwidth]{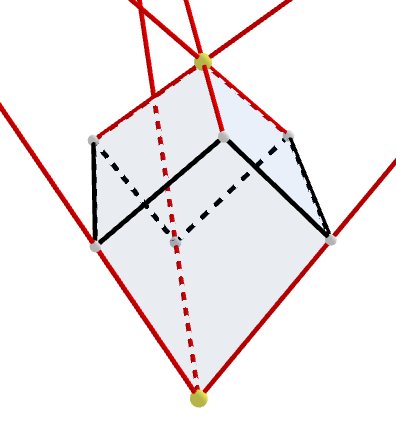}
     \includegraphics[width=0.3\textwidth]{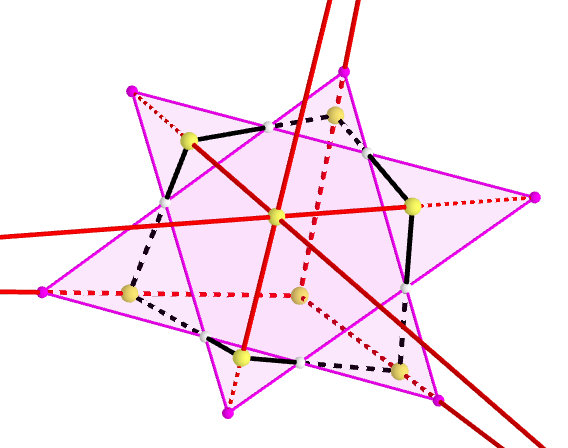}
     \includegraphics[width=0.3\textwidth]{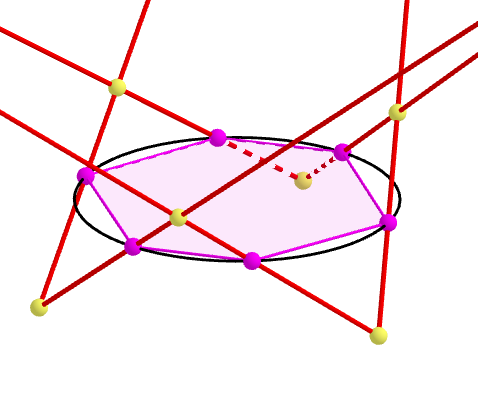}
     \caption{If three pairs of opposite edges on a hexahedron (a projective cube) are chosen, such that one line from each pair meet in two opposite vertices (left), a planar slice splits into two perspective triangles, proving Desargues' Theorem (middle). Otherwise the three edge pairs form a zigzag on the cube and yield a hexagon in a planar slice, proving Pascal's Theorem (right).}
     \label{fig:Desargues}
 \end{figure}
 \end{example}

\begin{figure}[ht!]
     \centering
     Ex.1 \includegraphics[width=0.24\textwidth]{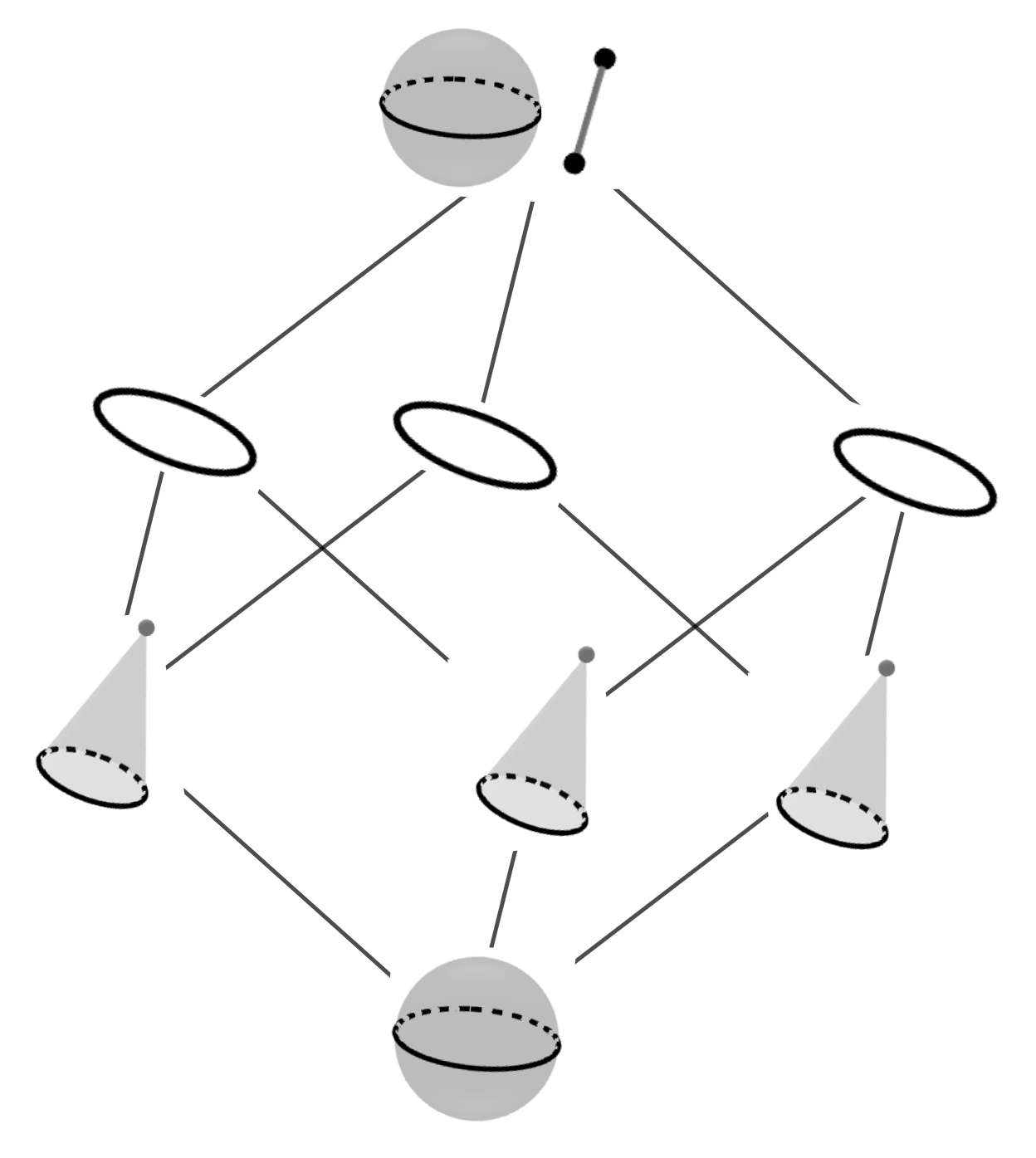}
     Ex.2 \includegraphics[width=0.25\textwidth]{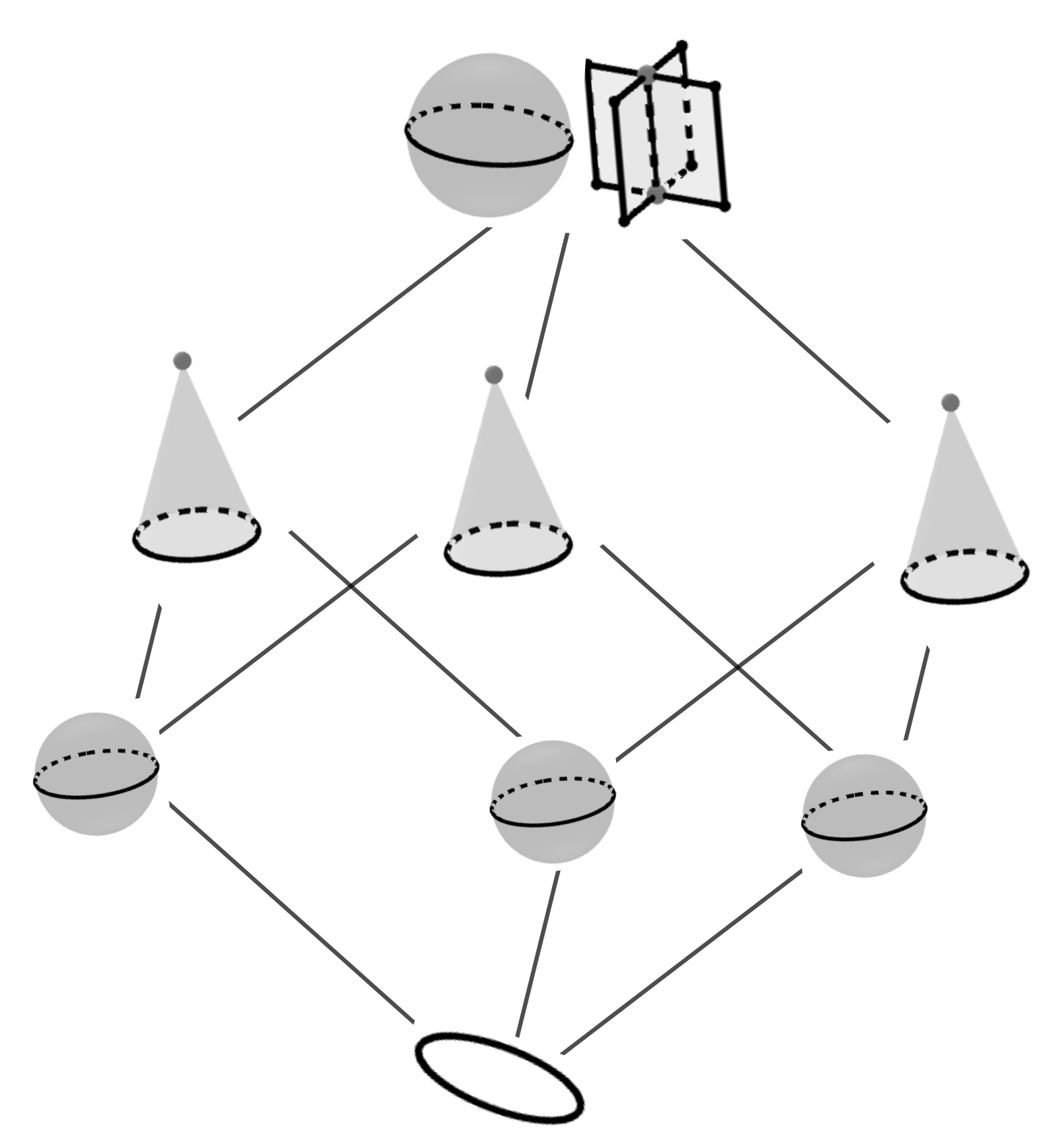}
     Ex.3 \includegraphics[width=0.25\textwidth]{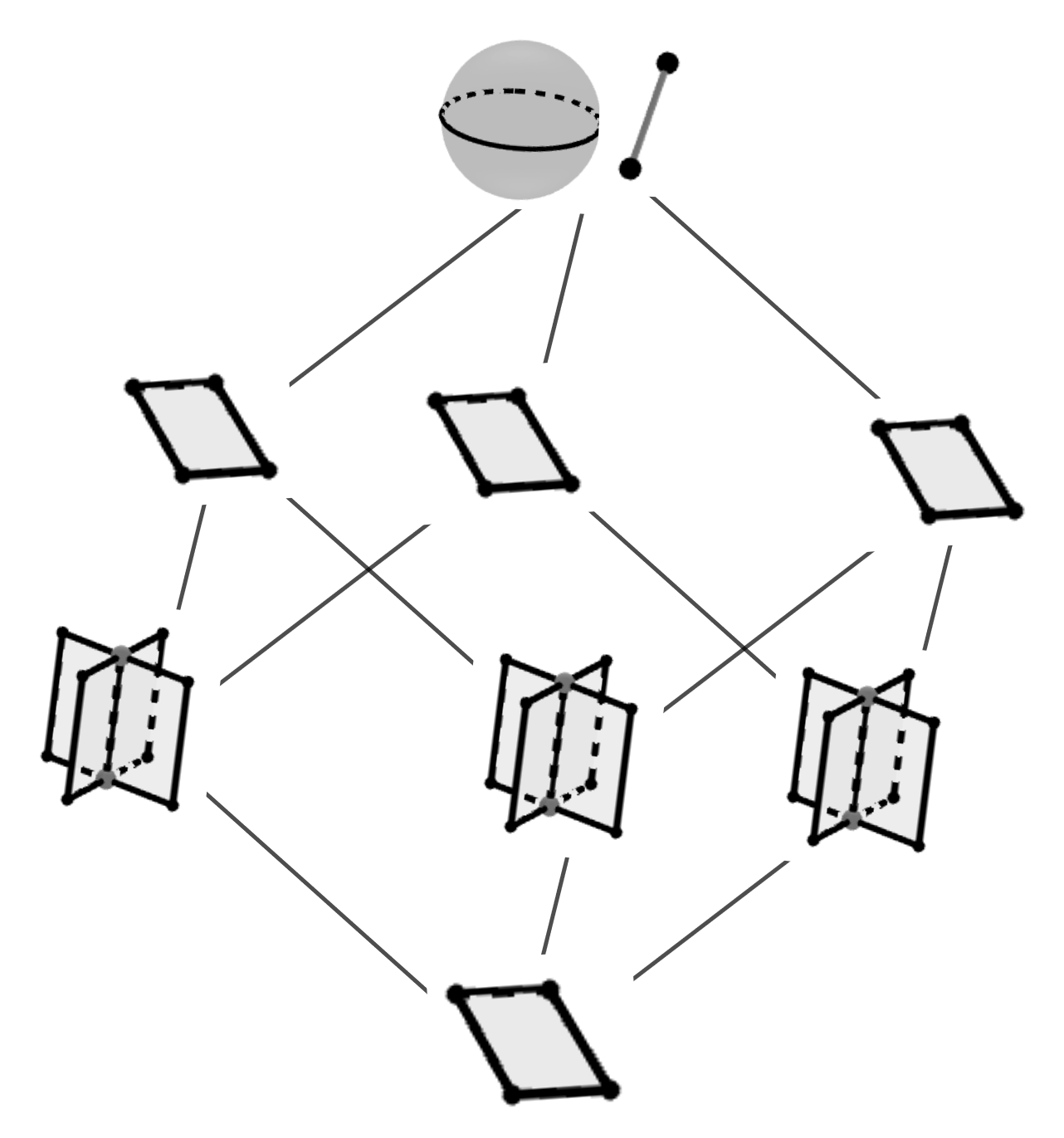} 
     \caption{Graphs of the 8-quadric configurations of the above examples. The corresponding 2D-theorems obtained by appropriate slices are: Ex.1 Poncelet/Brianchon, Ex.2 Original/Monge, Ex.3 Pascal/Desargues}
     \label{fig:cubegraphs}
 \end{figure}

The above examples show the following common structure: The three quadrics $\ptq{T}^1,\ptq{T}^2,\ptq{T}^3$ give rise to 6 degenerate quadrics either by pairwise intersections (Examples \ref{ex:3-cones} and \ref{ex:3-contact-elements}) or by pairwise envelopes (Example \ref{ex:3-spheres}). Depending on the choice of such a triplet $\ptq{S}_1,\ptq{S}_2,\ptq{S}_3$, one from each pair, the final quadric $\ptq{T}^0$ is either regular or a pair of points (or planes). This is summarized in Figure \ref{fig:cubegraphs}.

\charliex{

\section{Algebraic proof of the Penrose theorems}
\label{sec:alg-proof}

We now turn to a purely algebraic proof of Penrose's 8-conic theorem.
We begin with describing a function $\mathcal{P}(\ptc{S}_0, \myln{p}, \myln{q}, \myln{r}, a,b,c)$ whose value is  a complete Penrose 8-configuration $\{\ptc{S}_i, \ptc{T}^j\}$ whereby:
\begin{itemize}
    \item $\ptc{S}_0$ is an arbitrary complete conic, 
    \item $\{\myln{p},\myln{q},\myln{r}\}$ are three arbitrary weighted lines in $\RP{2}$, and
    \item $\{a,b,c\}$ are real parameters. 
\end{itemize}

Note that the parameter space of $\mathcal{P}$ is 17-dimensional:  $\ptc{S}_0$ depends on 5 parameters, each of the weighted lines $\{\myln{p},\myln{q},\myln{r}\}$ depends on 3 real parameters, and $\{a,b,c\}$ contributes 3 more real parameters for a total of 17.
 
The following discussion presents and verifies the explicit formulas provided by the algorithm for the cube $\{\ptc{S}_i, \ptc{T}^j\}$. 
It makes repeated use of the fact that two point-wise conics are in double contact with each other if and only if their pencil contains a double line. That is, there exists  $a$ and $b$ such that
\begin{equation} \label{eq:dblln}
    a\,\ptc{C}_1 + b\,\ptc{C}_2 = \abdb{m}
\end{equation}
where ${m}$ is the joining line of the contact points. This equation can obtained from Eq. \ref{dcp} by 
translating the latter into abbreviated notation (Section \ref{sec:abbr-nota}).  Note that $\abdb{m}$ has to be non-zero to establish double contact, so that we can normalize $a$ and $b$ to have $\abdb{m}$ alone on the right-hand side.

It focuses on the case of point-wise conics; similar formulas exist also for line-wise conics.

\subsection{The Penrose formulas}

We start with an arbitrary conic $\ptc{S}_0$ and a line $\myln{p}$. Consider the sum:
\begin{equation}\label{eq:t1}
     \ptc{T}^1 = \ptc{S}_0+ \abdb{p}
\end{equation}
Then $\ptc{T}^1$ and $\ptc{S}_0$ are in double contact because we have $\ptc{T}^1- \ptc{S}_0=\abdb{p}$ as in Eq. \ref{eq:dblln}.

In the same way, we form a conic $\ptc{T}^2$ and $\ptc{T}^3$ from $\ptc{S}_0$ and lines $\myln{q}$ and $\myln{r}$:
\begin{eqnarray*}
    \ptc{T}^2&=& \ptc{S}_0+ \abdb{q}\\
    \ptc{T}^3&=& \ptc{S}_0+ \abdb{r}
\end{eqnarray*}
As discussed in Section \ref{sec:facefigure}, the conic  $\ptc{S}_0$ can move along two conics in double contact with it, while staying in double contact with both. The family of such conics can be expressed using  $\ptc{S}_0$, the chords of contact $\myln{q}$ and $\myln{r}$, and a free parameter $a$:
\begin{eqnarray}\label{eq:Sijk-det}
    \ptc{S}_1&=& (1-a^2)\,\ptc{S}_0+\abdb{q} + 2\, a\, \abpr{q}{r} + \abdb{r}
\end{eqnarray}
And we find the conics in double contact with $\ptc{T}^1-\ptc{T}^3$ and $\ptc{T}^1-\ptc{T}^2$
\begin{eqnarray*}  
    \ptc{S}_2&=& (1-b^2)\,\ptc{S}_0+\abdb{p} + 2\, b\, \abpr{p}{r} + \abdb{r}\\
    \ptc{S}_3&=& (1-c^2)\,\ptc{S}_0+\abdb{p} + 2\, c\, \abpr{p}{q} + \abdb{q}
\end{eqnarray*}
where $b$ and $c$ are free constants. 

That $\ptc{S}_3$ is in double contact with $\ptc{T}^1$ we see by the fact that a linear combination of them gives a double line:
\begin{equation} \label{eq:level1-doubleLine}
     \ptc{S}_3-(1-c^2)\ptc{T}^1=(c\myln{p} + \myln{q})^2 
\end{equation}
where $c\myln{p} + \myln{q}$ is the chord of contact. 
The symmetry shows that there are such contacts between the other conics as well.

The last to show is that the three conics $\ptc{S}_1$, $\ptc{S}_2$ and $\ptc{S}_3$ are in double contact with an eighth conic. This one is totally determined by the other elements in the structure, and the formula is given by:
\begin{eqnarray} \label{eq:T0-det}
\ptc{T}^0&=&(1 - a^2 - b^2 - 2 a b c-c^2)\, \ptc{S}_0\\ \nonumber & +& (1-a^2)\, \abdb{p} + (1 -  b^2)\, \abdb{q}+ (1 - c^2)\, \abdb{r}\\ \nonumber & +& 2 (a + b\, c)\, \abpr{q}{r} + 2 (b + a c)\, \abpr{p}{r} + 2 (c+a b )\, \abpr{p}{q}
\end{eqnarray}

We find that $\ptc{T}^0$ and $\ptc{S}_1$ are in double contact because we have:
\begin{equation}
\begin{split}
    &(1-a^2)\,\ptc{T}^0-(1 - a^2 - b^2 - 2 a b c-c^2)\,\ptc{S}_1\\&=\{(1 - a^2) \,\myln{p} + (c + a b) \,\myln{q} + (b + a c) \,\myln{r}\}^2
    \end{split}
\end{equation}where
$$(1 - a^2) \,\myln{p} + (c + a b) \,\myln{q} + (b + a c) \,\myln{r}$$is the line through the contact points of $\ptc{T}^0$ and $\ptc{S}_1$. From the symmetry, we see that $\ptc{S}_2$ and $\ptc{S}_3$ are also in double contact with $\ptc{T}^0$.

\subsection{The algebraic proof}

The algebraic proof of the 8-conic theorem reverses this process. Given a Penrose 7-configuration, it calculates  arguments $(\myln{p}, \myln{q}, \myln{r}, a,b,c)$ which produces this configuration, and uses them, along with $\ptc{S}_0$, to calculate $\ptc{T}^0$.

Let $\ptc{S}_i$ and $\ptc{T}^j$ be given as a Penrose 7-configuration, without $\ptc{T}^0$.
With suitable normalization of $\ptc{T}_i$, we can apply Eq. \ref{eq:t1} to obtain:
\begin{equation}\label{eq:t1-2}
     \abdb{p} = \ptc{S}_0 - \ptc{T}^1,~~~~ \abdb{q} = \ptc{S}_0 - \ptc{T}^2,~~~~\abdb{r} = \ptc{S}_0 - \ptc{T}^3
\end{equation}

Each of the real parameters $\{a,b,c\}$ corresponds to a complete face of the given Penrose 7-configuration. For example, $a$ corresponds to face $\face^{23}_{01}$, whose chords of contact pass through $P^{23}_{01}$. We seek $a$ so that $\ptc{S}_1$ satisfies Eq. \ref{eq:Sijk-det}. By Section \ref{sec:facefigure}, this value can be directly obtained in terms of $\ptc{S}_0, \ptc{q},$ and $\ptc{r}$.  Similarly for $b$ and $c$.
Finally, we apply Eq. \ref{eq:T0-det} to obtain the desired eighth conic $\ptc{T}^0$, completing the cube and the proof.

This proof can also be applied to prove the 8-quadric theorem. One has only to replace double contact with ring contact, and the line pairs and double lines with plane pairs and double planes, resp.
}

\newcommand{\subd}[2]{\text{\scriptsize $\begin{vmatrix} #1 \\ #2  \end{vmatrix}$}}
\newcommand{\dsubd}[2]{\text{\scriptsize $\overline{\begin{vmatrix} #1 \\ #2 \end{vmatrix}}$}}

\section{Determinant structure and generalizations}
\label{sec:det-structure}

The astute reader may recognize in the formulas of the last section, the following determinant forms with polynomial coefficients (vertical bars denote determinants and, by homogeneous coordinates, we may ignore a scaling by $-1$):
\begin{equation} \label{eq:Tijk-det}
    \ptc{T}^1=\begin{vmatrix}
\ptc{S}_0 & \myln{p} \\
\myln{p} & -1
\end{vmatrix}, \qquad
\ptc{T}^2=\begin{vmatrix}
\ptc{S}_0 & \myln{q} \\
\myln{q} & -1
\end{vmatrix},\qquad
\ptc{T}^3=\begin{vmatrix}
\ptc{S}_0 & \myln{r} \\
\myln{r} & -1
\end{vmatrix},
\end{equation}

\begin{equation}  \label{eq:Sijk-det-3}
      \ptc{S}_1=\begin{vmatrix}
\ptc{S}_0 & \myln{q} & \myln{r} \\
\myln{q} & -1 & a \\
\myln{r} & a & -1
\end{vmatrix}, \qquad
\ptc{S}_2=\begin{vmatrix}
\ptc{S}_0 & \myln{p} & \myln{r} \\
\myln{p} & -1 & b \\
\myln{r} & b & -1
\end{vmatrix}, \qquad
\ptc{S}_3=\begin{vmatrix}
\ptc{S}_0 & \myln{p} & \myln{q} \\
\myln{p} & -1 & c \\
\myln{q} & c & -1
\end{vmatrix},
\end{equation}

\begin{equation} \label{eq:T0-det-2}
      \ptc{T}^0=\begin{vmatrix}
\ptc{S}_0 & \myln{p} & \myln{q} & \myln{r} \\
\myln{p} & -1 & c & b \\
\myln{q} & c & -1 & a \\
\myln{r} & b & a  & -1
\end{vmatrix},
\end{equation}
or in a notation more conducive of the general discussion that follows:
\begin{equation}
    \ptc{S}_{\{1\}}=\begin{vmatrix}
\ptc{S}_0 & \myln{p}_1 \\
\myln{p}_1 & d_1
\end{vmatrix}, \qquad
\ptc{S}_{\{2\}}=\begin{vmatrix}
\ptc{S}_0 & \myln{p}_2 \\
\myln{p}_2 & d_2
\end{vmatrix},\qquad
\ptc{S}_{\{3\}}=\begin{vmatrix}
\ptc{S}_0 & \myln{p}_3 \\
\myln{p}_3 & d_3
\end{vmatrix},
\end{equation}

\begin{equation}
      \ptc{S}_{\{2,3\}}=\begin{vmatrix}
\ptc{S}_0 & \myln{p}_2 & \myln{p}_3 \\
\myln{p}_2 & d_2 & a_{23} \\
\myln{p}_3 & a_{23} & d_3
\end{vmatrix}, \qquad
\ptc{S}_{\{1,3\}}=\begin{vmatrix}
\ptc{S}_0 & \myln{p}_1 & \myln{p}_3 \\
\myln{p}_1 & d_1 & a_{13} \\
\myln{p}_3 & a_{13} & d_3
\end{vmatrix}, \qquad
\ptc{S}_{\{1,2\}}=\begin{vmatrix}
\ptc{S}_0 & \myln{p}_1 & \myln{p}_2 \\
\myln{p}_1 & d_1 & a_{12} \\
\myln{p}_2 & a_{12} & d_2
\end{vmatrix},
\end{equation}

\begin{equation}
      \ptc{S}_{\{1,2,3\}}=\begin{vmatrix}
\ptc{S}_0 & \myln{p}_1 & \myln{p}_2 & \myln{p}_3 \\
\myln{p}_1 & d_1 & a_{12} & a_{13} \\
\myln{p}_2 & a_{12} & d_2 & a_{23} \\
\myln{p}_3 & a_{13} & a_{23}  & d_3
\end{vmatrix},
\end{equation}
where $d_j$ are homogeneous parameters.  

{\it We shall find that the proof of the Penrose theorem, explicit expressions for the chords of content, as well as even more geometric information about the configuration is contained in the subdeterminants of $\ptc{S}_{\{1,2,3\}}$, and furthermore that hypercube analogues naturally arise in this method of proof}.  In order to describe this we introduce the following notation (c.f. \cite{sylvester}).  A $k \times k$ subdeterminant is determined by choosing $k$ rows and $k$ columns of the original matrix. We write $\subd{I_r}{I_c}$ for this subdeterminant, where $I_r$ ($I_c$) is a subset of the 0-based indices of the chosen rows (columns).  The Desnanot–Jacobi identity (see e.g. \cite{Wikipedia_DodgsonCondensation}) states that
\begin{theorem}(The Desnanot–Jacobi identity)\label{theorem:desnanotjacobi} If $j \neq k$ and $j,k \notin I_r \cup I_c$ then
\begin{equation}
    \subd{I_r \cup \{k\}}{I_c \cup \{k\}}\subd{I_r \cup \{j\}}{I_c \cup \{j\}}-\subd{I_r}{I_c}\subd{I_r \cup \{j,k\}}{I_c \cup \{j,k\}}=\subd{I_r \cup \{j\}}{I_c \cup \{k\}}\subd{I_r \cup \{k\}}{I_c \cup \{j\}}
\end{equation}
\end{theorem}

In order to apply this to conics, we henceforth fix the full matrix to be (i) symmetric, (ii) have a second order homogeneous polynomial in the top left corner, (iii) have first order homogeneous polynomials in the remainder of the first row and first column (both $0$ indexed) and (iv). have scalars as the remainder of the entries, i.e.,
\begin{equation}
    \ptc{S}_{\{1,2,...,n\}}:=\subd{012...n}{012...n}=
    \begin{vmatrix}
        \ptc{S_0} & \myln{p}_1 & \myln{p}_2 & \cdots & \myln{p}_{n-1} & \myln{p}_n\\
\myln{p}_1 & d_1 & a_{12} & \cdots & a_{1,n-1} & a_{1n}\\
\myln{p}_2 & a_{12} & d_2 & \cdots & a_{2,n-1} & a_{2n}\\
\vdots & \vdots & \vdots & \ddots  & \vdots     & \vdots \\
\myln{p}_{n-1} & a_{1,n-1} & a_{2,n-1} & \cdots & d_{n-1} & a_{n-1,n}\\
\myln{p}_n & a_{1n} & a_{2n} & \cdots & a_{n-1,n} & d_n
    \end{vmatrix}.
\end{equation}
  The subdeterminants can be classified as follows:
\begin{itemize}
    \item if $0 \in I_c \cap I_r$ then $\subd{I_c}{I_r}$ is a conic
    \item if $0 \in I_c$ but $0 \notin I_r$ then $\subd{I_c}{I_r}$ is a line
    \item if $0 \notin I_c$ and $0 \notin I_r$ then $\subd{I_c}{I_r}$ is a scalar.
\end{itemize}

\begin{figure}[h!]
    \centering
    \includegraphics[width=0.6\textwidth]{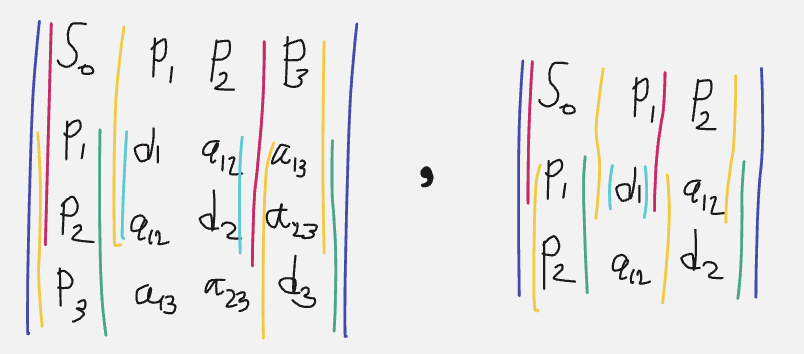} 
    \includegraphics[width=0.6\textwidth]{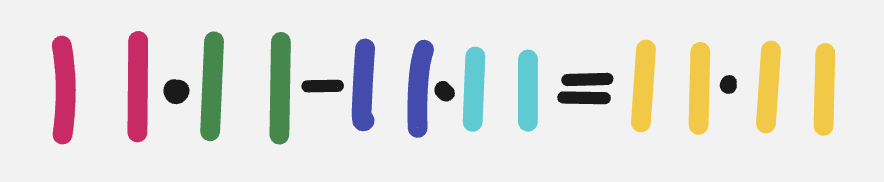}

    \caption{\emph{Above}: Color-coding of the five determinants involved in Lemma \ref{lemma:contactsubset} for the first (\emph{right}) and second (\emph{left}) levels of the Penrose cube, $n=3$, in the respective cases $\Omega=\{1,2\}$, $k=3$ and $\Omega=\{1\}$, $k=2$.  \emph{Below}: The lemma in graphical form. Note the two yellow determinants are by symmetry identical.}
    \label{fig:five-det-figure}
\end{figure}

Take $\Omega \subset \{1,2,...,n\}$.  To proceed, we give names to some special subdeterminants (note: this is not exhaustive, an exhaustive description of the case $n=3$ will be given in Sec. \ref{sec:subdet}):
\begin{equation}
    \ptc{S}_{\Omega}:=\subd{\Omega \cup \{0\}}{\Omega \cup \{0\}},
\end{equation}
\begin{equation}
    \myln{p}_{\Omega,k}:=\subd{\Omega \cup \{k\}}{\Omega \cup \{0\}},
\end{equation}
\begin{equation}
    f_{\Omega}:=\subd{\Omega}{\Omega},\: \Omega \neq \emptyset, \quad f_{\emptyset}=1.
\end{equation}

Note that the symmetry of the full determinant implies that
\begin{equation}
    \subd{\Omega \cup \{k\}}{\Omega \cup \{0\}}=\subd{\Omega \cup \{0\}}{\Omega \cup \{k\}}.
\end{equation}
Therefore, specializing Theorem \ref{theorem:desnanotjacobi} to the case $I_c=I_r=\Omega$, $j=0$ gives
\begin{lemma} \label{lemma:contactsubset}
\begin{equation}\label{eq:contactsubset}
    f_{\Omega \cup \{k\}}\ptc{S}_{\Omega}-f_{\Omega}\ptc{S}_{\Omega \cup \{k\}}=\myln{p}_{\Omega,k}^2,
\end{equation}
\end{lemma}
In the case where $f_{\Omega}, f_{\Omega \cup \{k\}},\myln{p}_{\Omega,k} \neq 0$, this shows that the subset lattice $\ptc{S}_{\Omega}$, ${\Omega}\subset \{1,2,...,n\}$ can be interpreted as a double contact hypercube graph.  Figures \ref{fig:cubesubsets} and \ref{fig:hcubesubsets} illustrate this in the cases $n=3$ and $n=4$ respectively.  Since these restrictions hold generically, any set of parameters not satisfying them are given as a limit of configurations which do.  Thus we refer to the graph in the case where $f_{\Omega}, f_{\Omega \cup \{k\}}$ or $\myln{p}_{\Omega,k}$ do vanish as a \emph{generalized double contact graph}.

\begin{figure}[h!]
     \centering
\begin{tikzpicture}[every node/.style={circle,draw,inner sep=2pt}, scale=1.5]

\node (S0) at (1,4)  {$\ptc{S}_{\emptyset}$};  
\node (S1) at (0,3)  {$\ptc{S}_{\{1\}}$};  
\node (S2) at (1,3)  {$\ptc{S}_{\{2\}}$};  
\node (S3) at (2,3)  {$\ptc{S}_{\{3\}}$};  

\node (S12) at (0,2)  {$\ptc{S}_{\{1,2\}}$};  
\node (S13) at (1,2)  {$\ptc{S}_{\{1,3\}}$};  
\node (S23) at (2,2)  {$\ptc{S}_{\{2,3\}}$};  

\node (S123) at (1,1)  {$\ptc{S}_{\{1,2,3\}}$};  

\draw (S0) -- (S1);
\draw (S0) -- (S2);
\draw (S0) -- (S3);

\draw (S1) -- (S12);
\draw (S1) -- (S13);
\draw (S2) -- (S12);
\draw (S2) -- (S23);
\draw (S3) -- (S13);
\draw (S3) -- (S23);

\draw (S12) -- (S123);
\draw (S13) -- (S123);
\draw (S23) -- (S123);

\end{tikzpicture}
\caption{A graph showing the subset structure of the Penrose cube.}
     \label{fig:cubesubsets}
 \end{figure}

 \begin{figure}[h!]
     \centering
\begin{tikzpicture}[every node/.style={circle,draw,inner sep=2pt}, scale=1.5]


\node (S0) at (2.5,5)  {$\ptc{S}_{\emptyset}$};  
\node (S1) at (1,4)  {$\ptc{S}_{\{1\}}$};  
\node (S2) at (2,4)  {$\ptc{S}_{\{2\}}$};  
\node (S3) at (3,4)  {$\ptc{S}_{\{3\}}$};  
\node (S4) at (4,4)  {$\ptc{S}_{\{4\}}$};  

\node (S12) at (0,3)  {$\ptc{S}_{\{1,2\}}$};  
\node (S13) at (1,3)  {$\ptc{S}_{\{1,3\}}$};  
\node (S14) at (2,3)  {$\ptc{S}_{\{1,4\}}$};  
\node (S23) at (3,3)  {$\ptc{S}_{\{2,3\}}$};
\node (S24) at (4,3)  {$\ptc{S}_{\{2,4\}}$};
\node (S34) at (5,3)  {$\ptc{S}_{\{3,4\}}$};  

\node (S123) at (1,2)  {$\ptc{S}_{\{1,2,3\}}$};  
\node (S124) at (2,2)  {$\ptc{S}_{\{1,2,4\}}$};  
\node (S134) at (3,2)  {$\ptc{S}_{\{1,3,4\}}$};  
\node (S234) at (4,2)  {$\ptc{S}_{\{2,3,4\}}$};  

\node (S1234) at (2.5,1)  {$\ptc{S}_{\{1,2,3,4\}}$};  

\draw (S0) -- (S1);
\draw (S0) -- (S2);
\draw (S0) -- (S3);
\draw (S0) -- (S4);

\draw (S1) -- (S12);
\draw (S1) -- (S13);
\draw (S1) -- (S14);
\draw (S2) -- (S12);
\draw (S2) -- (S23);
\draw (S2) -- (S24);
\draw (S3) -- (S13);
\draw (S3) -- (S23);
\draw (S3) -- (S34);
\draw (S4) -- (S14);
\draw (S4) -- (S24);
\draw (S4) -- (S34);

\draw (S12) -- (S123);
\draw (S13) -- (S123);
\draw (S23) -- (S123);
\draw (S12) -- (S124);
\draw (S14) -- (S124);
\draw (S24) -- (S124);
\draw (S13) -- (S134);
\draw (S14) -- (S134);
\draw (S34) -- (S134);
\draw (S23) -- (S234);
\draw (S24) -- (S234);
\draw (S34) -- (S234);

\draw (S123) -- (S1234);
\draw (S124) -- (S1234);
\draw (S134) -- (S1234);
\draw (S234) -- (S1234);

\end{tikzpicture}
\caption{A graph showing the subset structure of the Penrose hypercube.}
     \label{fig:hcubesubsets}
 \end{figure}

\subsection{Vanishing parameters}

Vanishing of the chord of contact $\myln{p}_{\Omega,k}$, or the simultaneous vanishing $f_{\Omega}=f_{\Omega \cup \{k\}}=0$, in equation \eqref{eq:contactsubset} leads to a situation where equation \eqref{eq:dblln} is not strictly satisfied.  However, by construction, the resulting configuration is the limit of configurations where equation \eqref{eq:dblln} does hold.  Before considering these more extreme cases, we consider what the vanishing of each of the scalar entries of the matrix means.

\subsubsection{Vanishing of $a_{jk}$}
\begin{figure}[h!]
    \centering
    \includegraphics[width=0.4\textwidth]{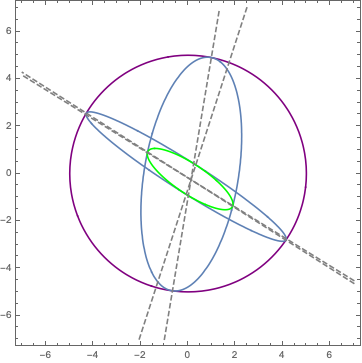} 
    \includegraphics[width=0.4\textwidth]{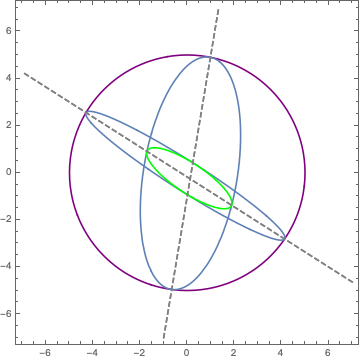}
    \caption{A Penrose face with $\ptc{S}_{\emptyset}$ in purple, $\ptc{S}_{\{j\}}$, $\ptc{S}_{\{k\}}$ in blue, $\ptc{S}_{\{j,k\}}$ in green, and chords of contact gray and dashed.  \emph{Left:} as $a_{jk}$ becomes small, the pairs of chords of contact approach each other.  \emph{Right:} when $a_{jk}$ vanishes the four chords of contact of the face collapse to two.}
    \label{fig:fusingchords}
\end{figure}

Taking $a_{jk}=0$ and expanding the determinants for the second layer chords of contact yields
\begin{equation}
    \myln{p}_{\{k\},j}=-d_k\myln{p}_j
\end{equation}
and
\begin{equation}
    \myln{p}_{\{j\},k}=-d_j\myln{p}_k.
\end{equation}
In summary, $a_{jk}=0$ implies that the pairs of chords of contact of the face $\{\ptc{S}_{\emptyset},\ptc{S}_{\{j\}},\ptc{S}_{\{k\}},\ptc{S}_{\{j,k\}}\}$ as shown in Figure \ref{fig:fusingchords}.

\subsubsection{Vanishing of $d_j$}

When $d_j=0$, equation \eqref{eq:contactsubset} implies that
\begin{equation}
    \ptc{S}_{\{j\}}=-\myln{p}_j^2
\end{equation}
since $f_{\{j\}}=d_j$, and
\begin{equation}
    f_{\{j,k\}}\ptc{S}_{\{j\}}=\myln{p}_{\{j\},k}^2.
\end{equation}
Thus the contact chords for both edges adjacent to $\ptc{S}_{\{j\}}$ collapse together, and $\ptc{S}_{\{j\}}$ itself becomes the same line squared.  Furthermore, this line, $\myln{p}_j$, intersects $\ptc{S}_{\emptyset}$ and $\ptc{S}_{\{j,k\}}$ in the same pair of points, because
\begin{equation} \label{eq:sjkexp}
    \ptc{S}_{\{j,k\}}=(d_jd_k-a_{jk}^2) \ptc{S}_{\emptyset}+2a_{jk}\myln{p}_k \myln{p}_j-d_k\myln{p}_j^2-d_j\myln{p}_k^2.
\end{equation}
and when $d_j=0$, this becomes
\begin{equation}
    \ptc{S}_{\{j,k\}}=-a_{jk}^2 \ptc{S}_{\emptyset}+\myln{p}_j(2a_{jk}\myln{p}_k-d_k\myln{p}_j).
\end{equation}
Therefore, \emph{as far as degenerate double contact is concerned, $\ptc{S}_j$ can be thought of as much as a point pair as a double line}, in spite of the fact that the description of the current section treats all conics as point-wise and thus, to begin with, only yields to an interpretation of $\ptc{S}_j$ as a double line. 

\begin{figure}[h!]
    \centering
    \includegraphics[width=0.4\textwidth]{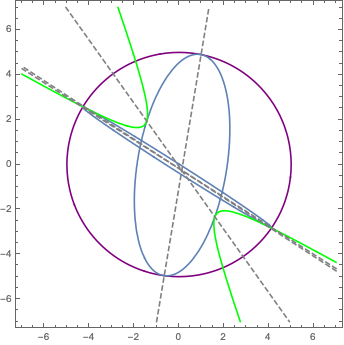}
    \caption{The coloring is the same as Figure \ref{fig:fusingchords}. As $d_j$ goes to $0$, $\ptc{S}_{\{j\}}$ goes to $\myln{p}_j^2$ and $\myln{p}_j$, and $\myln{p}_{\{j\},2}$ come together.}
    \label{fig:d1tozero}
\end{figure}

\subsubsection{Vanishing of $f_{\{j,k\}}$}\label{sec:fjk}

The equation 
\begin{equation}
    f_{\{j,k\}}=d_jd_k-a_{jk}^2=0
\end{equation}
holds exactly when 
\begin{equation}
    a_{jk}=\pm \sqrt{d_jd_k}.
\end{equation}
putting this into equation \eqref{eq:sjkexp} yields
\begin{equation}
    \ptc{S}_{\{j,k\}}=2a_{jk}\myln{p}_j\myln{p}_k-d_j \myln{p}_k^2-d_k \myln{p}_j^2=-(\sqrt{d_j}\myln{p}_k\pm \sqrt{d_k}\myln{p}_j)^2
\end{equation}
which can be seen to intersect $\ptc{S}_{\{j\}}$ and $\ptc{S}_{\{k\}}$ in the same pair of points because
\begin{equation}
    d_k \ptc{S}_{\{j\}}-d_j \ptc{S}_{\{k\}}=(\sqrt{d_j}\myln{p}_k-\sqrt{d_k}\myln{p}_j)(\sqrt{d_j}\myln{p}_k+\sqrt{d_k}\myln{p}_j).
\end{equation}
Thus \emph{when $f_{\{j,k\}}$ vanishes (but $d_j, d_k \neq 0$) $\ptc{S}_{\{j,k\}}$ becomes a double line which intersects $\ptc{S}_{\{j\}}$ and $\ptc{S}_{\{k\}}$ in the same pair of points.}

\subsubsection{Vanishing of $f_{\Omega}$ in general}

Application of Sylvester's determinant identity (\cite{sylvester}), which generalizes the Desnanot-Jacobi identity, allows for relations to be expressed describing any face $\{\ptc{S}_{\Omega},\ptc{S}_{\Omega \cup \{j\}},\ptc{S}_{\Omega \cup \{k\}},\ptc{S}_{\Omega \cup \{j,k\}}\}$ in a way analogous to the first layer face $\Omega=\emptyset$ considered above, provided that $f_{\Omega} \neq 0$:
\begin{equation}
    f_{\Omega}^2 \ptc{S}_{\Omega \cup \{j,k\}}=\begin{vmatrix}
        \ptc{S}_{\Omega} & \myln{p}_{\Omega,j}& \myln{p}_{\Omega,k}\\
        \myln{p}_{\Omega,j} & f_{\Omega \cup \{j\}} & g_{\Omega,j,k}\\
        \myln{p}_{\Omega,k} & g_{\Omega,j,k} & f_{\Omega \cup \{j\}}
    \end{vmatrix}
\end{equation}
where $g_{\Omega,j,k}=\subd{\Omega \cup \{j\}}{\Omega \cup \{k\}}$, the $2 \times 2$ symmetric subdeterminants of which are $\ptc{S}_{\Omega \cup \{j\}}$ and $\ptc{S}_{\Omega \cup \{k\}}$ by equation \eqref{eq:contactsubset}.  Therefore, 
\begin{itemize}
    \item when $f_{\Omega} \neq 0$ and $f_{\Omega \cup \{j\}}=0$, $\ptc{S}_{\Omega \cup \{j\}}$ becomes a double line intersecting $\ptc{S}_{\Omega}$ and $\ptc{S}_{\Omega \cup \{j,k\}}$ in the same pair of points and
    \item when $f_{\Omega},f_{\Omega \cup \{j\}}, f_{\Omega \cup \{k\}} \neq 0$ and $f_{\Omega \cup \{j,k\}}=0$, $\ptc{S}_{\Omega \cup \{j,k\}}$ becomes a double line intersecting $\ptc{S}_{\Omega \cup \{j\}}$ and $\ptc{S}_{\Omega \cup \{k\}}$ in the same pair of points.
\end{itemize}
Putting these together, we obtain the general statement:
\begin{lemma}
    if $f_{\Omega}=0$ but none of the previous-layer-values $f_{\Gamma}$, $\Gamma \subset \Omega$, vanish, then $\ptc{S}_{\Omega}$ is a double line intersecting its neighbors $\ptc{S}_{\Omega\setminus \{j\}}$ and $\ptc{S}_{\Omega\cup \{k\}}$ in the same pair of points.
\end{lemma}
We leave a more thorough treatment of the cases where $f_{\Omega}$ vanishes together with previous layer values $f_{\Gamma}$ to future work, but will focus on one particular case where this occurs since it allows us to explicitly express the parameter constraints corresponding to the configurations of the classical theorems of projective geometry mentioned in the Introduction. 

\begin{remark}
    The circumstance where a double line can also be interpreted as a point pair is essentially related to the condition $f_{\Omega}=0$.  An example of a case where this need not occur is when $\ptc{S}_0=\myln{q}^2$, i.e., it is a double line.  Then each $\ptc{S}_{\{j\}}$ is a line pair, with the lines of the pair crossing at an arbitrary point on $\myln{q}$.
\end{remark}

\subsection{Special parameter values and special cases: Dual Salmon, Brianchon, Desargues, etc.}

As pointed out in section \ref{sec:fjk},
\begin{equation}
    f_{\{i,j\}}=0 \iff a_{ij}=\sigma_{ij}\sqrt{d_i d_j}
\end{equation}
where $\sigma_{ij}=\pm1$, and
\begin{equation}
    S_{\{i,j\}}=-(\sqrt{d_j}\myln{p}_i-\sigma_{ij}\sqrt{d_i}\myln{p}_j)^2.
\end{equation}
In summary, the $\ptc{S}_{\{i,j\}}$ are double lines/point pairs.  

We now show that $f_{\{1,2,3\}}=0$ on top of the vanishing of the three $f_{\{i,j\}}$ implies the concurrency of the three lines.

\begin{align}
    f_{\{1,2,3\}}&=\begin{vmatrix}
        d_1 & \sigma_{12}\sqrt{d_1 d_2} & \sigma_{13}\sqrt{d_1 d_3} \\
\sigma_{12}\sqrt{d_1 d_2} & d_2 &  \sigma_{23}\sqrt{d_2 d_3}\\
\sigma_{13}\sqrt{d_1 d_3} & \sigma_{23}\sqrt{d_2 d_3} & d_3 
    \end{vmatrix}\\
    &=\sqrt{d_1d_2d_3}\begin{vmatrix}
        \sqrt{d_1} & \sigma_{12}\sqrt{d_1} & \sigma_{13}\sqrt{d_1} \\
\sigma_{12}\sqrt{d_2} & \sqrt{d_2} &  \sigma_{23}\sqrt{d_2}\\
\sigma_{13}\sqrt{d_3} & \sigma_{23}\sqrt{d_3} & \sqrt{d_3} 
    \end{vmatrix}\\
    &=d_1d_2d_3\begin{vmatrix}
        1 & \sigma_{12} & \sigma_{13} \\
\sigma_{12} & 1 &  \sigma_{23}\\
\sigma_{13} & \sigma_{23} & 1 
    \end{vmatrix}\\
    &=d_1d_2d_3(1-\sigma_{12}^2-\sigma_{13}^2-\sigma_{23}^2+2\sigma_{12}\sigma_{13}\sigma_{23})\\
    &=2d_1d_2d_3(\sigma_{12}\sigma_{13}\sigma_{23}-1).
\end{align}
This vanishes exactly when
\begin{equation}
    \sigma_{23}=\sigma_{12}\sigma_{13}.
\end{equation}
The dependency of the three lines then follows:
\begin{equation}
\sqrt{d_3}(\sqrt{d_2}\myln{p}_1-\sigma_{12}\sqrt{d_1}\myln{p}_2)-\sqrt{d_2}(\sqrt{d_3}\myln{p}_1-\sigma_{13}\sqrt{d_1}\myln{p}_3)=-\sqrt{d_1}\sigma_{12}(\sqrt{d_3}\myln{p}_2-\sigma_{23}\sqrt{d_2}\myln{p}_3).
\end{equation}
This is exactly the dual Salmon Theorem shown in Figure \ref{fig:SalmonInIntro} on the right, in the introduction.  Special cases include the Brianchon Theorem, when the first layer conics $\ptc{S}_{\{j\}}$ are line pairs (pairs of tangents to $\ptc{S}$) and the Desargues Theorem when $\ptq{S}$ is a double line.  We also note that if we take $\ptq{S}$ to be a double line but drop the condition $f_{\{1,2,3\}}=0$, the final conic will in general be nondegenerate and the Braikenridge–Maclaurin theorem (the converse of Pascal's theorem) results because the final conic $\ptc{S}_{\{1,2,3\}}$ must go through all three of the point pairs for the previous layer conics, see Table \ref{tab:specializations}.

\begin{table}[h]
    \centering
    \begin{tabular}{|c|c|c|c|c|}
    \hline
    Theorem & $f_{\{1,2,3\}}$ & $\{\ptc{S}_{\{1\}}, \ptc{S}_{\{2\}}, \ptc{S}_{\{3\}}\}$ & $\ptc{S}_{\emptyset}$ & Figure \\ \hline \hline
    Dual Salmon & 0 & -  & - & Figure \ref{fig:SalmonInIntro} right \\ \hline
    Brianchon & 0 & line pairs & - & Figure \ref{fig:PascalInIntro} right \\ \hline
    Pappos & 0 &line pairs & line pair & Figure \ref{fig:PappusInIntro} right \\ \hline
    Braikenridge-Maclaurin & - & - & double line & Figure \ref{fig:PascalInIntro} left\\ \hline
    Desargues & 0 & - & double line  & - \\ \hline
    \end{tabular}
    \caption{For all entries, we assume $f_{\{1,2\}}=f_{\{1,3\}}=f_{\{2,3\}}=0$ which implies, as is shown in \ref{sec:fjk}, that $\ptc{S}_{\{1,2\}},\ptc{S}_{\{1,3\}}$ and $\ptc{S}_{\{2,3\}}$ are double lines. A '-' indicates no special condition. The table gives the names of well-known special cases, the accompanying degeneracy conditions as well as the corresponding figure in the introduction.}
    \label{tab:specializations}
\end{table}

\begin{remark}
    It can also be checked that $\ptc{S}_{\{1,2,3\}}=0$ in the case we just considered, however the behavior we saw does not hold generally when one of the $\ptc{S}_{\Omega}$ vanishes, e.g., when $\ptc{S}_{\emptyset}=0$ there is no requirement that the double lines of the next layer are concurrent.  Thus the forced concurrency of $\myln{p}_{\{1,2\},3}$, $\myln{p}_{\{1,3\},2}$ and $\myln{p}_{\{2,3\},1}$, which can be interpreted line-wise as $\ptc{S}_{\{1,2,3\}}$ being a double point, is not reducible to the vanishing of the point-wise conic but is essentially connected to the simulteneous vanishing of the $f_{\{i,j\}}$ and $f_{\{1,2,3\}}$.
\end{remark}

\subsection{Face conics} \label{sec:face-conics}

In order to further interpret these cases as well as to generally learn more about the geometry of the configuration, we name the following subdeterminant conics which, as we shall see, have remarkable properties:
\begin{equation} \label{eq:face-conics-2}
    \ptc{H}_{\Omega,j,k}=\subd{\Omega \cup \{0,j\}}{\Omega \cup \{0,k\}}.
\end{equation}
Taking $I_c=I_r=\Omega \cup \{0\}$ in the Desnanot–Jacobi identity (Theorem \ref{theorem:desnanotjacobi}) gives
\begin{lemma}
\begin{equation} \label{eq:face-conics}
    \ptc{S}_{\Omega \cup \{k\}}\ptc{S}_{\Omega \cup \{j\}}-\ptc{S}_{\Omega}\ptc{S}_{\Omega \cup \{j,k\}}=\ptc{H}_{\Omega,j,k}^2
\end{equation}
and therefore the conic $\ptc{H}_{\Omega,j,k}$ vanishes on each of the 4 pairs of contact points between conics in the face $\{\ptc{S}_{\Omega},\ptc{S}_{\Omega \cup \{k\}},\ptc{S}_{\Omega \cup \{j\}},\ptc{S}_{\Omega \cup \{j,k\}}\}$.  Thus we call the $\ptc{H}_{\Omega,j,k}$ {\it face conics.}
\end{lemma}

\subsection{Face Diagonals} \label{sec:face-diag}

The quantities we have considered thus far all pertain to a particular face, edge or node of the (generalized) double contact graph.  This was because neither $I_r$ nor $I_c$ had more than one element not in the other.  In this subsection we discuss lines whose geometric interpretation requires considering an entire cube.  

\begin{figure}[ht!]
    \centering
    \includegraphics[width=0.6\textwidth]{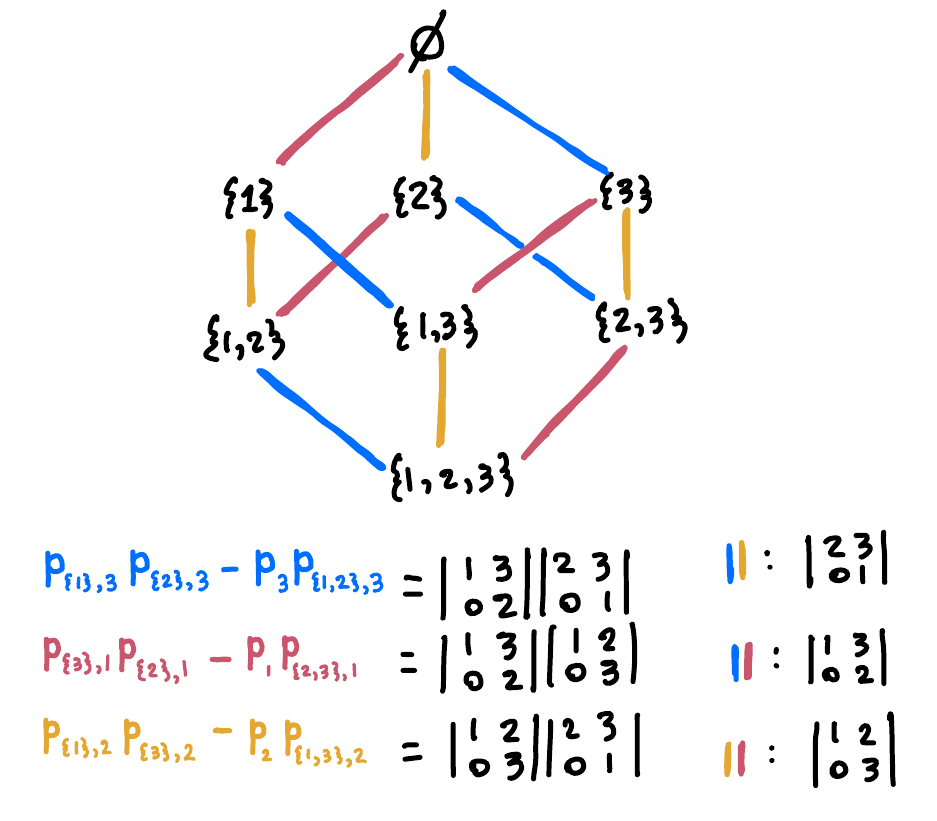}
    \caption{Illustration of the algebra of Section \ref{sec:face-diag}}
    \label{fig:facediag}
\end{figure}

Fix $n=3$.  Then the Desnanot-Jacobi identity applied to $I_r=\{m\}$, $I_c=\{0\}$, becomes
\begin{equation}\label{eq:face-diag}
    \myln{p}_{\{j\},m}\myln{p}_{\{k\},m}-\myln{p}_m \myln{p}_{\{j,k\},m}=\subd{mj}{0k}\subd{mk}{0j}.
\end{equation}
Letting $m,j$ and $k$ take on all the possible index values $1,2,3$, we get the equations shown color-coded in figure \ref{fig:facediag}.  The blue and yellow equations have the common factor $\subd{23}{01}$, thus the line this factor describes goes through the face points of the opposite faces bordered by the blue and yellow edges.  Similar statements hold for the faces bordered in blue and red, and red and yellow, as indicated by the legend in figure \ref{fig:facediag}.  Each face conic can be associated with the index which appears together with $0$:
\begin{equation}
    \myln{q}_k=\subd{mj}{0k}.
\end{equation}
Consulting figure \ref{fig:facediag}, we can see that $\myln{q}_k$ connects the face each of whose node labels is a subset containing $k$ to the face whose nodes all exclude $k$.

\subsection{Relationship of face conics and face diagonals}

In this subsection we demonstrate the following property of the face conics: \emph{the face conics for the 4 faces remaining when a pair of opposite faces are removed coincide in two points of the face diagonal corresponding to the two excluded faces.}  For short, we refer to the face whose index subsets at its vertices all contain $k$ as $\mathcal{F}_k$ and the opposite face, none of whose index subsets contain $k$, as $\mathcal{F}_{\bar{k}}$. Consult Figure \ref{fig:facediag}. Note that, in this notation, the face conic for $\mathcal{F}_k$ is $\ptc{H}_{\{k\},\ell,m}$ and the face conic for $\mathcal{F}_{\bar{k}}$ is $\ptc{H}_{\emptyset,\ell,m}$ where $\{k,\ell,m\}=\{1,2,3\}$.

First, we apply the Desnanot-Jacobi identity with $I_r=\{\ell\}$, $I_c=\{m\}$ and $k=0$ to obtain
\begin{equation}\label{faceconicsanddiag1}
    g_{\{j\},\ell,m}\ptc{H}_{\emptyset,\ell,m}-g_{\emptyset,\ell,m}\ptc{H}_{\{j\},\ell,m}= \myln{q}_{\ell}\myln{q}_{m},
\end{equation}
where $\{j,\ell,m\}=\{1,2,3\}$.  This implies that \emph{the four meeting points of the two face conics for opposite faces $\mathcal{F}_{\bar{j}}$ and $\mathcal{F}_j$ lie on the line pair $\myln{q}_{\ell}\myln{q}_{m}$, two on each, since each line of the line pair intersects each of the conics in exactly two points.}

Two face conics for adjacent faces meet in the two contact points for their shared double contact edge.  To discover where the other two intersection points lie, we apply Sylvester's identity to
\begin{equation}
    \ptc{H}_{\{j\},\ell,m}=\begin{vmatrix}
        \ptc{S}_0 & \myln{p}_j & \myln{p}_m\\
        \myln{p}_j& d_j & a_{jm}\\
        \myln{p}_{\ell} & a_{j\ell} & a_{m\ell}
    \end{vmatrix},
\end{equation}
obtaining
\begin{align}
\label{faceconicsanddiag2}
    a_{jm}\ptc{H}_{\{j\},m,\ell}&=\begin{vmatrix}
        \begin{vmatrix}
        \ptc{S}_0 & \myln{p}_m \\
        \myln{p}_j& a_{jm} 
    \end{vmatrix} & \begin{vmatrix}
        \myln{p}_j & \myln{p}_m\\
        d_j & a_{jm}\\
    \end{vmatrix}\\
    \rule{0pt}{.00001ex} & \\
    \begin{vmatrix}
        \myln{p}_j& a_{jm}\\
        \myln{p}_{\ell} & a_{m\ell}
    \end{vmatrix} & 
    \begin{vmatrix}
        d_j & a_{jm}\\
        a_{j\ell} & a_{m\ell}
    \end{vmatrix}
    \end{vmatrix}\nonumber \\
    &=\begin{vmatrix}
        \ptc{H}_{\emptyset,j,m} & \myln{p}_{\{j\},m}\\
        \myln{q}_{m} & g_{\{j\},\ell,m}
    \end{vmatrix} \nonumber \\
    &=   g_{\{j\},\ell,m} \ptc{H}_{\emptyset,j,m} - \myln{p}_{\{j\},m} \myln{q}_{m}
\end{align}
yielding
\begin{equation}
    g_{\{j\},\ell,m} \ptc{H}_{\emptyset,j,m} - a_{jm}\ptc{H}_{\{j\},m,\ell} =  \myln{p}_{\{j\},m} \myln{q}_{m}  
\end{equation}
As with Eq. \ref{faceconicsanddiag1}, this shows that $\ptc{H}_{\{j\},m,\ell}$ has two intersection points with $\ptc{H}_{\emptyset,j,m}$ on $\myln{q}_m$.

We now prove the statement at the beginning of this section for the case of the four faces $\mathcal{F}_1$, $\mathcal{F}_{\bar{1}}$, $\mathcal{F}_3$, $\mathcal{F}_{\bar{3}}$ which remain when the top and bottom faces, $\mathcal{F}_{\bar{2}}$, $\mathcal{F}_2$, are excluded.  Putting $j=3$ into equation \eqref{faceconicsanddiag1} yields the conclusion that the face conics for $\mathcal{F}_3$ and $\mathcal{F}_{\bar{3}}$ meet in two points on $\myln{q}_2$ (as well as $\myln{q}_1$).  Then, putting $j=3$, $m=2$ in equation \eqref{faceconicsanddiag2}, we further conclude that the face conics for $\mathcal{F}_3$ and $\mathcal{F}_{\bar{1}}$ meet in two puts on $\myln{q}_2$ as well.  Finally, plugging $j=1$ into equation \eqref{faceconicsanddiag1}, we find that the face conic for $\mathcal{F}_1$ and $\mathcal{F}_{\bar{1}}$ meet in $\myln{q}_2$.  Putting these incidences together, we arive at the conclusion stated at the beginning of this subsection.

\subsection{Subdeterminants of $\ptc{S}_{\{1,2,3\}}$ and the Penrose cube}
\label{sec:subdet}

Here we take the results of the previous section and apply them to the specifics of a Penrose 7-conic configuration.  Given $\ptc{S}_0$ and the three chords of contact $\myln{p}$, $\myln{q}$, and $\myln{r}$, Eqs. \ref{eq:Tijk-det}-\ref{eq:T0-det-2} show that the remaining 7 entries in the Penrose cube can be written as symmetric subdeterminants of $\ptc{S}_{\{1,2,3\}}$ involving $\ptc{S}_0$. As shown above, the non-symmetric subdeterminants express other aspects of the Penrose cube. Table \ref{tab:determinants} presents all possible cases, their geometric interpretation, and example equations from Sec. \ref{sec:alg-proof}. 

In the $2\times2$ subdeterminants: the indices $\{i,j,k\}$ that appear in a single expression in the following table take on distinct values in $\{1,2,3\}$.  The total number of $2 \times 2$ subdeterminants can be calculated by taking into account the underlying symmetry $\subd{I_r}{I_c} = \subd{I_c}{I_r}$, yielding a total of 21.   Similar considerations show that there are 10 $3 \times 3$ subdeterminants.

The table mixes the original notation introduced for the Penrose cube with the more uniform one used in Sec. \ref{sec:det-structure}. Recall that the two notations are related as follows:
\begin{align*}
    \{\myln{p},\myln{q}, \myln{r}, -1, -1, -1, a, b, c\} &\Leftrightarrow \{\myln{p}_1,\myln{p}_2,\myln{p}_2, d_1, d_2, d_3, a_{23}, a_{13}, a_{12}\} \\
 \{\ptc{S}_{\emptyset}, \ptc{S}_{\{i\}}, \ptc{S}_{\{i,j\}},\ptc{S}_{\{1,2,3\}}\} &\Leftrightarrow \{\ptc{S}_0, \ptc{T}^i, \ptc{S}_k, \ptc{T}^0\} 
\end{align*}

Most of the entries in the \emph{Description} column were introduced in Section \ref{sec:alg-proof}. There are two exceptions. For a discussion of the 8-point conic see Section \ref{sec:face-conics}; for a discussion of face diagonals, see Section \ref{sec:face-diag}.

\begin{table}[ht]
    \centering
    \begin{tabular}{|c|c|c|c|c|}
    \hline
    $\subd{-}{-}$ & \# & Formula & Description & Example \\ \hline \hline
    $\subd{0i}{0i}$& 3 & $d_i \ptc{S}_0 - \myln{p}_i^2$ & $\ptc{T}^i$ & Eq. \ref{eq:t1}\\ \hline
    $\subd{0i}{0j}$ & 3 & $a_{ij}\ptc{S}_0 - \myln{p}_i\myln{p}_j$ &  8-point conic of $\face_{0i}^{jk}$ & Eq. \ref{eq:face-conics} \\ \hline
    $\subd{0i}{ji}$ & 6 & $a_{ij}\myln{p}_i - d_i\myln{p}_j$& contact chords $\ptc{T}^i$-$\ptc{S}_k$ & Eq. \ref{eq:level1-doubleLine} \\ \hline
    $\subd{0i}{jk}$ & 3 & $a_{ik}\myln{p}_j - a_{ij}\myln{p}_k$ & face diagonal & Eq. \ref{eq:face-diag} \\ \hline
    $\subd{ij}{ij}$ & 3 &  $a_{ij}^2 - d_i d_j$ & weight of $\ptc{S}_0$ in  $\ptc{S}_k$ & Eq. \ref{eq:Sijk-det} \\ \hline
    $\subd{ij}{ik}$ & 3 & $d_ia_{jk} - a_{ij}a_{ik}$ & weight of $\myln{p}_j\myln{p}_k$ in $\ptc{T}^0$ & Eq. \ref{eq:T0-det} \\ \hline \hline
    $\subd{0ij}{0ij}$& 3 & $(d_id_j-a_{ij}^2)\,\ptc{S}_0-d_j\abdb{p_i} + 2\, a_{ij}\, \abpr{p_i}{p_j} - d_i\abdb{p_j}$ & $\ptc{S}_k$ & Eq. \ref{eq:Sijk-det}\\ \hline
    $\subd{0ij}{0ik}$& 3 & \text{\tiny $ a_{ik} p_i p_j + a_{ij} p_i p_k - d_i p_j p_k -a_{jk} p_i^2 + (a_{jk}d_i - a_{ij} a_{ik}) S_0 $} & \text{\normalsize{ 8-point conic of $\face^{0k}_{ij}$}} & Eq. \ref{eq:face-conics}  \\ \hline
    $\subd{0ij}{ijk}$& 3 & \text{\tiny $ (a_{ij} a_{jk} - a_{ik} d_j) p_i + (a_{ij} a_{ik} - a_{jk} d_i) p_j +(d_i d_j- a_{ij}^2) p_k $} &\text{\normalsize{ contact chords $\ptc{S}_k$-$\ptc{T}^0$}} &  \\ \hline
    $\subd{ijk}{ijk}$& 1 & $1 - a_{ij}^2 - a_{ik}^2 -a_{jk}^2- 2 a_{ij}a_{ik}a_{jk}$ & weight of $\ptc{S}_0$ in $\ptc{T}^0$ & Eq. \ref{eq:T0-det} \\ \hline \hline
    $\subd{0123}{0123}$& 1 & ... & $\ptc{T}^0$ & Eq. \ref{eq:T0-det}  \\ \hline
    \end{tabular}
    \caption{Subdeterminants of $\ptc{S}_{\{0,1,2,3\}}$ and their meaning.}
    \label{tab:determinants}
\end{table}

\section{Related research and a brief history}
\label{sec:rere_history}

We recently became aware of \href{https://abelie.no/articles/morfologi.pdf}{an independent discovery of the 8-conic theorem} \cite{eide2013}, Sect. 9.2. (\href{https://abelie.no/articles/morphology.pdf}{English translation here}.) The author, Morten Eide, has in the meantime become a co-author of this article.

The seed for this article was planted when one of the authors heard \href{https://www.youtube.com/watch?v=JiDWGbsVEno\&t=1317s}{a podcast featuring Roger Penrose} \cite{numPenrose2020}, in which he describes formulating the 8-conic theorem in the early 1950's, but not having found the time to publish it. 

In May, 2024 some of the current authors 
contacted Roger Penrose to report that they had obtained a proof of the theorem. A collaboration ensued in which Roger Penrose shared details of his original approaches to the proof of the theorem.
He recalled that as an undergraduate he first recognized the cube structure of the theorem from examples in which most of the conic sections are degenerate. He then developed a method for gluing together several such cubes to generate new cubes that contain more regular conics. In this way, he finally arrived at the formulation of the general theorem.  Later, as a graduate student, he discovered an approach to prove the 8-conic theorem using an analogous 8-quadric theorem in 3-D.
It is this approach that is reconstructed in the geometric proof in Sec. \ref{sec:the_main_theorems}. The algebraic proof presented in Sec. \ref{sec:alg-proof} is due to Morton Eide.

Our initial proof of the 8-conic theorem, in contrast, worked in the 5-D parameter space of conics. After the proof of the 8-quadric theorem in 3-D (presented in this article) was obtained, we also developed a proof in the 9-D parameter space of quadrics. 
The parameter-space proofs have a similar structure to each other, but just as in the proofs presented in the first part of this article, the parameter-space proof of the 8-quadric theorem is shorter and easier than the parameter-space proof of the 8-conic theorem. We hope to present these proofs in the near future in one or two further preprints.

\charliex{
\section{Recapitulation and Outlook}
\label{sec:outlook}

Beginning with a series of examples of familiar theorems in projective geometry we formulated the Penrose 8-conic theorem, that includes these examples as special cases, and proved it, first geometrically and then algebraically. Here we attempt to evaluate the significance of the theorem and indicate some directions for future research.  Much but not all of the following discussion can be transferred also to the associated Penrose 8-quadric theorem.

We saw that the theorem has a striking universality, so that many theorems of projective geometry can be positioned as special cases.  Furthermore, the determinant formulas in Section \ref{sec:det-structure} provide an explicit parametrization of all valid Penrose cubes, yielding a (17-dimensional) space whose points each represent a Penrose cube. Since such a cube can be interpreted as a theorem, any path between two points in the space represents a continuous metamorphosis of one theorem to the the other, in which all intermediate stages are also theorems.  This is a certainly a remarkable space. 

\ifthenelse{\boolean{isCoda}}
{
After we finish our recapitulation, we return to this theme in the Coda (Section \ref{sec:coda}).
}{}


One general question is: Which projective theorems correspond to points in this \emph{Penrose space} and which not? That is, how much of projective geometry can be seen as a consequence of the Penrose 8-conic theorem?  How many different ``theorems'' does it contain? To answer this, observe that it is divided into regions by sub-manifolds where the rank of a matrix of one of the 8 conics changes. Each region can be considered to represent one equivalence class of theorems. Factoring out by different labeling of the same cube can be used to further reduce the number of equivalence classes. The task remains daunting.

Our approach has been more ``bottom-up'' by beginning with specific theorems and searching for patterns there. For example, the Penrose cubes for many of known theorems have a 3-fold symmetry, whereby the 4 levels of a Penrose cube (balanced on its bottom vertex $\ptc{T}^0$) each have the same rank and type (point-wise or line-wise). Consult Figure \ref{fig:cubegraphs}.  We have carried out a survey of 13 or 14 possible types, producing the beginning of a classification that will be published separately.

Besides the rank of the conics, other factors that produce variation in the specialization include:
\begin{description}
    \item[Distinction of point-wise and line-wise conics] The conics of the Penrose cube are complete conics (Section \ref{sec:nota-for-comp-conics}), including both a point-wise and line-wise partner.  Note that dualizing a Penrose theorem is obtained by dualizing its Penrose cube.  The formulas of Section \ref{sec:alg-proof} necessarily apply to either point-wise or line-wise conics.  The interaction between these two aspects can be subtle and deserves further attention, especially the rare cases that one or the other of the partners vanishes.
    \item[Real and complex elements] Our theory considers conics in $\RP{2}$, whose polynomials have real coefficients. But the zero sets of these polynomials can become imaginary in manifold ways. Tracking this transition is an important aspect that deserves further attention.  The presence of imaginary conics also allows connections to euclidean geometry, as the next item explains.
    \item[Interaction with metric geometry] The Cayley-Klein construction of metric spaces within projective space can find expression in a Penrose cube since the Cayley-Klein metric absolute is a conic.  For example, a euclidean circle is a conic that has double contact with the euclidean absolute imaginary points $(I,J)$. Including this point pair in a Penrose cube guarantees that its neighbors (when regular) will be circles. The Monge Circle theorem discussed above is one of many examples.  On the other hand, euclidean confocal conics are line-wise conics with a particular relationship to the euclidean absolute that allows Penrose cubes to be contructed expressing Euclidean theorems about confocal conics and focal points. Such examples raise the question to what extent the domain of the Penrose 8-conic theorem extends also into metric geometry. 
\end{description}

We have also explored extensions of the 8-conic theorem. The possibility of gluing several Penrose cubes together, mentioned above, remains an interesting option to obtain further theorems. For example, this method can be used to construct Pocelet configurations for general polygons out of those for triangles. Another example is Neville's theorem about three ellipses that have a common focal point in pairs: Here, a first application of the eight-conic theorem can be used to show that all three ellipses must have double contact with another conic. This then serves as a starting point for further cubes, attached to the first one. We are currently working on examples of this type of iterated application of the theorem.


We also saw that the determinant theory of Section \ref{sec:det-structure} can be directly extended to other configurations, such as hyper-cubes. 
There are also outstanding questions regarding the finer points of the subdeterminants of $\ptc{S}_{\{1,2,3\}}$; for example, to work out in detail how the simultaneous vanishing of several subdeterminants can effect in subtle ways the dual partners of particular conics in the cube. 
Another, more general direction is, in the underlying matrices, to replace the linear forms $\myln{p}_i$  with $n-$order curves; conics will be replaced with $2n$-order curves, whereby the number of contact points is $2n^2$. 

\ifthenelse{\boolean{isCoda}}
{}
{A more general direction of research involves the place of the Penrose 8-conic theorem within projective geometry. It is an archetypal theorem, containing within itself innumerable known and not-yet-known theorems. Based on our experience in the process of preparing this article, it has the potential to enrich how projective geometry is practiced and perceived, expanding the focus from isolated theorems to a continuum of metamorphoses of which the theorems themselves are static snapshots.  
}

\ifthenelse{\boolean{isCoda}}
{
\section{Coda} \label{sec:coda}

\emph{Note:} The views expressed here are not necessarily shared by all the authors of the article. They are included here to stimulate the reader's imagination and engagement with the contents of the article. 

The fact that the Penrose 8-conic theorem includes many classical theorems might well evoke a sense of \emph{deja vu} in students of projective geometry.  Without Desargues' introduction of ideal points -- where parallel lines meet -- any projective theorem breaks up into a collection of specialized theorems, depending on which sets of lines of the figure are parallel; each possibility requires an extra theorem to handle.  The projective version unifies all the special cases in one. The remarkable thing is, the Penrose 8-conic theorem then unifies these unities into a higher unity.

Something similar had emerged in the field of botany, coincident with the rebirth of projective geometry at the beginning of the 19th century.  Johann Wolfgang Goethe, in his studies of plant morphology (\cite{goethe1791}, \cite{Bortoft1996WholenessOfNature}), arrived at the notion of the \emph{archetypal plant} out of which all specific plant species are specializations. In contrast to his predecessor Linnaeus, who organized the plant kingdom in a tree structure based on external, measurable properties, Goethe took his point of departure from the archetype, the One, and described how each of the Many plant species is a specialization of it.  Each species, in turn, is itself a unity that manifests in the many individual earthly plants belonging to the species. The species are a dynamic unity, and can be brought into a metamorphosis transforming one species into another. Goethean knowledge is inseparable from a training of the perceptual and thinking faculties of the scientist in order to mentally grasp these aspects. 

We can say that an approach that emphasizes the Many at the expense of the One loses itself in diversity (Linnaeus), while an approach that one-sidedly emphasizes the One tends to Platonism (or the reductionism of molecular biology).  For Goethe, the healthy balance was found in understanding how these two, the One and the Many (or, the whole and the parts), condition each other. Goethe was convinced that only in the context of such a polarity can one speak of a  \emph{living organism} rather than a \emph{machine}.  The active participation of the scientist in the act of knowing is characteristic of this approach. A detailed discussion of this theme lies outside the scope of this article. 

The formal similarities to our findings on the Penrose 8-conic theorem lie close at hand.  The 8-conic theorem serves as an archetype for the vast profusion of specialized theorems described above. We also noted that each of these projective theorems is a species in itself, split into many by the fall from the projective into the euclidean world.  Moving within the parameter space calls forth a continual metamorphosis of one projective theorem into another. 


These similarities alone suggest the aptness of the term \emph{organic geometry} for this field of research. We are not the first to suggest it. Already in 1870, William Clifford, a shooting star in the mathematical heavens, strongly recommended the use of this term instead of the then-popular \emph{synthetic geometry} (\cite{Clifford1870MiquelTheorem}), citing as precedent Heinrich Gretschel \cite{Gretschel1868OrganischeGeometrie}, who in turn credited the pioneer of synthetic geometry, Jakob Steiner \cite{Steiner1832SystematischeEntwicklung}, as his inspiration. 

Is this startling resemblance just a coincidence? Or are there deeper forces at work?

To answer this, we first characterize the place of Goethe's scientific method in the history of science. Roughly speaking, it is an extension of the scientific method as practiced since the Renaissance by bringing together the ``object'' of study with the ``subject'' of study, that is, the human researcher who, through thinking, creates order in the natural or mathematical phenomena. In this paradigm, there is no ``object'' of study until the ``subject'' creates it, and more importantly, the resulting knowledge of the object is not the same when separated from its subject. We could characterize this as a transition from ``third person'' to ``first person'' knowledge.

We now call attention to three aspects of our current situation:

\begin{itemize}
    \item Artificial intelligence can be seen as the crowning achievement of third-person knowledge. For observers whose horizon is limited to third-person knowledge, the anxious question arises, ``What is left for the human being to do?''.The difficulty of finding an answer has led to transhumanism, that envisions merging the ``out-smarted'' humanity with the machines.
    \item The emergence of AI has been remarkably aligned with the severe phase of the global climate crisis. While the diagnosis of this crisis is clear, the therapy is not.  To the extent that human thinking, acting through science and technology, is itself responsible the crisis, is it reasonable to expect that healing can be obtained by applying more of the same ``third person'' thinking?
    \item Projective geometry has established itself within mathematics as the main branch (or archetype) of geometry, while its influence on the natural sciences is virtually nil. In contrast, many physicists of the 20th century have acknowledged the key role played by mathematics in their discoveries. For example, Eugene Wigner (1960) wrote: "The miracle of the appropriateness of the language of mathematics for the formulation of the laws of physics is a wonderful gift which we neither understand nor deserve."
\end{itemize}

We close our considerations with an attempt to provide a narrative bringing these three disparate threads together:

    \emph{The evolution of human consciousness continues. The era of ``third person'' science and society has fulfilled its function and is giving way, in fits and starts, to that of the ``first person''.  In dying, however, it has given us in artificial intelligence a powerful tool that externalizes third person knowledge. This frees us to turn with an open mind to the pressing crises surrounding us, crises that have their origin in our dealing with the organic world, both in Nature and society, with thinking fitted to the inorganic world.  In this transition from physics to biology, from social hierarchy to social organism, organic geometry can perhaps provide, using Wegner's words, the appropriate mathematical language for formulating the laws of \emph{organisms} within a participatory, ``first person'', Goethean-inspired science. }
}
{} 



}
\appendix
\section{Appendix}
\label{app:binique}

In this appendix, we show the existence of two quadrics completing the cube in the case where all six ring planes concur in a common axis.  By the results of sections \ref{sec:alg-proof} and \ref{sec:det-structure}, only one of the two will be attainable as the limit of cubes arising in the general case, where the ring planes don't have a common axis (see remark below).

In the case under consideration, when all ring planes meet in a common axis $\Omega$, all of the quadrics will be in double contact, i.e., share two contact elements: $\ptq{T}^1$ and $\ptq{T}^2$ share the points 
\begin{equation}
    \ptc{R}_0^1 \cap \ptc{R}_0^2=\ptc{R}_0^1 \cap \Omega
\end{equation}
and, by marching around the cube in the same fashion as Lemma \ref{lem:faconcur}, these two points are common to all the quadrics of the cube:
\begin{equation}
    \ptc{R}_i^j \cap \Omega=
    \ptc{R}_k^m \cap \Omega,
\end{equation}
for all indices $i,j,k,m$, as are the tangent planes at these points (the case where $\Omega$ is tangent to the rings provides some interesting variations but we omit consideration of it as it doesn't differ in the essentials).

In an appropriate basis, the matrices of the family of quadrics sharing these two contact elements can be written
\begin{equation}\label{eq:dcschoords}
    \ptc{Q}=a_0 \ptq{Q}_0+a_1 \ptq{Q}_1+a_2 \ptq{Q}_2+a_3 \ptq{Q}_3
\end{equation}
where
\begin{equation*}
    \ptq{Q}_0=\begin{pmatrix}
        0 & 1 & 0 & 0\\
        1 & 0 & 0 & 0\\
        0 & 0 & 0 & 0\\
        0 & 0 & 0 & 0\\
    \end{pmatrix}, \qquad
    \ptq{Q}_1=\begin{pmatrix}
        0 & 0 & 0 & 0\\
        0 & 0 & 0 & 0\\
        0 & 0 & 0 & 1\\
        0 & 0 & 1 & 0\\
    \end{pmatrix},
\end{equation*}
\begin{equation*}
    \ptq{Q}_2=\begin{pmatrix}
        0 & 0 & 0 & 0\\
        0 & 0 & 0 & 0\\
        0 & 0 & 1 & 0\\
        0 & 0 & 0 & 0\\
    \end{pmatrix}, \qquad
    \ptq{Q}_3=\begin{pmatrix}
        0 & 0 & 0 & 0\\
        0 & 0 & 0 & 0\\
        0 & 0 & 0 & 0\\
        0 & 0 & 0 & 1\\
    \end{pmatrix}.
\end{equation*}
In order for the contact elements to be contained in the literal sense, we must have $a_0=0$, thus this parametric space can be viewed as an affine $3$-space with the quadrics $\ptq{Q}$ given by equation \eqref{eq:dcschoords} as `points'.  We shall use single quotes when talking about geometric objects in this space in order to distinguish them from geometric elements in the space in which the $\ptq{Q}$ appear as surfaces, e.g., we shall refer to a pencil of quadrics lying in this parametric space as a `line'.  It makes sense to include also the `plane' $a_0=0$, making the parametric space into a projective space.  This `plane' consists of the plane pairs which meet in $\Omega$.  The double planes form a `conic' $\mathtt{V}$, the locus of `points' $\ptc{Q}$ satisfying
\begin{equation}
    a_0=0, \qquad a_1 a_2-a_3^2=0.
\end{equation}

\begin{figure}[ht!]
    \centering
    \includegraphics[width=0.6\textwidth]{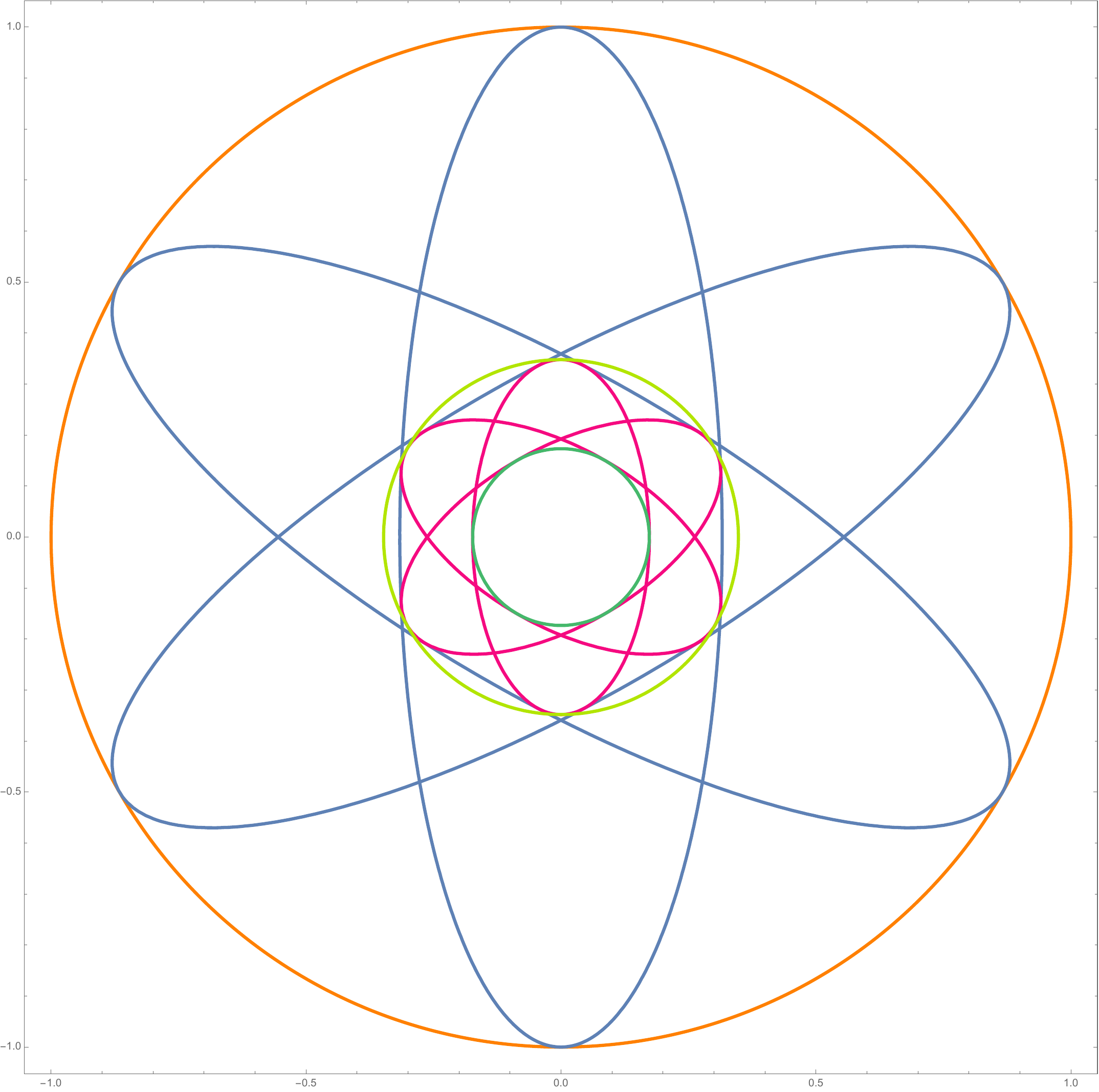}
    \caption{Penrose configuration where all face points coincide and, as a result, the cube has two completions.  $\ptc{S}_0$ is orange, $\ptc{T}^1$, $\ptc{T}^2$, $\ptc{T}^3$ are blue, $\ptc{S}_1$, $\ptc{S}_2$, $\ptc{S}_3$ are magenta, and the two possible completions $\ptc{T}^0$ are dark green and lime green.}
    \label{fig:binique}
\end{figure}

In summary, when all the ring planes share an axis, all the quadrics of the cube can be viewed as belonging to a parametric $3$-space, including any possible choice of $\ptq{T}^0$ that completes the cube.  Any quadric in this space and in ring contact with a given quadric $\ptq{Q}$ lies on a `line' containing a `point' of $\mathtt{V}$.  Therefore, the family of quadrics in ring contact with $\ptq{Q}$ is the `cone' $\mathtt{C}_{\ptq{Q}}$ made up of the `lines' connecting $\ptc{Q}$ to any `point' of $\mathtt{V}$.  Any possible $\ptq{T}^0$ must therefore be in the intersection
\begin{equation}
    \mathtt{C}_{\ptq{S}_1} \cap \mathtt{C}_{\ptq{S}_2} \cap \mathtt{C}_{\ptq{S}_3}.
\end{equation}
The intersection of two `cones' in general forms a fourth degree curves, however, in the present case, they meet in the `conic' $\mathtt{V}$ and thus there intersection breaks down into a pair of `conics':
\begin{equation}
    \mathtt{C}_{\ptq{S}_i} \cap \mathtt{C}_{\ptq{S}_j}=\mathtt{V} \cup \mathtt{W}_{ij}.
\end{equation}
Therefore
\begin{equation}
    \mathtt{C}_{\ptq{S}_1} \cap \mathtt{C}_{\ptq{S}_2} \cap \mathtt{C}_{\ptq{S}_3}=\mathtt{V} \cup (\mathtt{W}_{12} \cap \mathtt{W}_{13}).
\end{equation}
In general, $\mathtt{W}_{12} \cap \mathtt{W}_{13}$ consists of two `points' either of which can be taken as a regular quadric completing the cube.  Figure \ref{fig:binique} illustrates this for the case of conics, where the condition of the planes sharing a common axis is replaced by that of the face points of all faces coinciding (see \cite{penthm24-2d} for more on these parametric three spaces in the conic case).

\begin{remark}
    The approach of sections \ref{sec:alg-proof} and \ref{sec:det-structure} always provides a unique completion regardless of whether or not the face axes/face points collapse to one.  Therefore the two completions provided in the present appendix can be distinguished: one that arises from the determinant formulas and is thus stable under parametric perturbations which disturb the coincident face-axes/face-point property, and one which is only present thanks to this coincidence.  In figure \ref{fig:binique}
    the former is the dark green central circle and the latter is the lime green circle surrounding the magenta conics.
\end{remark}

\bibliography{main}
\bibliographystyle{alpha}

\end{document}